\pdfoutput=1
\RequirePackage{ifpdf}
\ifpdf 
\documentclass[pdftex]{sigma}
\else
\documentclass{sigma}
\fi

\usepackage{tikz}
\numberwithin{equation}{section}

\usepackage{mathtools}
\usepackage[all]{xy}

\newtheorem{Theorem}{Theorem}[section]
\newtheorem{Corollary}[Theorem]{Corollary}
\newtheorem{Lemma}[Theorem]{Lemma}
\newtheorem{Proposition}[Theorem]{Proposition}
 { \theoremstyle{definition}
\newtheorem{Definition}[Theorem]{Definition}

\newtheorem{Example}[Theorem]{Example}
\newtheorem{Remark}[Theorem]{Remark} }

\def\mbb #1{\mathbb{#1}}
\def\mbf #1{\mathbf{#1}}

\def\Cal #1{\mathcal{#1}}

\def\C{\mbb{C}}
\def\CP{\mbb{CP}}
\def\R{\mbb{R}}
\def\Z{\mbb{Z}}
\def\N{\mbb{N}}
\def\GL{{\rm GL}}
\def\Mat{\operatorname{Mat}}
\def\tr{\operatorname{tr}}
\def\Re{\operatorname{Re}}
\def\Im{\operatorname{Im}}

\def\res{{\rm Res}}
\def\cM{\mathcal{M}}
\def\cE{\mathcal{E}}
\def\sX{\mathsf{X}}

\def\sM{\mathsf{M}}
\def\sR{\mathsf{R}}
\def\sS{\mathsf{S}}
\def\rS{\mathsf{\check S}}
\def\rX{\mathsf{\check X}}
\def\rM{\mathsf{\check M}}
\def\rXO{{\mathsf{\check X}_O}}
\def\rXI{{\mathsf{\check X}_I}}

\def\widebar #1{\mkern 2mu\overline #1}
\def\Barbar #1{\Bar{\Bar #1}}
\def\p{\mathsf{\mkern1.5mu P} }

\newcommand*\hg{\begingroup
 \dopFq}
\def\dopFq#1#2#3{\left(\genfrac{}{}{0pt}{1}{#1}{#2} ; #3\right)%
 \endgroup}

\begin{document}
\allowdisplaybreaks

\newcommand{\arXivNumber}{1301.5228}

\renewcommand{\PaperNumber}{006}

\FirstPageHeading
	
\ShortArticleName{Unfolding of a Resonant Irregular Singularity}
	
\ArticleName{Analytic Classification of Families\\ of Linear Differential Systems Unfolding\\ a Resonant Irregular Singularity}
	
\Author{Martin KLIME\v{S}}
	
\AuthorNameForHeading{M.~Klime\v{s}}
\Address{Independent Researcher, Prague, Czech Republic}
\Email{\href{mailto:klmm@seznam.cz}{klmm@seznam.cz}}
	
\ArticleDates{Received March 13, 2019, in final form December 21, 2019; Published online January 23, 2020}
	
\Abstract{We give a complete classification of analytic equivalence of germs of parametric families of systems of complex linear differential equations unfolding a generic resonant singularity of Poincar\'e rank 1 in dimension $n = 2$ whose leading matrix is a Jordan bloc. The moduli space of analytic equivalence classes is described in terms of a tuple of formal invariants and a single analytic invariant obtained from the trace of monodromy, and analytic normal forms are given. We also explain the underlying phenomena of confluence of two simple singularities and of a turning point, the associated Stokes geometry, and the change of order of Borel summability of formal solutions in dependence on a complex parameter.}
	
	
\Classification{34M03; 34M35; 34M40}

\section{Introduction}

A \emph{system of meromorphic linear differential equations} with a singularity at the origin can be written locally as $ \Delta_0(x) y=0$, $y(x)\in\C^n$, with
\begin{gather}\label{eq:RC-Delta0}
\Delta_0(x)=x^{k+1}\frac{{\rm d}}{{\rm d}x}-A_0(x), \qquad x\in(\C,0),
\end{gather}
where $A_0(x)$ is a $(n \times n)$-matrix with analytic entries in a neighborhood of~0, $A_0(0)\neq 0$, $(\C,0)$ stands for a germ of a neighborhood of the origin in~$\C$,
and~$k$ is a non-negative integer, called the \emph{Poincar\'e rank}. An \emph{unfolding} of~\eqref{eq:RC-Delta0} is a germ at~0 of a parametric family of systems $\Delta(x,m) y=0$, $y(x,m)\in\C^n$, with
\begin{gather}\label{eq:RC-1}
\Delta(x,m)=h(x,m)\frac{{\rm d}}{{\rm d}x}-A(x,m), \qquad (x,m)\in\big(\C\times\C^l,0\big),
\end{gather}
such that $\Delta(x,0)=\Delta_0(x)$, where the scalar function $h(x,m)$ and the $(n \times n)$-matrix function $A(x,m)$
depend analytically on both the variable $x$ and the parameter~$m$.
Two families of linear systems~\eqref{eq:RC-1} depending on the same parameter $m$ are \emph{analytically equivalent}
if there exists an invertible analytic linear gauge transformation bringing solutions of the first system to solutions of the second one.

The analytic classification of singularities of single systems \eqref{eq:RC-Delta0} is now well known in full generality (see, e.g., \cite{BaVa,Balser, BJL12-I,BJL12-II,IlYa}).
While for singularities of Poincar\'e rank $k=0$, called \emph{Fuchsian singularities}, the analytic classification coincide with the formal one: a formal power series transformation between two Fuchsian singularities always converges,
this is no longer true in the case of \emph{irregular singularities}\footnote{A singularity of a system \eqref{eq:RC-Delta0} is \emph{Fuchsian} if it has Poincar\'e rank 0, and
	it is \emph{regular} if it is meromorphically gauge equivalent to a Fuchsian singularity, otherwise it is \emph{irregular}.
	A Fuchsian singularity is \emph{non-resonant} if no two eigenvalues of the leading matrix $A_0(0)$ differ by a (non-zero) integer, while 	
	an irregular singularity is \emph{non-resonant} if the eigenvalues of the leading matrix $A_0(0)$ are all distinct.} of higher Poincar\'e rank $k>0$, for which formal gauge transformations are generally divergent.

The geometric description of this divergence is known as the \emph{Stokes phenomenon}.
The essence of this phenomenon is the following: formal gauge transformations to a normal form are always asymptotic to some true analytic ones,
which also conjugate the systems, but only over certain sectors in the $x$-space, the union of which covers a full neighborhood of the singularity. If continued to larger sectors, these transformations would in general explode near the singularity and loose their asymptoticity. The mismatch between the different sectoral transformations on the overlapping of the sectors constitutes an obstruction to the convergence of the formal transformation.
This obstruction may be expressed in terms of the transition automorphisms between the different sectoral gauge transformations, and represented by, so called, Stokes matrices.
Roughly speaking, the set of these Stokes matrices forms together with a set of formal invariants a complete analytic invariant of the irregular singularity.

Since the discovery of the Stokes phenomenon, it has always been tempting to try to understand it from the point of view of confluence of singularities in perturbed families \eqref{eq:RC-1}.
The investigation of parametric unfoldings of singularities has several goals:
\begin{enumerate}\itemsep=0pt
\item[1)] to provide the modulus of analytic equivalence and/or analytic normal forms for germs of parametric systems,
\item[2)] to explain the Stokes phenomenon of irregular singularities through confluence of several Fuchsian ones,
\item[3)] to understand the degeneration of linear problems (such as isomonodromic systems associated to Painlev\'e equations), and the limits in moduli spaces of meromorphic connections over Riemann surfaces.
\end{enumerate}

It has been conjectured independently by V.I.~Arnold, A.A.~Bolibruch and J.-P.~Ramis,
that Stokes matrices of the limit problem can be obtained as limits of transition matrices between the canonical solution bases at each of the regular singular points of a generically perturbed system with respect to which the corresponding local monodromy is diagonal.
This was later proved by A.~Glutsyuk for non-resonant~\cite{G0}
and certain resonant irregular singularities~\cite{Glu}.
But such approach covers only the sectors in the parameter space, on which the deformation is generic: where all the singularities (roots of $h(x,m)$) and non-resonant Fuchsian.
Being able to cover also the parameter values for which some of the singularities may be resonant Fuchsian or irregular is essential for analytic classification.
This problem was recently resolved by J.~Hurtubise, C.~Lambert and C.~Rousseau \cite{HLR, HR, LR}, and independently by L.~Parise~\cite{Parise},
for germs of parametric families of systems unfolding a non-resonant irregular singularity, that is one for which the leading matrix coefficient $A_0(0)$ in~\eqref{eq:RC-Delta0} has all eigenvalues distinct.
Their approach is based on forming, so called, \emph{mixed solution bases} associated to pairs of different singularities in the unfolding \eqref{eq:RC-1},
a method previously known to be fruitful when applied to hypergeometric systems and their generalizations \cite{Duv,LR1,Sch,Zha}.
The natural domains on which such mixed bases are defined, while rather complicated, are closely related to the real dynamics of the complex vector field
${\rm e}^{{\rm i}\theta}h(x,m)\frac{\partial}{\partial x}$ for some $\theta\in\R$, or equivalently, to the \emph{horizontal foliation} of the meromorphic Abelian differential
${\rm e}^{-{\rm i}\theta}\frac{{\rm d}x}{h(x,m)}$. As in the nonparametric case, on the overlapping of these domains the different mixed bases can be related by \emph{unfolded Stokes matrices}, the set of which constitute the analytic part of the invariant~\cite{HLR, LR}. Furthermore, as a byproduct, such description provides a canonical decomposition of the monodromy matrices into products of these unfolded Stokes matrices (which are unipotent) and diagonal matrices determined by the formal invariant.

This article provides the first results on analytic classification of parametric families unfolding a \emph{resonant irregular singularity}.
We consider germs of parametric families of $(2 \times 2)$-systems $\Delta(x,m)$ in a neighborhood of $(x,m)=0$, that unfold a system $\Delta_0(x)=\Delta(x,0)$ of the form
\begin{gather}\label{eq:RC-2}
\Delta_0(x)=x^2\frac{{\rm d}}{{\rm d}x}-A_0(x),\qquad\text{with}\quad A_0(0)=\begin{pmatrix} \lambda_0^{(0)} & 1 \\ 0 & \lambda_0^{(0)} \end{pmatrix},
\end{gather}
which has a resonant singularity of Poincar\'e rank 1 at the origin, under a \emph{generic condition} that the element $a_{21}^{(1)}$ on the position $(2,1)$ of the
matrix $\mathop{\res}\limits_{x=0}\frac{A_0(x)}{x^2}=\tfrac{{\rm d}}{{\rm d}x}A_0(x)\big|_{x=0}$ is non-zero:
\begin{gather}\label{eq:RC-2a}
a_{21}^{(1)}=-\frac{{\rm d}}{{\rm d}x}\det\big(A_0(x)-\lambda_0^{(0)}I\big)\big|_{x=0}\neq 0.
\end{gather}
\emph{No restriction will be imposed on the nature of the analytic deformation $\Delta(x,m)$ of $\Delta_0(x)$ or on the complex parameter $m\in\big(\C^l,0\big)$.}

We will provide a complete analytic classification of all germs of parametric systems $\Delta(x,m)$ unfolding such a $\Delta_0(x)$ (Theorem~\ref{theorem:RC-I})
in terms of a set of formal invariants and a single analytic invariant obtained from the trace of monodromy.
We also provide an explicit analytic normal form, i.e., a universal unfolding for any system $\Delta_0$ \eqref{eq:RC-2} satisfying~\eqref{eq:RC-2a} (Theorem~\ref{theorem:RC-II}), which is closely related to the biconfluent hypergeometric equation and to the modified Bessel's equation (Example~\ref{example:RC-Bessel}).

More importantly, we will explain both the phenomena of confluence of singularities, and of confluence of the eigenvalues of the principal matrix $A(0,m)$ of the system, resulting in change of order of Borel summability of formal solutions, and corresponding to the apparition of a~resonance in the irregular singularity.
To the best of our knowledge, this phenomenon has not been studied before. The geometric explanation that we shall offer is that of a~confluence of a~singularity and of a ``turning point'' (Section~\ref{sec:RC-1.3}).

Our approach follows the same footsteps as that of Hurtubise, Lambert and Rousseau \cite{HLR,HR, LR}, by constructing a set of canonical mixed solution bases on certain domains in the $(x,m)$-space (Theorem~\ref{theorem:RC-IV}), which, in turn, is equivalent to constructing ``sectoral'' gauge transformations on these domains between formally equivalent families (Corollary~\ref{theorem:RC-III}).
While the form of these domains is again rather complicated, they are determined by the demand on the sectoral gauge transformations to be bounded, a condition that is related to the asymptotic behavior of the solutions near the singularities.
In the end, the form of these domains is determined by the geometry of the horizontal foliation associated to the polar part (principal part) of the meromorphic quadratic differential (see Section~\ref{sec:RC-1.3})
\[\frac{\det\big(A(x,m)-\frac{1}{2}\tr A(x,m)\big)}{h(x,m)^2}({\rm d}x)^2.\]

The unfolded Stokes matrices defined on the overlappings of these domains can be in our situation expressed almost explicitly (Theorem~\ref{theorem:RC-IV}) as functions of formal invariants and of the trace of monodromy.
In fact, it turns out that aside of the formal invariants there is a single additional analytic invariant given by the trace of the monodromy around all the singularities (this would no longer be true for more complicated parametric systems~\eqref{eq:RC-1} of Poincar\'e rank $k>1$ or dimension $n>2$).

Let us remark that while the class of systems investigated here is relatively special, the geometric intuition behind our treatment
provides a first glimpse of a general description of unfoldings of $2 \times 2$ meromorphic linear differential systems or connections,
that should be developed in a future study.

Finally, it is worth mentioning, that confluences of the kind studied here appear naturally not only in the aforementioned biconfluent hypergeometric equation, but also in other important settings such as in the isomonodromic deformation problems associated to certain Painlev\'e equations. Namely, a system~\eqref{eq:RC-1} with $h(x,m)=x^2$
unfolding \eqref{eq:RC-2}, \eqref{eq:RC-2a}, appears in the degeneration of the traceless $2 \times 2$ isomonodromic problem associated to the Fifth Painlev\'e equation $P_{\rm V}\to P_{\rm V}^{\deg}$, and also in the degeneration Sixth Painlev\'e equation $P_{\rm VI}\to P_{\rm V}$, where it appears as a formal $(2 \times 2)$-bloc of a $3\times3$ isomonodromic problem in Birkhoff normal form \cite[Appendix]{Kl}.

\section{Statement of results}

\begin{Definition}
Throughout the text $\Delta(x,m)$ will denote a germ at $0$ of a parametric family of systems \eqref{eq:RC-1} unfolding \eqref{eq:RC-2} satisfying \eqref{eq:RC-2a}.
For brevity, we call it a \emph{parametric system}. We denote
\[\Delta_m:=\Delta(\cdot,m)\]
the restriction of $\Delta$ to a fixed parameter $m$.
\end{Definition}

\begin{Definition}[gauge transformations]
Let $y=T(x,m) y'$ be a linear transformation of the dependent variable. Let us define a transformed system
\begin{gather}\label{eq:RC-gauge}
 T^* \Delta:=h\frac{{\rm d}}{{\rm d}x}-\left[T^{-1}AT-h T^{-1}\frac{{\rm d} T}{{\rm d}x}\right],
\end{gather}
which satisfies $ (T^* \Delta) y'=0$ if and only if $\Delta y=0$.

Two parametric systems $\Delta(x,m)=h(x,m)\frac{{\rm d}}{{\rm d}x}-A(x,m)$ and $\Delta'(x,m)=h'(x,m)\frac{{\rm d}}{{\rm d}x}-A'(x,m)$,
depending on the same parameter $m$, are \emph{analytically equivalent},
if there exists a~germ of invertible linear gauge transformation $T(x,m)$, depending analytically on $(x,m)$, such that $ h'^{-1}\cdot\Delta'=h^{-1}\cdot T^* \Delta$.
\end{Definition}

\begin{Definition}[the invariants]\label{definition:RC-3} \quad
\begin{enumerate}\itemsep=0pt
\item[(i)] Applying the Weierstrass preparation theorem on the function $h(x,m)$ in \eqref{eq:RC-1}, we can assume that,
up to multiplying $\Delta(x,m)$ by a non-vanishing germ of scalar function,
\begin{gather*}
 h(x,m)=x^2+h^{(1)}(m)x+h^{(0)}(m),\qquad h^{(1)}(0)=h^{(0)}(0)=0,
\end{gather*}
where the coefficients $h^{(0)}(m)$, $h^{(1)}(m)$ are analytic. \emph{We shall suppose that $h$ is in this form from now on.}
We then define invariant polynomials $\lambda(x,m)$, $\alpha(x,m)$ by
\begin{gather}
\lambda(x,m)=\tfrac{1}{2} \tr A(x,m) \mod h(x,m) \nonumber\\
\hphantom{\lambda(x,m)}{} = \lambda^{(1)}(m) x+\lambda^{(0)}(m), \nonumber\\
\alpha(x,m)=-\det\big[A(x,m) - \lambda(x,m)I\big] \mod h(x,m) \nonumber\\
\hphantom{\lambda(x,m)}{} =\alpha^{(1)}(m) x+\alpha^{(0)}(m),\qquad \alpha^{(0)}(0)=0,\label{eq:RC-lambdaalpha}
\end{gather}
using the Weierstrass division theorem. The generic condition \eqref{eq:RC-2a} means that $\alpha^{(1)}(0)\neq 0$.
We call the triple $h(x,m)$, $\lambda(x,m)$, $\alpha(x,m)$ \emph{formal invariants} of $\Delta$ (this terminology will be justified by Proposition~\ref{proposition:RC-formal}).

\item[(ii)] We define an \emph{analytic invariant} $\gamma(m)$ by
\begin{gather}\label{eq:RC-gamma}
\gamma(m)={\rm e}^{-2\pi{\rm i} \lambda^{(1)}(m)}\cdot\tr M(m),
\end{gather}
where $M(m)$ is a monodromy matrix of some fundamental solution matrix $Y(x,m)$ of the system $\Delta(\cdot,m)$
around the two zeros of $h(x,m)$ in the positive direction:
\[Y({\rm e}^{2\pi{\rm i}}x_0,m)=Y(x_0,m)M(m).\]
The value of $\gamma(m)$ is independent of the choice of the fundamental solution $Y(x,m)$ or of the point $x_0$, and can be calculated point-wise for each value of $m$.
\end{enumerate}
\end{Definition}

\begin{Proposition}[prenormal form]\label{proposition:RC-prenormal} \quad
\begin{enumerate}
\item[$(i)$] The invariants $h(x,m)$, $\lambda(x,m)$, $\alpha(x,m)$ and $\gamma(m)$ are analytic in $m$, and are invariant under analytic equivalence of systems.

\item[$(ii)$] Up to an analytic gauge transformation, the system $\Delta(x,m)$ can be written as
\begin{gather}\label{eq:RC-prenormal1}
\Delta(x,m)=h(x,m)\frac{{\rm d}}{{\rm d} x}-\begin{pmatrix} \lambda(x,m) & 1 \\ \alpha(x,m)+h(x,m)r(x,m) & \lambda(x,m) \end{pmatrix}
\end{gather}
for some analytic germ $r(x,m)$.
\end{enumerate}
\end{Proposition}

\begin{proof}
(i) Elementary from the formula \eqref{eq:RC-gauge}, the Weierstrass division theorem, and analytic dependence of solutions on parameters.
(ii) Since gauge transformations commute with scalar matrices, the same gauge transformation works both for $\Delta$ and $\Delta+\lambda I$, and one can suppose that the trace invariant $\lambda(x,m)=0$.

Let $A=(a_{ij})$ be the matrix of the system $\Delta$,
$A(0,0)=\left(\begin{smallmatrix} 0 & 1 \\ 0 & 0 \end{smallmatrix}\right)$.
The gauge transformation $T_1\colon y\mapsto T_1(x,m)y$,
with $T_1=\left(\begin{smallmatrix} 1 & 0 \\ -\frac{a_{11}}{a_{12}} & \frac{1}{a_{12}} \end{smallmatrix}\right)$, which is analytic since $a_{12}(0,0)=1$,
brings $\Delta$ to
$ \Delta_1=:T_1^* \Delta=h(x,m)\frac{{\rm d}}{{\rm d}x}-\left(\begin{smallmatrix} 0 & 1 \\ b_{21} & b_{22} \end{smallmatrix}\right)$ for some $b_{ii}(x,m)$, $i=1,2$.
Now use $ T_2(x,m)={\rm e}^{\int\tfrac{b_{22}(x,m)}{2h(x,m)}{\rm d}x} \left(\begin{smallmatrix} 1 & 0 \\ \frac{1}{2}b_{22} & 1 \end{smallmatrix}\right)$
to get rid of the term $b_{22}$, the trace of the matrix of $\Delta_1$, which is divisible by $h(x,m)$ by the assumption that $\lambda(x,m)=0$.
Then $T_2^* \Delta_1$ is in the demanded form \eqref{eq:RC-prenormal1}.
\end{proof}

\begin{Example}[biconfluent hypergeometric equation]\label{example:RC-Bessel}
The hypergeometric equation is given by the second order linear differential operator
\begin{gather*}
D\hg{a,b}{1+c}{t}:=t(\delta_t+a)(\delta_t+b)-\delta_t(\delta_t+c),
\end{gather*}
where $\delta_t=t\frac{{\rm d}}{{\rm d}t}$ is the Euler operator.
The change of variable $t\mapsto\frac{t}{ab}$ gives
\begin{gather}\label{eq:RC-biconfluent}
\widetilde D\hg{a,b}{1+c}{t}:=\tfrac{t}{ab}(\delta_t+a)(\delta_t+b)-\delta_t(\delta_t+c),
\end{gather}
whose limit $a,b\to\infty$ is the \emph{biconfluent hypergeometric equation}
\begin{gather}\label{eq:RC-biconfluent0}
\widetilde D\hg{\infty,\infty}{1+c}{t}:=t-\delta_t(\delta_t+c).
\end{gather}
This confluence was studied (for $a=b+\frac{1}{2}$) by A.~Duval \cite{Duv}.
A particular solution to the equation \eqref{eq:RC-biconfluent0} is known to be given by the hypergeometric function
\[{}_0F_1\hg{-}{1+c}{t}:=\sum_{n=0}^{+\infty}\frac{t^n}{n!(1+c)_n},\]
$()_n$ denoting the Pochhammer symbol.
Moreover, if $\phi(t)$ is any solution to $\widetilde D\hg{\infty,\infty}{1+c}{t}\phi=0$, then the function
$u(s)=s^c\phi\big(\frac{s^2}{4}\big)$ satisfies the \emph{modified Bessel equation}
\begin{gather*}
\big[c^2+s^2-\delta_s^2\big]u=0.
\end{gather*}

Let $y_2$ be a solution to \eqref{eq:RC-biconfluent}, and let $y_1=-(\delta_t+c)y_2$,
then the vector function $y=\left(\begin{smallmatrix}y_1\\ y_2\end{smallmatrix}\right)$
satisfies the associated system
\begin{gather*}
\big(1-\tfrac{t}{ab}\big)\delta_t y=\begin{pmatrix} t\frac{a+b-c}{ab} & -t\frac{(a-c)(b-c)}{ab}\vspace{1mm}\\ \frac{t}{ab}-1 & c(\frac{t}{ab}-1) \end{pmatrix},
\end{gather*}
which in the local coordinate $x=t^{-1}$ is written as a parametric system
\begin{gather*}
x\big(x-\tfrac{1}{ab}\big)\frac{{\rm d}y}{{\rm d}x}=\begin{pmatrix} -\frac{a+b-c}{ab} & \frac{(a-c)(b-c)}{ab}\vspace{1mm}\\ x-\frac{1}{ab} & c(x-\frac{1}{ab}) \end{pmatrix},
\end{gather*}
which is of the considered form for any $c\in\C$, with a deformation parameter $m=\big(\frac{1}{a},\frac{1}{b}\big)\in\big(\C^2,0\big)$. The invariants of this system are
\begin{gather*}
h(x,m)=x^2-\tfrac{1}{ab}x,\\
\lambda(x,m)=\tfrac{c}{2}x-\tfrac12\left(\tfrac1{a}+\tfrac1{b}\right),\\
\alpha(x,m)=\left[1-c\left(\tfrac1{a}+\tfrac1{b}\right)+\tfrac{9}{2}c^2\tfrac1{ab}\right]x+\tfrac14\left(\tfrac1{a}-\tfrac1{b}\right)^2+2c^2\left(\tfrac1{ab}\right)^2,\\
\gamma(m)=2\cos(\pi c),
\end{gather*}
(see Lemma~\ref{lemma:RC-1} for the invariant $\gamma$).
\end{Example}

\subsection{Analytic theory}\label{sec:RC-1.1}

Analytic classification of germs of single systems $\Delta_0(x)$ is classical and was originally given in~\cite{JLP2}. Here we provide an analytic classification of parametric systems $\Delta(x,m)$ unfolding~$\Delta_0(x)$.

\begin{Theorem}[analytic classification]\label{theorem:RC-I} \quad
\begin{enumerate}\itemsep=0pt
\item[$(a)$] Two germs of parametric systems $\Delta(x,m)$, $\Delta'(x,m)$ are analytically gauge equivalent if and only if their invariants $h$, $\lambda$, $\alpha$, $\gamma$ are the same:
\begin{alignat*}{3}
& h(x,m)= h'(x,m),\qquad &&\lambda(x,m) =\lambda'(x,m),&\\
& \alpha(x,m)= \alpha'(x,m),\qquad &&\gamma(m)=\gamma'(m).&
 \end{alignat*}
\item[$(b)$] Any four germs of analytic functions $ h(x,m)$, $\lambda(x,m)$, $\alpha(x,m)$, $\gamma(m)$ with $h(x,0)=x^2$, $\alpha^{(0)}(0)=0$ and $\alpha^{(1)}(0)\neq 0$,
are realizable as invariants of some parametric system $\Delta(x,m)$.
\end{enumerate}
\end{Theorem}

\begin{Corollary} \label{corollary:RC-1}Two germs of parametric systems $\Delta(x,m)$, $\Delta'(x,m)$ are analytically equivalent if and only if there exists a product neighborhood $\sX\times\sM$ of 0 in $\C\times\C^l$ such that for each $m\in \sM$ the restricted systems $\Delta_m(x)$, $\Delta'_m(x)$ are analytically equivalent on $\sX$.
\end{Corollary}

The following Theorem~\ref{theorem:RC-II} provides an \emph{essentially unique normal form} for any germ of parametric system unfolding $\Delta_0$.

\begin{Theorem}[universal unfolding]\label{theorem:RC-II} Let $\Delta(x,m)$ be a germ of parametric system, and $h(x,m)$, $\lambda(x,m)$, $\alpha(x,m)$, $\gamma(m)$ its invariants.
\begin{enumerate}\itemsep=0pt
\item[$(i)$] If $\gamma(0)\neq 2$, then $\Delta(x,m)$ is analytically equivalent to a germ at~$0$ of a parametric system $\widetilde\Delta(x,m)$ given by
\begin{gather}\label{eq:RC-fn1}
\widetilde\Delta(x,m)=h(x,m)\frac{{\rm d}}{{\rm d}x}-\begin{pmatrix} \lambda(x,m) & 1 \\ \alpha(x,m)+q(m)h(x,m) & \lambda(x,m) \end{pmatrix},
\end{gather}
where $q(m)$ is an analytic germ such that
\begin{gather}\label{eq:RC-gammaQ}
\gamma(m)=-2\cos \pi\sqrt{1+4q(m)}.
\end{gather}
Let us remark that $\widetilde\Delta$ is meromorphic in $x\in\CP^1$ and has a regular singular point at infinity $($which is not Fuchsian unless $q(m)=0)$ when considered on the trivial vector bundle over $\CP^1$.

\item[$(ii)$] If $\gamma(0)\neq -2$, then $\Delta(x,m)$ is analytically equivalent to a germ at~$0$ of a parametric system $\widetilde\Delta'(x,m)$ given by
\begin{gather}\label{eq:RC-fn2}
\widetilde\Delta'(x,m)=h(x,m)\frac{{\rm d}}{{\rm d}x}-\begin{pmatrix} \lambda(x,m) & 1+b(m)x \\ \beta^{(0)}(m)+x\beta^{(1)}(m) & \lambda(x,m) \end{pmatrix},
\end{gather}
with
\begin{gather}\label{eq:RC-betab}
\beta^{(0)}(m)=\alpha^{(0)}+bh^{(0)}\beta^{(1)},\qquad
\beta^{(1)}(m)=\frac{\alpha^{(1)}-b\alpha^{(0)}}{1-bh^{(1)}+b^2h^{(0)}},
\end{gather}
where $b(m)$ is an analytic germ such that
\begin{gather}\label{eq:RC-gammab}
\gamma(m)=2\cos 2\pi\sqrt{b(m)\beta^{(1)}(m)}.
\end{gather}
Let us remark that $\widetilde\Delta'$ is meromorphic in $x\in\CP^1$ and has a Fuchsian singular point at infinity when considered on the trivial vector bundle over $\CP^1$.
\end{enumerate}
\end{Theorem}

\begin{Remark}In Theorem~\ref{theorem:RC-II}(i) $\gamma(0)=-2$ does not pose a problem since the function $t\mapsto -2\cos\pi\sqrt t$ appearing in the right side of \eqref{eq:RC-gammaQ} is analytically invertible near $t=0$
and one is free to choose $q(0)=-\frac{1}{4}$. Similarly, in {(ii)} for $\gamma(0)=2$ one can find an analytic germ $b$ with $b(0)=0$ satisfying \eqref{eq:RC-gammab}.
\end{Remark}

Theorem~\ref{theorem:RC-II} follows from Theorem~\ref{theorem:RC-I}(a) by a direct calculation of the invariants of the two parametric systems $\tilde\Delta$, $\tilde\Delta'$ using the following lemma.

\begin{Lemma}\label{lemma:RC-1}The analytic invariant $\gamma$ defined by~\eqref{eq:RC-gamma} of a system
\begin{gather}\label{eq:RC-4}
h(x)\frac{{\rm d}}{{\rm d}x}-\big[A^{(0)}+A^{(1)}x\big]=0,\qquad \text{with}\quad
A^{(k)}=\left(\begin{smallmatrix} a_{11}^{(k)} & a_{12}^{(k)} \\[3pt] a_{21}^{(k)} & a_{22}^{(k)} \end{smallmatrix}\right)
\end{gather}
and $ h(x)=x^2+h^{(1)}x+h^{(0)} $, is equal to
\begin{gather}\label{eq:RC-4g}
 \gamma=2\cos 2\pi\sqrt{\left(\frac{a_{11}^{(1)}-a_{22}^{(1)}}{2}\right)^2 +a_{12}^{(1)}a_{21}^{(1)}}.
\end{gather}
\end{Lemma}

\begin{proof}
This system considered on the trivial vector bundle over the Riemann sphere $\CP^1$ has singularities only at the zero points of $h(x)$ and at the point $x=\infty$.
Therefore the monodromy matrix $M$ in the formula \eqref{eq:RC-gamma}
\[\gamma={\rm e}^{-2\pi{\rm i}\lambda^{(1)}}\tr(M),\qquad\text{where}\quad \lambda^{(1)}=\frac{a_{11}^{(1)}+a_{22}^{(1)}}{2},\]
is also a matrix of monodromy around $x=\infty$ in the negative direction. In the coordinate $t=x^{-1}$ the system~\eqref{eq:RC-4} is equivalent to
\begin{gather}\label{eq:RC-4z}
t \big(1+h^{(1)}t+h^{(0)}t^2\big)\frac{{\rm d}}{{\rm d}t}+\big[A^{(1)}+A^{(0)}t\big]=0,
\end{gather}
which has a Fuchsian singularity at $t=0$.
The eigenvalues of its principal matrix $-A^{(1)}$ are $-\lambda^{(1)}\pm\sqrt D$ where $ D:=\big(\frac{a_{11}^{(1)}-a_{22}^{(1)}}{2}\big)^2+a_{12}^{(1)}a_{21}^{(1)}$.
Suppose first that the singularity is non-resonant, i.e., that $2\sqrt D\notin\Z$, in which case there exists a local analytic gauge transformation
$T(t)$ near $t=0$, that brings \eqref{eq:RC-4z} to the diagonal system
\begin{gather*}
t \frac{{\rm d}}{{\rm d}t}+\left(\begin{smallmatrix} \lambda^{(1)}+\sqrt D & 0 \\ 0 & \lambda^{(1)}-\sqrt D \end{smallmatrix}\right)=0,
\end{gather*}
(cf.~\cite[Chapter~16]{IlYa}), for which an associated diagonal fundamental solution has its monodromy matrix around $t=0$ in the negative direction equal to
\[M={\rm e}^{2\pi{\rm i}\lambda^{(1)}} \left(\begin{smallmatrix} {\rm e}^{2\pi{\rm i}\sqrt D} & 0 \\ 0 & {\rm e}^{-2\pi{\rm i}\sqrt D} \end{smallmatrix}\right).\]
Therefore $\gamma=2\cos 2\pi\sqrt D$.

The resonant case is a limit of non-resonant cases, and the formula~\eqref{eq:RC-4g} for~$\gamma$ remains valid, because the trace of monodromy
depends analytically on the coefficients of~$A$.
\end{proof}

\begin{proof}[Proof of Theorem~\ref{theorem:RC-II}] Use \eqref{eq:RC-lambdaalpha} to verify that $h(x)$, $\lambda(x)$ and $\alpha(x)$ are indeed the formal invariants of the system $\widetilde\Delta(x,m)$ \eqref{eq:RC-fn1} and $\widetilde\Delta'(x,m)$~\eqref{eq:RC-fn2}.

To verify \eqref{eq:RC-gammaQ}, set $Q:=\frac{1}{2}\big({-}1\pm\sqrt{1+4q}\big)$, so that $q=Q^2+Q, $ and let
$T(x):=\left(\begin{smallmatrix} 1 & 0 \\ Qx & 1 \end{smallmatrix}\right)$, then
\[T^*\widetilde\Delta(x,m)=h(x)\frac{{\rm d}}{{\rm d}x}-\begin{pmatrix} \lambda(x)+Qx & 1 \\ \alpha(x)+\big(h^{(0)}+h^{(1)}x\big)Q^2 & \lambda(x)-Qx \end{pmatrix}.\]
Now the system is in the form~\eqref{eq:RC-4}, and $ \gamma=2\cos 2\pi Q=-2\cos \pi\sqrt{1+4q}$ using \eqref{eq:RC-4g}.

The identity \eqref{eq:RC-gammab} follows directly from the formula~\eqref{eq:RC-4g}.
If $\gamma(0)\neq-2$, then the equation~\eqref{eq:RC-gammab} with $\beta^{(1)}(m)=\alpha^{(1)}(0)+O(m)$ given by~\eqref{eq:RC-betab}, $\alpha^{(1)}(0)\neq0$, has an analytic solution~$b(m)$ for small~$m$.
\end{proof}

\begin{Remark}A.~Bolibrukh showed that any irreducible system on a neighborhood of $0\in\C$ is analytically gauge equivalent on this neighborhood to a global system on
the trivial vector bundle over $\CP^1$ with a Fuchsian singularity at $\infty$ and no other additional singularities (see~\cite{Bo92,Bo94} or~\cite{Il}).
One can show that the restriction $\Delta_m$ of a parametric system $\Delta$ considered here to any parameter $m$ from some neighborhood of 0 is irreducible,
and that the global system to which it is then gauge equivalent by the Bolibrukh's theorem will necessarily have the form~\eqref{eq:RC-fn2}.
The problem of Theorem~\ref{theorem:RC-II}(ii) is to be able to do this analytically in~$m$ on a neighborhood of origin in the parameter space $\C^l$, which is where the condition $\gamma(0)\neq -2$ becomes necessary. Aside of the irreducibility condition, one of the essential ingredients in Bolibrukh's proof
is triangularity of the total monodromy matrix $M(m)$ (Definition~\ref{definition:RC-3}(ii)). If this matrix can be triangularized
analytically in $m$, then Bolibrukh's theorem holds also with local analytic dependence on $m$ (using the local rigidity of trivial vector bundles over~$\CP^1$ \cite[Corollary~5.4]{Sabbah}). Let us remark that, in particular, $\gamma(0)\neq\pm 2$ will insure analytic diagonalizability of~$M(m)$.

V.~Kostov has showed \cite{Ko} for a general system $\Delta_0(x)=\Delta(x,0)$ \eqref{eq:RC-Delta0}
in a Birkhoff normal form, that is a system $\Delta(x,0)=x^{k+1}\frac{{\rm d}}{{\rm d}x}-\big[A^{(0)}(0)+\dots+A^{(k)}(0)x^k\big]$,
whose eigenvalues of $A^{(k)}(0)=\mathop{\res}\limits_{x=0} \frac{A(x,0)}{x^{k+1}}$ do not differ by a non-zero integer,
that any its unfolding is analytically gauge equivalent to a parametric system in a generalized Birkhoff normal form $\Delta(x,m)=\big(x^{k+1}+h^{(k)}(m)x^k+\dots+h^{(0)}(m)\big)\frac{{\rm d}}{{\rm d}x}-\big[A^{(0)}(m)+\cdots+A^{(k)}(m)x^k\big]$.
Theorem~\ref{theorem:RC-II}(ii) confirms this for the parametric systems $\Delta$ studied here.
In the case of $\gamma(0)=-2$, if the system~$\Delta_0$ is in a Birkhoff normal form (which can always be assumed by the Bolibrukh's theorem),
then by~\eqref{eq:RC-4g} the eigenvalues of~$A^{(1)}$ differ by an odd integer and the condition of Kostov is violated.
Correspondingly, also the equation~\eqref{eq:RC-gammaQ} may fail to have an analytic solution with given~$q(0)$,
in which case the parametric family fails to be analytically equivalent to the generalized Birkhoff normal form~\eqref{eq:RC-fn2}.
\end{Remark}

\subsection{Formal theory}\label{sec:RC-1.2}

\begin{Proposition}[formal classification]\label{proposition:RC-formal}
A parametric system $\Delta(x,m)$ is formally equivalent to its formal normal form
\begin{gather}\label{eq:RC-model}
\widehat\Delta(x,m)=h(x,m)\frac{{\rm d}}{{\rm d}x}-\begin{pmatrix} \lambda(x,m) & 1 \\ \alpha(x,m) & \lambda(x,m) \end{pmatrix},
\end{gather}
by means of a unique formal gauge transformation in $(x,m)$
\[\hat T(x,m)=\sum_{j,|\mbf k|=0}^{+\infty} T^{(j,\mbf k)}x^j m^{\mbf k},\qquad m^{\mbf k}=m_1^{k_1}\cdots m_l^{k_l},\]
with $T^{(0,\mbf 0)}=I$. Generically, this series is divergent in both~$x$ and~$m$.

In this sense, two parametric systems $\Delta(x,m)$ and $\Delta'(x,m)$ are formally equivalent if and only if their formal invariants are the same:
$h=h'$, $\lambda=\lambda'$, $\alpha=\alpha'$.
\end{Proposition}

\begin{Remark}
Linear gauge transformations $T(x,m)$ commute with scalar matrix functions
\[T^*(\Delta-\lambda I)=T^* \Delta -\lambda I,\]
i.e., two systems $\Delta$, $\Delta'$ are analytically (resp.\ formally) gauge equivalent if and only if the systems $\Delta-\lambda I$, $\Delta'-\lambda I$ are.
Hence we can restrict our discussion to \emph{traceless systems} whose formal invariant
\[\lambda(x,m)=0.\]
\end{Remark}

\begin{Definition}[reduced invariants $\epsilon(m)$, $\mu(m)$]\label{definition:RC-reduction}
Let $\Delta(x,m)$ be a traceless parametric system whose formal invariants are $h(x,m)$, $\lambda(x,m)=0$ and $\alpha(x,m)$.
After an analytic translation and dilatation of the $x$-coordinate
\begin{gather*}
 x\mapsto \alpha^{(1)}x-\tfrac{h^{(1)}}{2}
\end{gather*}
and an introduction of new parameters
\begin{gather}\label{eq:RC-epsilonmu}
\epsilon(m)=\big(\tfrac{1}{\alpha^{(1)}}\big)^2 \big( \big(\tfrac{h^{(1)}}{2}\big)^2-h^{(0)} \big),\qquad
\mu(m)=\tfrac{\alpha^{(0)}}{(\alpha^{(1)})^2}-\tfrac{h^{(1)}}{2\alpha^{(1)}}
\end{gather}
we obtain a parametric system with whose formal invariants are
\begin{gather}\label{eq:RC-rfi}
h(x,m)=x^2 -\epsilon(m),\qquad \lambda(x,m)=0,\qquad \alpha(x,m)=\mu(m)+ x.
\end{gather}
By Proposition~\ref{proposition:RC-prenormal}, this system can be written up to analytic equivalence as
\begin{gather}\label{eq:RC-prenormal}
\Delta(x,m)=\big(x^2 -\epsilon\big)\frac{{\rm d}}{{\rm d}x}-\begin{pmatrix} 0 & 1 \\ \mu+x+\big(x^2 -\epsilon\big) r(x,m) & 0 \end{pmatrix}.
\end{gather}
\end{Definition}

\textit{To simplify the discussion, from now on we will assume that $\Delta(x,m)$ is in the form \eqref{eq:RC-prenormal},
and correspondingly its formal normal form of Proposition~{\rm \ref{proposition:RC-formal}} is}
\begin{gather}\label{eq:RC-model2}
 \widehat\Delta(x,m)=\big(x^2 -\epsilon\big)\frac{{\rm d}}{{\rm d}x}-\begin{pmatrix} 0 & 1 \\ \mu+x & 0 \end{pmatrix}.
\end{gather}

\begin{proof}[Proof of Proposition \ref{proposition:RC-formal}]
Let $\Delta(x,m)$ be a parametric system in the prenormal form~\eqref{eq:RC-prenormal}. We will show that there exists a formal gauge transformation $\hat T(x,m)$ in form of a power series in $(x,\mu,\epsilon)$ whose coefficients depends analytically on $m$, that brings $\Delta(x,m)$ to the reduced formal normal form $\widehat\Delta(x,m)$~\eqref{eq:RC-model2}. We shall be looking for $\hat T$ written as
\[\hat T(x,m)=a(x,m) I + b(x,m)\begin{pmatrix} 0 & 1 \\ \mu+x & 0 \end{pmatrix}+\big(x^2 -\epsilon\big)\begin{pmatrix} 0 & 0 \\ c(x,m) & d(x,m) \end{pmatrix}. \]
We want that $\widehat\Delta=\hat T^*\Delta$, which means
\[\left[ \begin{pmatrix} 0 & 1 \\ \mu+x & 0 \end{pmatrix}, \begin{pmatrix} 0 & 0 \\ c(x,m) & d(x,m) \end{pmatrix}\right]+\begin{pmatrix} 0 & 0 \\ r(x,m) & 0 \end{pmatrix}\cdot\hat T(x,m)=
\frac{{\rm d}\hat T(x,m)}{{\rm d}x},\]
where $[\cdot ,\cdot ]$ stands for the commutator of matrices. This gives a system of equations
\begin{gather}
c=a', \label{eq:RC-abcd1} \\
d=b', \label{eq:RC-abcd2} \\
-(\mu+x)d+ar=b+(\mu+x)b'+2xc+\big(x^2 -\epsilon\big)c', \label{eq:RC-abcd3} \\
-c+br=a'+2xd+\big(x^2 -\epsilon\big)d', \label{eq:RC-abcd4}
\end{gather}
where $'$ stands for the (formal) derivative w.r.t.~$x$.
Substituting \eqref{eq:RC-abcd1} and \eqref{eq:RC-abcd2} in~\eqref{eq:RC-abcd3} and~\eqref{eq:RC-abcd4} gives
\begin{gather}
b+2(\mu+x)b' =ar-2xa'-\big(x^2 -\epsilon\big)a'', \label{eq:RC-abcd5} \\
2a' =br-2xb'-\big(x^2 -\epsilon\big)b''. \label{eq:RC-abcd6}
\end{gather}
Writing
\begin{gather*}
a(x,m)=\sum_{(j,k,l)}a_{j,k,l}(m)\mu^j\epsilon^kx^l,\qquad
b(x,m)=\sum_{(j,k,l)}b_{j,k,l}(m)\mu^j\epsilon^kx^l, \\
r(x,m)=\sum_{l}r_l(m)x^l,
\end{gather*}
and identifying the coefficients of the term $\mu^j\epsilon^kx^l$ in~\eqref{eq:RC-abcd5} and~\eqref{eq:RC-abcd6} shows that
\begin{gather*}
 (2l+1) b_{j,k,l}+2(l+1) b_{j-1,k,l+1} \ \text{is a finite linear combination}\\
 \hspace*{53mm}~\text{of} \ a_{\tilde j, \tilde k, \tilde l} , \ {\big(\tilde j,\tilde k, \tilde l\big)}\leq_{\rm LEX} {(j,k,l)},\\
 2(l+1) a_{j,k,l+1} \ \text{is a finite linear combination of} \ b_{\tilde j, \tilde k, \tilde l} , \ {\big(\tilde j,\tilde k, \tilde l\big)}\leq_{\rm LEX} {(j,k,l)},
\end{gather*}
where $\leq_{\rm LEX}$ is the lexicographic ordering on $\N^3$.
There is no constraint on the coefficients $a_{j,k,0}$, which we choose $0$ for $(j,k)\neq(0,0)$, and $a_{0,0,0}=1$.
All the coefficients are now uniquely determined through a transfinite recursion with respect to the $\leq_{\rm LEX}$-ordering, which is a well-ordering on the index set $\N^3$.
\end{proof}

\subsection{Stokes phenomenon and confluence}\label{sec:RC-1.3}

Let us assume that the system $\Delta(x,m)$ is traceless in the form \eqref{eq:RC-prenormal}.
The reduced formal invariants $\epsilon(m)$ and $\mu(m)$ \eqref{eq:RC-epsilonmu} are responsible for two basic qualitative changes of the system:
\begin{itemize}\itemsep=0pt
\item[--] $\epsilon(m)$ corresponds to separation of the double singularity ($\epsilon=0$) into two simple (Fuchsian) ones ($\epsilon\neq0$),
\item[--] $\mu(m)$ corresponds, when $\epsilon(m)=0$ and the singularity is irregular, to separation of the double eigenvalue ($\mu=0$) of $A(0,m)$ into two simple ones $(\mu\neq 0)$, hence to the disappearance of resonance of the irregular singularity. This in turn leads to a change of order of Borel summability of formal normalizing transformations (Remark~\ref{remark:RC-sectoralnormalization} below).
\end{itemize}

In this section we shall describe the effect that these two qualitative changes have on the structure of the solution space and on the Stokes phenomenon.
We will explain them in terms of the ``Stokes geometry'' of the natural domains of normalization,
which are associated, as we shall see, to the meromorphic quadratic differential
\begin{gather}\label{eq:RC-quadraticdifferential}
\frac{\alpha(x,m)}{h(x,m)^2}({\rm d}x)^2=\frac{\mu+x}{(x^2-\epsilon)^2}({\rm d}x)^2.
\end{gather}
This differential is the negative of the determinant of the meromorphic ``Higgs field''
\[\frac{1}{x^2-\epsilon}\begin{pmatrix} 0 & 1 \\ \mu+x & 0 \end{pmatrix}{\rm d}x\]
associated to the formal normal form~\eqref{eq:RC-model2} on the trivial vector bundle,
the matrix of which has eigenvalues $\pm\frac{\sqrt{\mu+x}}{x^2-\epsilon}$. The point $x=-\mu$, at which the two eigenvalues merge
is a ``turning point'' in the terminology of exact WKB analysis~\cite{KT}. This point is not a singularity of the differential system, but it is a~``spectral'' singularity of the matrix of the formal normal form system, and will play an equally important role in our description.

Let the multivalued function
\begin{gather}\label{eq:RC-Theta}
\Theta(x,m):=\int_\infty^{x}\frac{ \sqrt{\alpha(x,m)}}{h(x,m)} {\rm d}x=\int_\infty^{x}\frac{ \sqrt{\mu+x}}{x^2-\epsilon} {\rm d}x,
\end{gather}
with ramification points at the zero locus of $h(x,m)=x^2-\epsilon$ and of $\alpha(x,m)=\mu+x$,
be the rectifying coordinate for the quadratic differential~\eqref{eq:RC-quadraticdifferential} which then becomes~$({\rm d}\Theta)^2$.
The different solutions of the system $\Delta(x,m)$ are expected to have an asymptotic behavior near the singular points of order that is of a combination of
$(\mu+x)^{-\frac14}{\rm e}^{\Theta(x,m)}$ and $(\mu+x)^{-\frac14}{\rm e}^{-\Theta(x,m)}$.
This will be made explicit in the following remark, which summarizes some classical results
on the local behavior of the solutions of~$\Delta_m(x)$ and local normal forms near each of its singular points for all fixed values of parameter~$m$.
The general theory of singularities of linear differential systems has been developed by G.D.~Birkhoff, W.J.~Trjitzinski, J.~Malmquist, M.~Hukuhara, H.~Turrittin, Y.~Sibuya and many others. Reader familiar with basics of this theory may skip this Remark~\ref{remark:RC-sectoralnormalization} and go straight to
Theorem~\ref{theorem:RC-IV} below, which shows how these disparate descriptions for different values of $\epsilon$, $\mu$ fit into a single parametric picture.

\begin{Remark}[canonical solution bases and sectoral normalization of $\Delta_m(x)$]\label{remark:RC-sectoralnormalization}
Let $\Delta(x,m)$ as in~\eqref{eq:RC-prenormal} be analytic on some polydisc $\sX\times\sM\subseteq\C\times\C^l$, where $\sM$ is small enough so
that both roots of $h(x,m)=x^2 -\epsilon$ are in $\sX$ for all $m\in\sM$.
As before, let $\Delta_m(x)$ denote the restriction of $\Delta(x,m)$ to the fixed value of~$m$.
Depending on $\epsilon(m)$ and $\mu(m)$, there are the following four possible situations (see, e.g., \cite{Balser,IlYa,Sib,Wa}):
\begin{itemize}\itemsep=0pt
\item[(a)] $\epsilon=\mu=0$:
The restricted system $\Delta_m$ has a \emph{resonant irregular singularity} at the origin
and a formal fundamental solution matrix
\[\hat Y_{O,m}(x)=\hat H_{O,m}(x)
\left(\begin{smallmatrix} (\mu+x)^{-\frac{1}{4}} & 0 \\ 0 & (\mu+x)^{\frac{1}{4} } \end{smallmatrix}\right)
\tfrac{{\rm i}}{\sqrt2}\left(\begin{smallmatrix} 1 & 1 \\ 1 & -1 \end{smallmatrix}\right)
\left(\begin{smallmatrix} {\rm e}^{\Theta(x,m)} & 0 \\ 0 & {\rm e}^{-\Theta(x,m)} \end{smallmatrix}\right),
\]
where $\hat H_{O,m}(x)$ is a formal power series in~$x$ with matrix coefficients, that is Borel $\frac{1}{2}$-summable in all directions except of the direction $\arg x=0+2\pi\Z$, which is tangent to the curves $\Im\Theta(x,m)=0$ at the origin.
Associated to $\hat H_{O,m}(x)$ is its Borel sum which is a unique bounded sectoral gauge transformation $H_{O,m}(x)$ defined on a ramified sector
\begin{gather}\label{eq:RC-sectorSO}
S_{O,m}=\{x\in\sX \,|\, |\arg x+\pi|<2\pi-\eta\},
\end{gather}
with $\eta>0$ arbitrarily small.
The system $\Delta_m$ is formally equivalent to $\widehat\Delta_m$ \eqref{eq:RC-model2} by means of some formal Borel $\frac{1}{2}$-summable gauge transformation $\hat T_{m}(x)$ whose Borel sum is defined on the same sector $S_{O,m}$~\cite{JLP2}.

\item[(b)] $\epsilon=0$, $\mu\neq 0$:
The restricted system $\Delta_m$ has a \emph{non-resonant irregular singularity} at the origin
and a formal fundamental solution matrix
\[\hat Y_{I,m}(x)=\hat H_{I,m}(x)
\left(\begin{smallmatrix} (\mu+x)^{-\frac{1}{4}} & 0 \\ 0 & (\mu+x)^{\frac{1}{4} } \end{smallmatrix}\right)
\tfrac{{\rm i}}{\sqrt2}\left(\begin{smallmatrix} 1 & 1 \\ 1 & -1 \end{smallmatrix}\right)
\left(\begin{smallmatrix} {\rm e}^{\Theta(x,m)} & 0 \\ 0 & {\rm e}^{-\Theta(x,m)} \end{smallmatrix}\right),
\]
where $\hat H_{I,m}(x)$ is a formal power series in~$x$ with matrix coefficients, that is Borel $1$-summable in all directions except of the directions $\arg x=\arg\sqrt\mu+\pi\Z$, which are tangent to the curves $\Im\Theta(x,m)=0$ at the origin.
Associated to $\hat H_{I,m}(x)$ are its Borel sums which are a unique pair of bounded sectoral gauge transformations $H_{I\pm,m}(x)$ defined on a pair of sectors
\begin{gather}\label{eq:RC-sectorSI}
S_{I\pm,m}=\big\{x\in\sX \,|\, |\arg x-\arg\sqrt\mu\pm\tfrac{\pi}{2}|<\pi-\eta\big\},
\end{gather}
with $\eta>0$ arbitrarily small.
The system $\Delta_m$ is formally equivalent to $\widehat\Delta_m$~\eqref{eq:RC-model2} by means of some formal Borel $1$-summable gauge transformation $\hat T_{I,m}(x)$ whose Borel sums are defined on the same pair of sectors $S_{I\pm,m}$~\cite{JLP1}.

\item[(c)] $\epsilon\neq 0$: The restricted system has \emph{two Fuchsian singularities} at $x_1=\sqrt\epsilon$ and $x_2=-\sqrt\epsilon$.
Supposing that $m$ is such that the Fuchsian singularity at $x_i$ is \emph{non-resonant}, i.e., that $ \frac{\sqrt{\mu+x_i}}{x_i}\notin\Z$,
then there exists a fundamental solution matrix
\[Y_{i,m}(x)= H_{i,m}(x)
\left(\begin{smallmatrix} (\mu+x)^{-\frac{1}{4}} & 0 \\ 0 & (\mu+x)^{\frac{1}{4} } \end{smallmatrix}\right)
\tfrac{{\rm i}}{\sqrt2}\left(\begin{smallmatrix} 1 & 1 \\ 1 & -1 \end{smallmatrix}\right)
\left(\begin{smallmatrix} {\rm e}^{\Theta(x,m)} & 0 \\ 0 & {\rm e}^{-\Theta(x,m)} \end{smallmatrix}\right),
\]
where $H_{i,m}(x)$ is a convergent formal power series in $x - x_i$ with matrix coefficients,
the sum of which is defined on a neighborhood of~$x_i$
\[S_{i,m}=\big\{x\in\sX \,|\, |x-x_i|<2\sqrt{|\epsilon|}\big\}.\]
The system $\Delta_m$ is locally equivalent to $\widehat\Delta_m$ \eqref{eq:RC-model2} by a convergent transformation $T_{I,m}(x)$ on $S_{i,m}$.
The local gauge transformations $H_{i,m}(x)$ and $T_{i,m}(x)$ depend analytically on $m$ as long as the singularity stays non-resonant Fuchsian.

\item[(d)] $\epsilon\neq 0$: If a Fuchsian singularity at $x_i$ is \emph{resonant}, $\frac{\sqrt{\mu+x_i}}{x_i}=k\in\Z\setminus\{0\}$,
then there exists a~fundamental solution matrix
\begin{gather*}\tilde Y_{i,m}=\begin{cases}
\tilde H_{i,m}(x)
\left(\!\begin{smallmatrix} (\mu+x)^{-\frac{1}{4}} & 0 \\ 0 & (\mu+x)^{\frac{1}{4} } \!\end{smallmatrix}\right)
\tfrac{{\rm i}}{\sqrt2}\left(\!\begin{smallmatrix} 1 & 1 \\ 1 & -1 \end{smallmatrix}\!\right)
\left(\!\begin{smallmatrix} {\rm e}^{\Theta(x,m)} & 0 \\ 0 & {\rm e}^{-\Theta(x,m)} \end{smallmatrix}\!\right)
(x-x_i)^{N},&\text{if}\ k\neq 0,\\
\tilde H_{i,m}(x)(x-x_i)^{N},&\text{if}\ k=0,
\end{cases}\end{gather*}
where $\tilde H_{i,m}(x)$ is a convergent formal power series in $x - x_i$ with matrix coefficients on the neighborhood $S_{i,m}$ of $x_i$,
and $N$ is an upper/lower triangular nilpotent matrix.
The system $\Delta_m$ is locally equivalent to $\widehat\Delta_m$~\eqref{eq:RC-model2} by a convergent transformation of the form $\tilde T_{I,m}(x)=\tilde T_{i,m}'(x)+(x-x_i)^k\log(x-x_i) \tilde T_{i,m}''(x)$.
\end{itemize}
\end{Remark}

\begin{figure}[t]\centering
 \begin{tikzpicture}
\node at (-5.4,0){\includegraphics[width=0.3\textwidth]{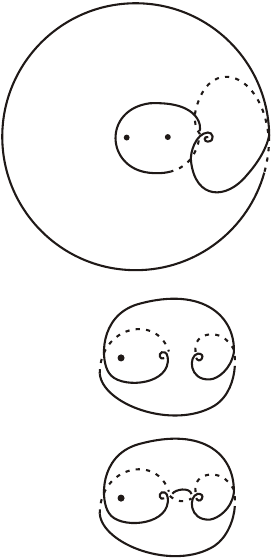}};

\node at (0,1.1){\includegraphics[width=0.3\textwidth]{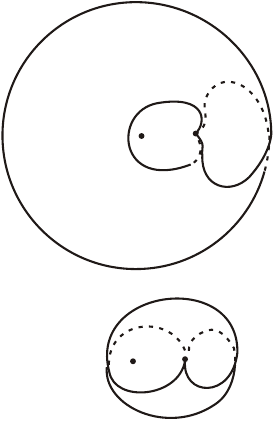}};

\node at (5.4,1.7){\includegraphics[width=0.3\textwidth]{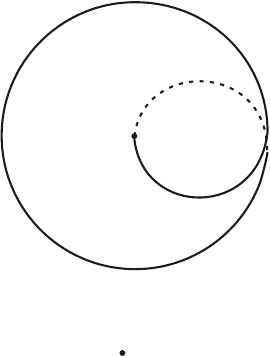}};
\node at (-5.4,-5.35){(a) $\mu>\sqrt{\epsilon}>0$};
\node at (0,-5.35){(b) $\mu>\sqrt{\epsilon}=0$};
\node at (5.4,-5.35){(c) $\mu=\sqrt{\epsilon}=0$};
\node at (-7,5){$\rXO(\check\mu,\check\epsilon)$};
\node at (-1.6,5){$\rXO(\check\mu,0)$};
\node at (3.8,5){$\rXO(0,0)$};
\node at (-5.3,2.3){$-\mu$};
\node at (-4.9,2.8){$-\sqrt{\epsilon}$};
\node at (-3.8,2.7){$\sqrt{\epsilon}$};
\node at (0.2,2.2){$-\mu$};
\node at (0.8,2.4){$0$};
\node at (5.2,2.4){$0$};
\node at (-6.6,-0.8){$\rXI(\check\mu,\check\epsilon)$};
\node at (-1.2,-0.8){$\rXI{}_+(\check\mu,0)$};
\node at (-1.2,-2.2){$\rXI{}_-(\check\mu,0)$};
\node at (4.3,-0.8){$\rXI(0,0)$};
\node at (-5.6,-1.55){$-\mu$};
\node at (-5.4,-1.2){$-\sqrt{\epsilon}$};
\node at (-4.0,-1.4){$\sqrt{\epsilon}$};
\node at (-0.2,-1.75){$-\mu$};
\node at (0.6,-1.5){$0$};
\node at (5.4,-1.4){$0$};
\node at (-6.6,-2.9){$\rXI{}_+(\check\mu,\check\epsilon)$};
\node at (-6.6,-4.7){$\rXI{}_-(\check\mu,\check\epsilon)$};
\node at (-5.6,-4.05){$-\mu$};
\node at (-5.4,-3.7){$-\sqrt{\epsilon}$};
\node at (-4.0,-3.9){$\sqrt{\epsilon}$};
\end{tikzpicture}
\caption{Examples of the outer and inner domains $\rXO(\epsilon,\mu)$, $\rXI(\epsilon,\mu)$ of Theorem~\ref{theorem:RC-III} for selected values of~$(\mu,\epsilon)$.}
	\label{figure:RC-OI}
\end{figure}
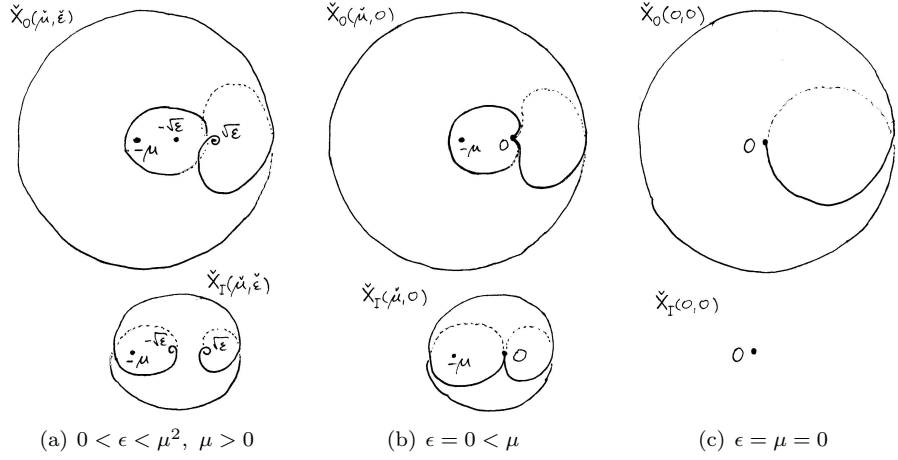

The \emph{change of order of summability} of the formal normalizing gauge transformations in between the cases (a) and (b) of Remark~\ref{remark:RC-sectoralnormalization} is a phenomenon that has not been studied previously.
In the following Theorem~\ref{theorem:RC-IV} and its Corollary~\ref{theorem:RC-III} it is explained by the form and organization of the domains on which there naturally exist certain canonical solution bases (``mixed bases'') as well as bounded normalizing gauge transformations.
What happens is that the sector $S_{O,m}$ \eqref{eq:RC-sectorSO} for $\epsilon=\mu=0$ unfolds to an ``\emph{outer domain}'' $\rXO(\mu,\epsilon)$, while also a new
pair of ``\emph{inner domains}'' $\rXI{}_\pm(\mu,\epsilon)$ appears for $(\mu,\epsilon)\neq (0,0)$ that for $\epsilon=0$ will locally agree with the pair of sectors~$S_{I\pm,m}$~\eqref{eq:RC-sectorSI}.
When $\epsilon=0$ and $\mu\to 0$, the inner domains $\rXI{}_\pm(\mu,\epsilon)$ will shrink and disappear (Fig.~\ref{figure:RC-OI}). As we shall see, their disappearance is caused by the coalescence of the turning point $x=-\mu$ and of the irregular singularity at $x=0$ of the quadratic differential~\eqref{eq:RC-quadraticdifferential}.

\begin{Definition}[bounded analytic functions on ramified parametric domains]\label{definition:RC-CalB}
A ramified parametric domain $\Omega$ over the $(x,m)$-space $\C\times\C^l$ with a ramification locus
\[\Sigma:=\big\{x^2-\epsilon(m)=0\big\}\cup\big\{\epsilon(m)\big(\mu(m)^2-\epsilon(m)\big)=0\big\}\]
is a connected topological set such that $\Omega\setminus\Sigma$ is an open subdomain of a covering space of $\big(\C\times\C^l\big)\setminus\Sigma$.
We will denote $(\check x,\check m)$ the ``ramified'' coordinates on $\Omega$ that are the lifting of the coordinates $(x,m)$,
and consider $\Omega$ as a parametric family of ramified domains over the $x$-plane depending on a parameter $\check m$
\[\Omega=\coprod_{\check m\in\rM}\Omega(\check m),\qquad \Omega(\check m)=\{\check x\,|\, (\check x,\check m)\in\Omega\},\]
where the interior of each
$\Omega(\check m)$ is a subdomain of a covering space of $\C\setminus\big\{x^2-\epsilon(m)=0\big\}$,
and the interior of $\rM$ is a subdomain of a covering space of $\C^l\setminus\big\{\epsilon(m)\big(\mu(m)^2-\epsilon(m)\big)=0\big\}$.
We allow for the ramifying locus $\Sigma$ (resp.~$\big\{\epsilon(m)\big(\mu(m)^2-\epsilon(m)\big)=0\big\}$) to be included in $\Omega$ (resp.~$\rM$) since we want to be able to cover a~full neighborhood of the origin in $\C\times\C^l$ (resp.~$\C^l$).

For a function $f\colon \Omega\to\C$, we write \[f\in{\mathcal B}(\Omega)\]
if $f$ is bounded continuous on $\Omega$ and analytic on its interior $\Omega\setminus\Sigma$ and at the same time
$f(\cdot,\check m)$ is analytic on the interior of $\Omega(\check m)$ for each $\check m\in \rM$.
\end{Definition}

\begin{figure}[t]\centering
 \begin{tikzpicture}
 \node at (0,0){\includegraphics{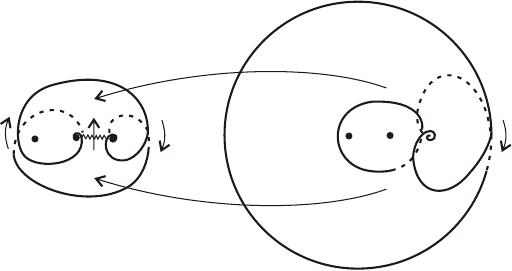}};
 \node at (-3.0,1.2){$Y_{I+}$};
 \node at (-3.0,-1.3){$Y_{I-}$};
 \node at (-1.0,1.25){$C_3$};
 \node at (-3.55,-0.25){$-\mu$};
\node at (-3.45,0.1){$-\sqrt{\epsilon}$};
\node at (-2.1,0.0){$\sqrt{\epsilon}$};
\node at (-2.6,-0.5){$N_1$};
\node at (-1.1,0.0){$C_1$};
\node at (-4.7,0.0){$NC_2$};
\node at (0.8,0.0){$Y_O$};
\node at (-1.0,-1.5){${\rm i}\left(\begin{smallmatrix}0& 1 \\ 1& 0\end{smallmatrix}\right)C_4$};
 \node at (1.9,-0.25){$-\mu$};
\node at (2.25,0.3){$-\sqrt{\epsilon}$};
\node at (3.4,0.0){$\sqrt{\epsilon}$};
\node at (5.9,0){$M=-{\rm i}C_0\left(\begin{smallmatrix}0& 1 \\ 1& 0\end{smallmatrix}\right)$};
 \end{tikzpicture}
\caption{Connection matrices between the fundamental solution matrices $Y_O$, $Y_{I\pm}$, where $N_1$, $N$ are as in~\eqref{eq:RC-N1N2}, and~$C_i$, $i=0,\dots,4$, are as in \eqref{eq:RC-Ci} in Lemma~\ref{lemma:RC-4} with $\kappa$ \eqref{eq:RC-kappax}.}	\label{figure:RC-YOI}
\end{figure}

\begin{Theorem}[fundamental solution matrix] \label{theorem:RC-IV}
Let $\Delta(x,m)$ be a germ of a parametric system with invariants $h(x,m)$, $\lambda(x,m)$, $\alpha(x,m)$ and $\gamma(m)$, and let $\epsilon(m)$, $\mu(m)$ be its reduced formal invariants~\eqref{eq:RC-epsilonmu}. There exists a covering of a neighborhood of~$0$ in the $(x,m)$-space $\C\times\C^l$ by three ramified parametric domains:
an outer domain $\rXO$ and a pair of inner domains $\rXI{}_\pm$ $($see Fig.~{\rm \ref{figure:RC-OI})} the form of which depends only on $(\epsilon(m),\mu(m))$,
and there exist bounded analytic matrix functions $H_\bullet\in\Mat_2\big({\mathcal B}\big(\rX_\bullet\big)\big)$, $\rX_\bullet=\rXO,\rXI{}_\pm$, on these domains, such that
the system~$\Delta(x,m)$ has canonical fundamental solution matrices~$Y_\bullet$ of the form
\begin{gather*}
Y_\bullet(x,m)=H_\bullet(x,m) \cdot
\left(\begin{smallmatrix} \alpha(x,m)^{-\frac{1}{4}} & 0 \\ 0 & \alpha(x,m)^{\frac{1}{4} } \end{smallmatrix}\right)
\tfrac{{\rm i}}{\sqrt2}\left(\begin{smallmatrix} 1 & 1 \\ 1 & -1 \end{smallmatrix}\right)
\left(\begin{smallmatrix} {\rm e}^{\Theta_\bullet(x,m)} & 0 \\ 0 & {\rm e}^{-\Theta_\bullet(x,m)} \end{smallmatrix}\right)
{\rm e}^{ \int_\infty^x \frac{\lambda(x,m)}{h(x,m)}{\rm d}x},
\end{gather*}
where $\Theta_\bullet$ is a branch of~\eqref{eq:RC-Theta} on $\rX_\bullet$, and $H_\bullet^{-1}$ has a pole at the zero of $\alpha(x,m)$ if $(\mu(m),\epsilon(m))\allowbreak \neq 0$.
The form and the construction of the domains $\rXO$, $\rXI{}_\pm$ will be detailed below.

The connection matrices between these fundamental solutions for each fixed $\check m\in\rM$ are as in Fig.~{\rm \ref{figure:RC-YOI}},
with~$N_1$, $N$ given by \eqref{eq:RC-N1N2}, $C_i$ as in~\eqref{eq:RC-Ci} with $\kappa$~\eqref{eq:RC-kappax}.
In particular, the monodromy matrix of~$Y_O(x,m)$ around both singular points equals $M=\left(\begin{smallmatrix} \gamma & -{\rm i} \\ -{\rm i} & 0 \end{smallmatrix}\right)$.
\end{Theorem}

The linear differential system satisfied by
\begin{gather*}
\left(\begin{smallmatrix} \alpha(x,m)^{-\frac{1}{4}} & 0 \\ 0 & \alpha(x,m)^{\frac{1}{4} } \end{smallmatrix}\right)
\tfrac{{\rm i}}{\sqrt2}\left(\begin{smallmatrix} 1 & 1 \\[4pt] 1 & -1 \end{smallmatrix}\right)
\left(\begin{smallmatrix} {\rm e}^{\Theta_\bullet(x,m)} & 0 \\ 0 & {\rm e}^{-\Theta_\bullet(x,m)} \end{smallmatrix}\right)
{\rm e}^{ \int_\infty^x \frac{\lambda(x,m)}{h(x,m)}{\rm d}x}
\end{gather*}
is
\begin{gather*}
\widebar\Delta(x,m)=h(x,m)\frac{{\rm d}}{{\rm d}x}-\left[\lambda(x,m)I+\left(\begin{smallmatrix} 0 & 1 \\ \alpha(x,m) & 0 \end{smallmatrix}\right)-
\frac{h(x,m)}{4\alpha(x,m)} \left(\begin{smallmatrix} 1 & 0 \\ 0 & -1 \end{smallmatrix}\right) \right],
\end{gather*}
which has an additional singularity at the zero point of $\alpha(x,m)$ and therefore does not belong to the considered class of parametric systems.
It does however play an intermediate role in comparing two parametric systems: if $H_\bullet$ and $H_\bullet'$ are the gauge transformations of Theorem~\ref{theorem:RC-IV} for two formally equivalent parametric systems $\Delta$ and $\Delta'$ on the same domain $\rX_\bullet$,
then their composition $T_\bullet=H_\bullet (H_\bullet')^{-1}$ is a gauge transformation such that $T^*\Delta=\Delta'$,
and such that $T_\bullet,T_\bullet^{-1}\in\GL_2\big({\mathcal B}\big(\rX_\bullet\big)\big)$ are both bounded and analytic on $\rX_\bullet$, including the zero point of $\alpha(x,m)$.

\begin{Corollary}[unfolded sectoral normalization] \label{theorem:RC-III}
Let $\Delta(x,m)$ be a germ of a parametric system \eqref{eq:RC-prenormal1} $($resp.~\eqref{eq:RC-prenormal}$)$,
and let $\rX_\bullet=\rXO, \rXI{}_\pm$ be the parametric domains of Theorem~{\rm \ref{theorem:RC-IV}} covering a full neighborhood of $0\in\C\times\C^l$.
There exist normalizing gauge transformations $T_\bullet\in\GL_2\big({\mathcal B}\big(\rX_\bullet\big)\big)$
bounded and analytic on these domains, that transform the parametric system $\Delta(x,m)$ to its formal normal form \eqref{eq:RC-model} $($resp.~\eqref{eq:RC-model2}$)$:
\[T_\bullet^*\Delta=\widehat\Delta,\qquad\bullet=O,I\pm.\]
When $\epsilon(m)\big(\mu(m)^2-\epsilon(m)\big)\neq 0$, the transformations $T_{I\pm}$ are both restrictions of the sane transformation $T_I\in\GL_2(\Cal B(\rXI))$
 defined on a domain $\rXI=\rXI{}_+\cup\rXI{}_-$ $($see Fig.~{\rm \ref{figure:RC-OI}(a))}.
\end{Corollary}

The form of the ramified parametric domains $\rX_\bullet=\rXO,\rXI{}_\pm$ depends only on the reduced formal invariants $\epsilon(m),\mu(m)$ \eqref{eq:RC-epsilonmu}.
Namely, they are parametric families of ramified domains $\rX_\bullet(\check\mu(\check m),\check\epsilon(\check m))$ over the $x$-plane,
\[\rX_\bullet=\coprod_{\check m\in\rM} \rX_\bullet(\check\mu(\check m),\check\epsilon(\check m)),\]
defined over a ramified domain $\rM$ covering a full neighborhood $\sM$ of $0$ in the parameter space of $m$
\[
\xymatrix{\rM\ni \check m \ \ar@{|->}[r]\ar@{|->}[d] & \ (\check\mu(\check m),\check\epsilon(\check m))\ar@{|->}[d] \\
	\sM\ni m \ \ar@{|->}[r] & \ (\mu(m),\epsilon(m))}
\]
(checked symbols denoting ramified variables/domains).
They are of two kinds (see Fig.~\ref{figure:RC-OI}):
\begin{itemize}\itemsep=0pt
	\item[(a)] The outer domain $\rXO(\check\mu(\check m),\check\epsilon(\check m))$ is doubly attached to the singularity $x_1=\sqrt\epsilon$. For $(\mu,\epsilon)=(0,0)$ it becomes a ramified sector $\rXO(0,0)$ at the origin of opening $>2\pi$,
	in which case $T_O(\cdot ,\check m)=T_{O,m}$ of Remark~\ref{remark:RC-sectoralnormalization}.
	\item[(b)] The inner domains $\rXI{}_\pm(\check\mu(\check m),\check\epsilon(\check m))$ are two parts of a single ramified parametric domain
	$\rXI={\rXI}{}_+\cup{\rXI}{}_-$ split in two by a cut between the singularities $x_1=\sqrt\epsilon$ and $x_2=-\sqrt\epsilon$, to which they are both attached.
	 For $\epsilon = 0$, $\mu \neq 0$, they become a pair of sectors $\rXI{}_\pm(\check\mu(\check m),0)$ of opening $>\pi$ at the origin, of Remark~\ref{remark:RC-sectoralnormalization}.
	For $\mu^2=\epsilon$ they shrink to a single point $\rXI{}_\pm(\check\mu,\check\epsilon)=\{-\mu\}$.
\end{itemize}

They are a close analogue of the ``Stokes domain'' in exact WKB analysis of second-order linear ODEs~\cite{KT}. Their construction is roughly the following:

\begin{figure}[t]\centering
 \begin{tikzpicture}
\node at (-6.1,0){\includegraphics[width=0.23\textwidth]{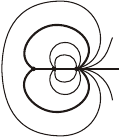}};
\node at (-2.1,0){\includegraphics[width=0.23\textwidth]{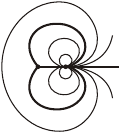}};
\node at (2.1,0){\includegraphics[width=0.23\textwidth]{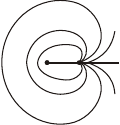}};
\node at (6.1,0){\includegraphics[width=0.23\textwidth]{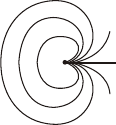}};
\node at (-7.3,0){\footnotesize $-\mu$};
\node at (-6.75,-0.2){\footnotesize $-\sqrt{\epsilon}$};
\node at (-5.9,-0.2){\footnotesize $\sqrt{\epsilon}$};
\node at (-3.2,0.0){\footnotesize$-\mu$};
\node at (-1.9,-0.23){\footnotesize $0$};
\node at (1.74,-0.2){\footnotesize$-\mu$};
\node at (2.45,-0.25){\footnotesize $\sqrt{\epsilon}$};
\node at (6.2,-0.2){\footnotesize $0$};
\node at (-6.2,-2.4){(a) $\mu>\sqrt{\epsilon}>0$};
\node at (-2.3,-2.4){(b) $\mu>\sqrt\epsilon=0$};
\node at (2.0,-2.4){(c) $\mu=\sqrt{\epsilon}>0$};
\node at (5.9,-2.4){(d) $\mu=\sqrt\epsilon=0$};
 \end{tikzpicture}
\caption{Examples of the horizontal foliation of the quadratic differential \eqref{eq:RC-quadraticdifferential}, with the separating trajectories between the outer and inner half-zones in bold, for selected values of $(\mu,\epsilon)$. In (a) the dashed bold trajectory between $-\sqrt{\epsilon}$ and $\sqrt{\epsilon}$ splits the inner zone into 2 half-zones corresponding to those in (b).} \label{figure:RC-zones1}
\end{figure}

The \emph{horizontal foliation} associated to the quadratic differential
\begin{gather}\label{eq:RC-quaddiffomega}
{\rm e}^{-2{\rm i}\omega}\frac{\alpha(x)}{h(x)^2}({\rm d}x)^2,\qquad \omega\in \left] -\frac{\pi}{2},\frac{\pi}{2}\right[,
\end{gather}
is one whose leaves are the real time trajectories of
\begin{gather}\label{eq:RC-horizontal}
\frac{{\rm d}x}{{\rm d}t}={\rm e}^{{\rm i}\omega}\frac{h(x)}{\sqrt{\alpha(x)}}={\rm e}^{{\rm i}\omega}\frac{x^2 -\epsilon}{\sqrt{\mu+x}},\qquad t\in\R,
\end{gather}
which are well defined up to orientation. When considered for $x\in\CP^1$, the dynamics near the points $x=\infty$, resp.~$x=-\mu$ if $\mu^2\neq\epsilon$, is of hyperbolic type with~$3$, resp.~$1$, hyperbolic sectors, separated by~$3$, resp.~$1$, separatrices. For generic values of $\omega$, these separatrices all land at one of the singularities $x=\pm\sqrt{\epsilon}$, and divide a fixed neighborhood of the origin in~$\CP^1$ into
\begin{itemize}\itemsep=0pt
	\item an \emph{outer zone} attached to a single singularity (chosen as $x=\sqrt{\epsilon}$) and bounded by the separatrix of $x=\infty$ and a pair of separatrices of $x=-\mu$, all of them landing at $x=\sqrt{\epsilon}$, and
	\item an \emph{inner zone} attached to both singularities and bounded by the 3 separatrices of $x=-\mu$: if $\epsilon=0$ this zone splits into two zones, and if $\mu^2=\epsilon$ it disappears. It will be convenient to split the inner zone into a pair of \emph{inner half-zones} along a certain trajectory of \eqref{eq:RC-horizontal}, so that the description is uniform for both $\epsilon\neq 0$ and $\epsilon=0$ (see Fig.~\ref{figure:RC-zones1}).
\end{itemize}
The inner and outer zones vary continuously with $\epsilon$, $\mu$ and $\omega$ as long as the topology of the phase portrait does not bifurcate.

Similarly, one defines the \emph{outer} and \emph{inner zone} of the quadratic differential~\eqref{eq:RC-quaddiffomega}
 \emph{relative to a neighborhood} $\sX=\{x\in\C\,|\, |x|<\delta_x\}$ of the origin, as consisting of complete real trajectories of~\eqref{eq:RC-horizontal} inside $\sX$, i.e., trajectories that stay in $\sX$ for all $t\in\R$.

The domains $\rXI{}_\pm(\check\mu,\check\epsilon)$, resp.~$\rXO(\check\mu,\check\epsilon)$, are constructed as union of these inner half-zones, resp.\ outer zones, relative to $\sX$ over varying $\omega$, such that they are continuously varying with the parameters. More precisely, they are first defined as the zones for $\omega=0$ outside of the values of $(\mu,\epsilon)$ for which the horizontal foliation bifurcates, and then enlarged in the $\check x$-space and continued in the $(\check\mu,\check\epsilon)$-space by varying $\omega$ a bit. Moreover the singular points $x=\pm\sqrt{\epsilon}$,
resp.~$x=\sqrt{\epsilon}$, and the turning point $x=-\mu$ that are in the adherence of the zones are added to the domains.

The details of this construction are left to Section~\ref{sec:RC-2.22}, where it is performed in a new coordinate $s=\sqrt{\mu+x}$.

The construction of the canonical solution bases of Theorem~\ref{theorem:RC-IV} on the domains~$\rXO$ and~$\rXI{}_\pm$ is similar to~\cite{HLR}. It is based on the following two propositions.

\begin{Proposition}[subdominant solutions]\label{proposition:RC-IVa}
For any real half-trajectory of \eqref{eq:RC-horizontal}, that is a~trajectory $x(t)$ with $t\in\pm\R_{>0}$,
inside any of the domains $\rXI{}_\pm(\check\mu,\check\epsilon)$ or $\rXO(\check\mu,\check\epsilon)$, tending to one of the singularities,
there is a unique $1$-dimensional subspace of the solution space of~$\Delta(x,m)$ consisting of those solutions that are bounded along the half-trajectory and have vanishing limit at the singularity. This subspace, called the space of subdominant solutions, is uniquely defined and independent of the homotopy class of the half-trajectory with fixed end-point at the singularity within the domain.
\end{Proposition}

\begin{Proposition}[mixed bases]\label{proposition:RC-IVb} For any complete real trajectory of \eqref{eq:RC-horizontal} inside any of the domains $\rXI{}_\pm(\check\mu,\check\epsilon)$ or $\rXO(\check\mu,\check\epsilon)$, the two spaces of subdominant solutions associated to the positive $(t>0)$ and the negative~$(t<0)$ halves of the trajectory, are linearly independent. Therefore, by choosing a generator for each of the two subspaces one forms a~mixed basis for the whole solution space. The canonical solution bases of Theorem~{\rm \ref{theorem:RC-IV}} are of this form.
\end{Proposition}

The fact that the pair of subdominant solutions that form the mixed basis of Proposition~\ref{proposition:RC-IVa} is linearly independent is of utmost importance here,
and, in the case of the inner domains $\rXI{}_\pm(\check\mu,\check\epsilon)$ that vanish at the limit, its proof is far from trivial.

Related to the statement of Corollary~\ref{theorem:RC-III}, we have the following result on convergence of the normalizing gauge transformations of Remark~\ref{remark:RC-sectoralnormalization}(b) $T_{I\pm, m}$ and (c) $T_{i, m}$ to (a) $T_{O,m}$.

\begin{Theorem}[convergence of the local normalizing transformations]\label{theorem:RC-V}
Following the notation of Remark~{\rm \ref{remark:RC-sectoralnormalization}}:
\begin{enumerate}\itemsep=0pt
\item[$(i)$] For $\epsilon(m)\!=\!0$, the normalizing gauge transformations $T_{I+, m}$ $($resp.~$T_{I-, m})$ converge to~$T_{O,m}$, as $\mu(m)\to 0$ radially,
for each $m$ with $0<\arg\mu(m)<2\pi$ $($resp.~$0>\arg\mu(m)>-2\pi)$.
The convergence is uniform on compact sets in $S_{I+, m}$ $($resp.~$S_{I-, m})$.
\item[$(ii)$]
The normalizing gauge transformation $T_{2,m}$, analytic on a neighborhood of $x_2 = -\sqrt\epsilon$, and its analytic continuation,
converges to $T_{O, m}$, when $\mu(m)=O(\epsilon(m))$ and $\epsilon(m)\to 0$ radially
with $\arg\sqrt{\epsilon(m)}=\beta$ for $|\beta|<\frac{\pi}{2}$.
The convergence is uniform on compact subsets of the sector
\[\sX\setminus\{|\arg x-\beta|<\nu\},\qquad\text{with} \quad \nu>0\ \text{arbitrarily small.}\]
\end{enumerate}
\end{Theorem}

The first statement is closely related to a theorem on convergence of subdominant solutions by F.E.~Mullin \cite[Theorem II]{Mullin} (see also \cite[Chapter 3]{Sibuya}).
The second statement was originally established by A.~Glutsyuk \cite{Glu} in more general setting.

\section{Proofs}

Without loss of generality, we can always assume that the parametric system $\Delta(x,m)$ has the form \eqref{eq:RC-prenormal} with formal invariants \eqref{eq:RC-rfi}.

Our strategy will be the following. The part (b) of Theorem~\ref{theorem:RC-I} is a direct consequence of Theorem~\ref{theorem:RC-II}.
To prove part (a) of Theorem~\ref{theorem:RC-I}, we will first construct the canonical fundamental matrix solutions of Theorem~\ref{theorem:RC-IV} together with their natural domains $\rXO$, $\rXI{}_\pm$. Loosely speaking, the modulus of analytic equivalence can be identified with a certain conjugacy class of the set of the connection matrices (Stokes matrices) between the canonical fundamental solutions.
We will express these matrices explicitly, and show that in our situation the modulus reduces to the single analytic invariant $\gamma$.

We will prefer to do all this in a new ramified coordinate
\[s=\sqrt{\mu+x}.\]
The lifting to this $s$-coordinate produces a two-fold symmetry of the systems as well as of their normalizing gauge transformations.
After establishing the analytic equivalence of the lifted systems in the $s$-coordinate, one uses this symmetry to push it back down to the $x$-coordinate.

While everything, all the transformations and connection matrices, will depend on the parameter~$m$, we will often drop it from our notation,
and think of it as implicitly present; for example, we will often write $(\mu,\epsilon)$ rather than $(\mu(m),\epsilon(m))$.

\subsection[Systems in the $s$-coordinate]{Systems in the $\boldsymbol{s}$-coordinate}\label{sec:RC-2.3}
Let $\Delta(x,m)$ be a parametric system in the prenormal form \eqref{eq:RC-prenormal}.
We want to prove that two such systems with the same $\mu$, $\epsilon$ are analytically equivalent if and only if they have the same trace of monodromy~$\gamma$.

Let $s$ be a new coordinate defined by
\begin{gather}\label{eq:RC-xs}
x=s^2 -\mu,
\end{gather}
and let
\begin{gather}\label{eq:RC-SV}
S(s)=\left(\begin{smallmatrix} s^{-\frac12} & 0 \\ 0 & s^{\frac12} \end{smallmatrix}\right),\qquad
V=\tfrac{{\rm i}}{\sqrt 2}\left(\begin{smallmatrix} 1&1 \\ 1&-1\end{smallmatrix}\right).
\end{gather}
Then in the $s$-coordinate
\begin{gather*}
\Delta=\tfrac{(s^2 -\mu)^2 -\epsilon}{2s }\tfrac{{\rm d} }{{\rm d} s } +
\left[\left(\begin{smallmatrix} 0 & 1 \\ s^2 & 0 \end{smallmatrix}\right)
+\big(\big(s^2 -\mu\big)^2 -\epsilon\big)r \left(\begin{smallmatrix} 0 & 0 \\ 1 & 0 \end{smallmatrix}\right) \right],
\end{gather*}
and the transformed parametric system $ \Delta^{ s}:=\frac{1}{s}\cdot(SV)^*\Delta $ is equal to
\begin{gather}\label{eq:RC-d7a}
\Delta^{ s}(s,m)=\tfrac{(s^2 -\mu)^2 -\epsilon}{2s^2 }\tfrac{{\rm d} }{{\rm d} s} -
\left[\left(\begin{smallmatrix} 1 & 0 \\ 0 & -1 \end{smallmatrix}\right) +\tfrac{ (s^2 -\mu )^2 -\epsilon}{4s^3} \left(\begin{smallmatrix} 0 & 1 \\ 1 & 0 \end{smallmatrix}\right)
+\tfrac{ (s^2 -\mu )^2 -\epsilon}{2 s^2}r \left(\begin{smallmatrix} 1 & 1 \\ -1 & -1 \end{smallmatrix}\right) \right].
\end{gather}
The advantage of this new system is, as it will turn out, that it is diagonalizable on some domains $\Omega$ in the $(s,m)$-space:
We will be looking for analytic \emph{normalizing gauge transformations} $F_\Omega(s,m)$ on $\Omega$, bringing $\Delta^{ s}$
to a diagonal system $F_\Omega{}^*\Delta^{ s}=\widebar\Delta^{s}$, where
\begin{gather}\label{eq:RC-d7c}
\widebar\Delta^{s}(s,m)=\tfrac{(s^2 -\mu)^2 -\epsilon}{2s^2 \vphantom{\rX}}\tfrac{{\rm d} }{{\rm d} s \vphantom{\rX}}- \left(\begin{smallmatrix} 1 & 0 \\[3pt] 0 & -1 \end{smallmatrix}\right).
\end{gather}
This diagonal system $\widebar\Delta^{s}$ will serve as a \emph{model system} in the $s$-coordinate
for which one easily calculates a canonical fundamental matrix solution, denoted $\Psi$ \eqref{eq:RC-Psi}.
Then each lifted system $\Delta^{ s}(s,m)$ will have a canonical fundamental matrix solution $\Phi_\Omega=F_\Omega\Psi_\Omega$ on the domain $\Omega$,
where $\Psi_\Omega$ is a restriction of $\Psi$ \eqref{eq:RC-Psi} to $\Omega$, and therefore the original system $\Delta(x,m)$ will have a canonical fundamental solution matrix $Y_\Omega=SVF_\Omega\Psi_\Omega$ on the image of the domain $\Omega$ in the $x$-coordinate.

Also the quadratic differential \eqref{eq:RC-quadraticdifferential} becomes in the $s$-coordinate
\begin{gather*}
\left(\tfrac{2s^2}{(s^2 -\mu)^2 -\epsilon \vphantom{\rX}}{\rm d}s\right)^2,
\end{gather*}
which is the negative of the determinant of the meromorphic ``Higgs field'' $\frac{2s^2}{(s^2 -\mu)^2 -\epsilon}\left(\begin{smallmatrix} 1 & 0 \\ 0 & -1 \end{smallmatrix}\right){\rm d}s$ associated to~\eqref{eq:RC-d7c}.

The system $\widebar\Delta:=s\cdot\big(V^{-1}S^{-1}\big)^*\widebar\Delta^{s}$ in the $x$ coordinate, corresponding to $\widebar\Delta^{s}$,
\begin{gather*}
 \widebar\Delta(x,m)=\big(x^2 -\epsilon\big)\tfrac{{\rm d}}{{\rm d}x}-\left[\left(\begin{smallmatrix} 0 & 1 \\ \mu+x & 0 \end{smallmatrix}\right)-
 \tfrac{x^2 -\epsilon}{4(\mu + x)} \left(\begin{smallmatrix} 1 & 0 \\ 0 & -1 \end{smallmatrix}\right) \right],
\end{gather*}
has however an additional singularity at the point $x= -\mu$, hence does not belong to the considered class of systems. So instead, in the $x$-coordinate, one shall take the \emph{formal normal form} $\widehat\Delta(x,m)$~\eqref{eq:RC-model2} as the model. Now, if $E_\Omega(s,m)$ is the diagonalizing gauge transformation ``$F_\Omega(s,m)$'' for $\widehat\Delta^s=s^{-1}\cdot(SV)^*\widehat\Delta$ on the same domain~$\Omega$, then $\widehat\Delta(x,m)$ will have a canonical fundamental solution matrix $SVE_\Omega\Psi_\Omega$, and the composed transformation
\begin{gather}\label{eq:RC-T}
T_{\Omega}(x,m)=S(s)VF_\Omega(s,m) E_\Omega(s,m)^{-1}V^{-1}S(s)^{-1},
\end{gather}
defined on the ramified projection of the domain $\Omega$ into the $x$-coordinate,
will be non-singular at the point $x= -\mu$, $\mu^2\neq\epsilon$, and will bring $\Delta$ to ${T_{\Omega} }^*\Delta=\widehat\Delta$.
This is how one obtains the gauge transformations of Corollary~\ref{theorem:RC-III}.

The matrix functions
\begin{gather*} 
Y_\Omega=SVF_\Omega\Psi_\Omega,\qquad(\text{resp.}~SVE_\Omega\Psi_\Omega)
\end{gather*}
will be the canonical \emph{fundamental solutions} of Theorem~\ref{theorem:RC-IV} for the parametric systems $\Delta(x,m)$ (resp.~$\widehat\Delta(x,m)$).

\subsubsection[Fundamental solution of $\widebar\Delta^s(s,m)$]{Fundamental solution of $\boldsymbol{\widebar\Delta^s(s,m)}$}
On a neighborhood of $\infty$ on the Riemann sphere $\CP^1=\C\cup\{\infty\}$, define the function $\theta(s,\mu,\epsilon)$ by
\begin{gather*}
\frac{{\rm d}}{{\rm d}s}\theta(s,\mu,\epsilon)=\frac{2s^2}{(s^2 -\mu)^2 -\epsilon},\qquad \theta(\infty,\mu,\epsilon)=0.
\end{gather*}
We have
\begin{gather} \label{eq:RC-theta}
\hskip6pt \theta(s,\mu,\epsilon)
 =
 \begin{cases}
 \frac{\sqrt{\mu+\sqrt\epsilon}}{2\sqrt\epsilon}\log\frac{s-\sqrt{\mu+\sqrt\epsilon}}{s+\sqrt{\mu+\sqrt\epsilon}}-
	\frac{\sqrt{\mu-\sqrt\epsilon}}{2\sqrt\epsilon}\log\frac{s-\sqrt{\mu-\sqrt\epsilon}}{s+\sqrt{\mu-\sqrt\epsilon}}, &
	\text{if $\epsilon\big(\mu^2 - \epsilon\big)\neq 0$,}\vspace{1mm}\\
 -\frac{s}{s^2-\mu}-\frac{1}{2 \sqrt\mu}\log\frac{s+ \sqrt\mu}{s- \sqrt\mu}, & \text{if $\epsilon=0$,}\vspace{1mm}\\
 \frac{1}{\sqrt{2\mu}}\log\frac{s-\sqrt{2\mu}}{s+\sqrt{2\mu}}, & \text{if $\mu^2 =\epsilon$,}\vspace{1mm}\\
 -\frac{2}{s}, & \text{if $\mu,\epsilon=0$,}
 \end{cases}
\end{gather}
which is analytic in $ s\in\CP^1 \setminus \bigcup_{i=1}^4 [0,s_i]$,
if each $[0,s_i]$ denotes the closed segment between the origin and a zero point $s_i(\mu,\epsilon)$
of $ x^2(s)-\epsilon=\big(s^2 -\mu\big)^2 -\epsilon$.
The function $\theta(s,\mu,\epsilon)$ is continuous in $(\mu,\epsilon)\in\C^2$ and analytic for
$(\mu,\epsilon)\in\C^2\setminus\big\{\epsilon \big(\mu^2 - \epsilon\big)\neq 0\big\}$. It is odd in~$s$
\[\theta(-s,\mu,\epsilon)=-\theta(s,\mu,\epsilon),\]
and it satisfies
\[ \theta(s,\mu,\epsilon)=\theta\big(s,{\rm e}^{2\pi{\rm i}}\mu,\epsilon\big)=\theta\big(s,\mu,{\rm e}^{2\pi{\rm i}}\epsilon\big) \]
for each $s$ in its domain.

The matrix-valued function
\begin{gather}\label{eq:RC-Psi}
\Psi(s,\mu,\epsilon)=\left(\begin{smallmatrix} {\rm e}^{ \theta(s,\mu,\epsilon)} & 0 \\ 0 & {\rm e}^{-\theta(s,\mu,\epsilon)}
 \end{smallmatrix}\right)
\end{gather}
is a \emph{fundamental solution} for the diagonal model system $\widebar\Delta^s(s,m)$ \eqref{eq:RC-d7c}.

\subsubsection[$\Z_2$-symmetry]{$\boldsymbol{\Z_2}$-symmetry}
Let us remark that if $\Psi_\Omega(s,\mu,\epsilon)$ is a fundamental solution of $\widebar\Delta^{s}$ on a domain $\Omega(m)$ in the $s$-plane,
then so is $ \Psi_\Omega^\p (s,m):=\left(\begin{smallmatrix} 0 & 1 \\ 1 & 0 \end{smallmatrix}\right) \Psi_\Omega (-s,m) \left(\begin{smallmatrix} 0 & 1 \\ 1 & 0 \end{smallmatrix}\right)$,
this time on a rotated domain $\Omega(m)^\p:=-\Omega(m)$. The same is true for the system $\Delta^{ s}$.
Consequently, if $F_\Omega$ is a normalizing transformation for~$\Delta^{ s}$ on a domain $\Omega$, $ F_\Omega{}^*\Delta^{ s}=\widebar\Delta^{s}$, then so is $F_\Omega^\p$ on~$\Omega^\p$.
The following definition gives the~$()^\p$ notation precise meaning.

\begin{Definition}[rotation action of $\Z_2$]\label{notation:RC-p}
If $g(s)$ is a function on some domain $Y$ in the $s$-space, denote
\[g^\p(s):=g\big({\rm e}^{-\pi{\rm i}}s\big),\qquad s\in Y^\p:={\rm e}^{\pi{\rm i}}Y\]
the rotated function on the rotated domain. For a $(2 \times 2)$-matrix function $G(s)$, denote
\[G^\p(s):=\left(\begin{smallmatrix} 0 & 1 \\ 1 & 0 \end{smallmatrix}\right) G\big({\rm e}^{-\pi{\rm i}}s\big) \left(\begin{smallmatrix} 0 & 1 \\ 1 & 0 \end{smallmatrix}\right),\]
and for a constant matrix $C$,
\[C^\p:=\left(\begin{smallmatrix} 0 & 1 \\ 1 & 0 \end{smallmatrix}\right) C \left(\begin{smallmatrix} 0 & 1 \\ 1 & 0 \end{smallmatrix}\right).\]
\end{Definition}

\subsection[Domains $\Omega$ and diagonalizing transformations $F_\Omega$]{Domains $\boldsymbol{\Omega}$ and diagonalizing transformations $\boldsymbol{F_\Omega}$}\label{sec:RC-2.6}

A diagonalizing transformation $F_\Omega$ for the system $\Delta^s(s,m)$~\eqref{eq:RC-d7a} on a domain $\Omega$ gives rise to a~canonical fundamental solution matrix
of the system
\[\Phi_\Omega=F_\Omega\Psi_\Omega,\]
where $\Psi_\Omega$ is a branch of the fundamental solution $\Psi$ \eqref{eq:RC-Psi} of the diagonal model~\eqref{eq:RC-d7c} on the domain $\Omega$.
The shape of the domains $\Omega$ of such a bounded gauge transformation $F_\Omega$
is related to the rate of growth of $\Psi$ in the $s$-space, or equivalently
to the real dynamics of the (rotating family) of vector fields
\begin{gather}\label{eq:RC-chi}
{\rm e}^{{\rm i}\omega}\chi(s,\mu,\epsilon):={\rm e}^{{\rm i}\omega}\frac{\big(s^2 -\mu\big)^2 -\epsilon}{2 s^2} \frac{\partial}{\partial s}={\rm e}^{{\rm i}\omega}\frac{\partial}{\partial \theta},
\qquad \omega\in \big] {-} \tfrac{\pi}{2},\tfrac{\pi}{2}\big[ ,
\end{gather}
which will be studied in details in Section~\ref{sec:RC-vectorfield}.
By a \emph{real trajectory} of this vector field, we mean a solution curve $s(t)$ of the real ODE
\[\frac{{\rm d}s}{{\rm d}t}={\rm e}^{{\rm i}\omega}\frac{\big(s^2 -\mu\big)^2 -\epsilon}{2 s^2},\qquad t\in\R.\]
A positive, resp.\ negative, half-trajectory is one with $t>t_0$, resp.~$t<t_0$, for some $t_0\in\R$.
The real trajectory of ${\rm e}^{{\rm i}\omega}\chi$ through a point $s_0\in\C$ corresponds to the line $\theta(s_0,\mu,\epsilon)+{\rm e}^{{\rm i}\omega}\R$ in the
coordinate $\theta$ \eqref{eq:RC-theta}.

\subsubsection{Subdominant solutions}

\begin{Lemma}[subdominant solutions and bounded isotropies]\label{lemma:RC-boundedsolutions}
For a fixed parameter $m$, let~$s_i$ be a singular point, that is a zero of $\big(s^2 -\mu\big)^2 -\epsilon$, such that $s_i\neq 0$ unless $(\mu,\epsilon)=0$,
and let~$\sigma$ be a real positive $($resp.\ negative$)$ half-trajectory of
${\rm e}^{{\rm i}\omega}\chi$, for some $\omega\in \big] {-} \frac{\pi}{2},\frac{\pi}{2}\big[$,
with forward $t\to\infty$ $($resp.\ backward $t\to-\infty)$ limit at~$s_i$, hence such that
\[\lim_{\sigma\ni s\to s_i} {\rm e}^{-\theta}=0,\qquad\text{resp.} \ \
\lim_{\sigma\ni s\to s_i} {\rm e}^{\theta}=0.\]
\begin{enumerate}\itemsep=0pt
\item[$(a)$]
Then there is a unique solution $\phi(s)$ to $\Delta_m^s$ $($to the system \eqref{eq:RC-d7a} restricted to the fixed~$m)$ that is bounded along $\sigma$ when $s\to s_i$ and such that
\[\lim_{\sigma\ni s\to s_i} {\rm e}^{\theta}\phi(s)=\left(\begin{smallmatrix} 0 \\ 1 \end{smallmatrix}\right),\qquad\text{resp.}\quad
\lim_{\sigma\ni s\to s_i} {\rm e}^{-\theta}\phi(s)= \left(\begin{smallmatrix} 1 \\ 0 \end{smallmatrix}\right).\]
This solution does not depend on the half-trajectory as long as $\sigma$ and $\omega$ are varied continuously; moreover when also the point $s_i$ varies continuously
with the parameter $m$, then the solution $\phi$ depends continuously on $m$ and analytically for those $m$ for which $s_i$ is Fuchsian, i.e., a simple zero of $\big(s^2 -\mu\big)^2 -\epsilon$.

\item[$(b)$] Any isotropy $($automorphism$)$
$A(s)$ of the diagonal model system $\widebar\Delta^{s}_{m}$, that is a gauge transformation preserving the system
$A^*\widebar\Delta^{s}_{m}=\widebar\Delta^{s}_{m}$, which is
bounded along the half-trajectory $\sigma$ is of the form
\[A(s)=\Psi(s)C\Psi(s)^{-1},\qquad \Psi\ \text{as in \eqref{eq:RC-Psi}},\]
with an upper-triangular $($resp.\ lower-triangular$)$ constant invertible matrix $C=(c_{kl})$, and
\[\lim_{\sigma\ni s\to s_i}A(s)=\left(\begin{smallmatrix} c_{11} & 0 \\ 0 & c_{22}\end{smallmatrix}\right).\]
\end{enumerate}

In particular, an isotropy of $\widebar\Delta^{s}$ bounded along a complete $($both forward and backward$)$ real trajectory $\sigma$ is just a constant diagonal matrix.
\end{Lemma}

\begin{Remark}\quad
\begin{enumerate}\itemsep=0pt
\item The subspace of the solution space of $\Delta_m^s$ spanned by the solution $\phi(s)$ of Lemma~\ref{lemma:RC-boundedsolutions}(a) is
characterized as containing the solutions that are the most ``flat'' (or in this case bounded) along the given incoming trajectory to a singular point $s_{i_+}$ (resp.\ outgoing trajectory from~$s_{i_-}$). In the terminology of Sibuya~\cite{Sibuya} they are called \emph{subdominant solutions} along the half-trajectory.
It determines a flag structure on the solution space (cf.~\cite{HLR}).
\item When $s_i$ is a Fuchsian singularity then ${\rm e}^{\theta}\phi(s)$ (resp.~${\rm e}^{-\theta}\phi(s)$) is in fact analytic at~$s_i$, even if the singularity is resonant (the subdominant solutions don't have a logarithmic term).
\end{enumerate}
\end{Remark}

\begin{proof} $(a)$ The existence and uniqueness of such properly normalized subdominant solution follows for each fixed $m$ directly from the usual theorems on existence of local normalizing transformations $F_i$ at the point $s_i$, $F_i(s_i)=I$, similar to the situation discussed in
Remark~\ref{remark:RC-sectoralnormalization}, depending on the type of the singularity,
by taking the second (resp.\ first) column of thus constructed canonical fundamental solution matrix $F_i\Psi$ (so called ``Levelt basis'' in the Fuchsian case).
While the local analytic dependence on $m$ is quite clear as long as the type of the singularity does not change, the fact that it is also analytic at the values for which $s_i$ is a~resonant Fuchsian singularity, and that it passes well to the limit when $s_i$ becomes irregular, is perhaps less clear.
One obtains it from a consequence of a parametric version of the Levinson's theorem \cite[Theorem~5.3]{HLR} (see \cite[Theorem~8.1 in Chapter~3]{CL} for the original non-parametric version):
\begin{Theorem}[Levinson's theorem]
Consider a system of linear differential equations on the real line of the form
\[\frac{{\rm d}y}{{\rm d}t}=\big[\Lambda_0(m)+\Lambda_1(t,m)+P(t,m)\big]y,\]
where $\Lambda_0(m)$ is diagonal with distinct real parts of the eigenvalues,
$\Lambda_1(t,m)$ is also diagonal with limit zero at $t=+\infty$, and
\begin{gather}\label{eq:RC-Levinson1}
\int_{0}^{+\infty}\left|\frac{{\rm d}}{{\rm d}t}\Lambda_1(t,m)\right|{\rm d}t<\infty,\qquad \int_{0}^{+\infty}|P(t,m)|{\rm d}t<\infty.
\end{gather}
Then for each eigenvalue $\lambda(t,m)$ of $\Lambda_0(m)+\Lambda_1(t,m)$
there exists $t_0>0$ and a solution $\phi_\lambda(t,m)$ on $]t_0,+\infty[$ such that
\begin{gather*}
\lim_{t\to+\infty} \phi_\lambda(t,m)\cdot \exp\left(-\int_{t_0}^t \lambda(t,m){\rm d}t\right)=v_\lambda(m),
\end{gather*}
where $v_\lambda(m)\neq 0$ is a given eigenvector of $\Lambda_0(m)$ corresponding to the eigenvalue $\lambda(+\infty,m)$.
If the system depends continuously $($resp.\ analytically$)$ on a parameter $m$ over compact sets in the $t$-space, with the integrals in~\eqref{eq:RC-Levinson1}
uniformly bounded, then the solutions can be chosen depending continuously $($resp.\ analytically$)$ on~$m$.
\end{Theorem}

In order to apply the Levinson's theorem, we shall restrict our system to an annular domain
\begin{gather}\label{eq:RC-annular}
|s|<\delta_s,\qquad \left|\tfrac{(s^2-\mu)^2-\epsilon}{2s^4}\right|,\left|s\tfrac{\partial}{\partial s}\tfrac{(s^2-\mu)^2-\epsilon}{2s^4}\right|<K,
\end{gather}
where $\delta_s>0$ determines some neighborhood of $0$ and $K>0$ is arbitrary large, in particular large enough so that the point $s_i$ belongs to this domain.
We take
$t={\rm e}^{-{\rm i}\omega}(\theta(s,m)-\theta(s_0,m))$ (resp.\ $t=-{\rm e}^{-{\rm i}\omega}(\theta(s,m)-\theta(s_0,m))$), where $\theta$ is the rectifying coordinate \eqref{eq:RC-theta} for the vector field $\chi$ and $s_0$ is such that the half-trajectory $\sigma$ starting in~$s_0$ is contained in the above annular domain.

In the coordinate $t$, the vector field \eqref{eq:RC-chi} is ${\rm e}^{{\rm i}\omega}\chi=\frac{\partial}{\partial t}$, and if
$H:=I+\frac{\sqrt{1+(sf)^2}-1}{sf}\left(\begin{smallmatrix} 0 & -1\\ 1& 0 \end{smallmatrix}\right)$, with $f:=\frac{(s^2-\mu)^2-\epsilon}{4s^4}$,
is a matrix consisting of eigenvectors for $\left(\begin{smallmatrix}1 & 0\\ 0&-1\end{smallmatrix}\right)+sf\left(\begin{smallmatrix}0&1\\ 1&0\end{smallmatrix}\right)$,
then the transformed system $H^*(\Delta^s_m)$ becomes
\[\frac{{\rm d}}{{\rm d}t}-\left[\Lambda_0(m)+P(t,m)\right],\]
with $\Lambda_0={\rm e}^{{\rm i}\omega}\left(\begin{smallmatrix} 1 & 0 \\ 0 & -1 \end{smallmatrix}\right)$, and
\begin{gather*}
{\rm e}^{-{\rm i}\omega}P(t,m)=\big(\sqrt{1+(sf)^2}-1\big)\left(\begin{smallmatrix} 1 & 0 \\ 0 & -1 \end{smallmatrix}\right)+
\tfrac{(s^2-\mu)^2-\epsilon}{2s^2}rH^{-1}\left(\begin{smallmatrix} 1 & 1 \\ -1 & -1 \end{smallmatrix}\right)H
-\tfrac{(s^2-\mu)^2-\epsilon}{2s^2}H^{-1}\tfrac{{\rm d}H}{{\rm d}s}.
\end{gather*}
To verify the condition \eqref{eq:RC-Levinson1} it is enough to show that $\frac{2s^2}{(s^2-\mu)^2-\epsilon}P(t,m)$ is uniformly bounded on the annular domain~\eqref{eq:RC-annular}. This follow from the boundedness of $\frac{2s^2}{(s^2-\mu)^2-\epsilon}(sf)^2=\frac{(s^2-\mu)^2-\epsilon}{2s^4}=2f$ and of
$\frac{\partial (2sf)}{\partial s}=2f+s\frac{\partial (2f)}{\partial s}$ on~\eqref{eq:RC-annular}.

In order to cover the whole parameter space, including the situation when a singularity $s_i$ approaches $0$ as $\mu^2\to\epsilon^2$, one needs to apply the parametric Levinson's theorem for increasing~$K\to+\infty$.

$(b)$ Let $\theta(s)$ be a branch of the function $\theta(s,\mu,\epsilon)$ in \eqref{eq:RC-theta} on~$U$.
We have
\[A(s)=\begin{pmatrix} c_{11} & {\rm e}^{2\theta(s)}c_{12} \\ {\rm e}^{-2\theta(s)}c_{21} & c_{22} \end{pmatrix}.\]
If $|\omega|<\frac{\pi}{2}$, then $\Re(\theta(s))\to+\infty$ (resp.~$-\infty$) as $\sigma\ni s\to s_i$,
which implies that $c_{12}=0$ (resp.~$c_{21}=0$), otherwise~$A$ would not be bounded.
\end{proof}

\subsubsection[Construction of the domains $\Omega$ and of the fundamental solution matrices on them]{Construction of the domains $\boldsymbol{\Omega}$ \\ and of the fundamental solution matrices on them}\label{paragraph:RC-Omega}
Let
\begin{gather}\label{eq:RC-sEsM}
\sS=\{0<|s|<\delta_s\},\qquad \cM=\{ |\mu|<\delta_\mu\},\qquad \cE=\{|\epsilon|<\delta_\epsilon\},
\end{gather}
be small discs, $\delta_s,\delta_\mu,\delta_\epsilon>0$. And let $s_i(\check\mu,\check\epsilon)$, $i=1,\dots,4$, be the zeros of $ \big(s^2 -\mu\big)^2 -\epsilon $
depending continuously on a ramified coordinate $(\check\mu,\check\epsilon)$ from the universal covering space of $\cM \times \cE$
ramified over the set $\big\{\epsilon \big(\mu^2 -\epsilon\big)=0\big\}$, where the ramification set is considered to be included in the covering space, and the topology is the lifted preimage of that on $\cM \times \cE$. We shall suppose that~$\delta_\mu$,~$\delta_\epsilon$ are small enough so that all the zero points $s_i(\check\mu,\check\epsilon)$ fall inside the disc of radius $\delta_s$,
\[\delta_\mu+\delta_\epsilon^{\frac{1}{2}}\ll \delta_s^2.\]
For each $(\check\mu,\check\epsilon)$, the disc $\sS$ has a universal ramified cover $\rS(\check\mu,\check\epsilon)$ with ramification at the zero points $s_i(\check\mu,\check\epsilon)$ of $\big(s^2 -\mu\big)^2 -\epsilon$. They glue up together
to form a ramified covering $\coprod_{(\check\mu,\check\epsilon)}\rS(\check\mu,\check\epsilon)$ of the $(s,\mu,\epsilon)$-space $\sS\times\cM\times\cE$.

\begin{Definition}\label{definition:RC-zones}	Let us consider the real phase portrait of the vector field ${\rm e}^{{\rm i}\omega}\chi$ \eqref{eq:RC-chi} inside a~pointed disc $\sS^*:=\sS\setminus\{0\}$.
	\begin{itemize}\itemsep=0pt
		\item[--] An angle $\omega\in \big] {-} \frac{\pi}{2},\frac{\pi}{2}\big[$ is \emph{admissible} (for given $(\mu,\epsilon)$) if for all $s_0\in\sS^*$ the real trajectory $s(t)$ of ${\rm e}^{{\rm i}\omega}\chi$ through $s_0=s(0)$ stays in $\sS^*$ for either all positive ($t>0$) or all negative ($t<0$) time.		
		\item[--] The \emph{zones of ${\rm e}^{{\rm i}\omega}\chi$ in $\sS^*$} are the connected components of the complement in $\sS^*\setminus\big\{\big(s^2-\mu\big)^2-\epsilon=0\big\}$ of all the real trajectories of ${\rm e}^{{\rm i}\omega}\chi$ that leave $\sS^*$.
\end{itemize}
\end{Definition}
By considering the vector field ${\rm e}^{{\rm i}\omega}\chi$ in the punctured disc $\sS^*$ the definition the zones stays the same whether $s=0$ is its pole ($\mu^2\neq\epsilon$) or not ($\mu^2=\epsilon$).
For an admissible $\omega$, there are no periodic orbits, and each zone is spanned by complete real trajectories of ${\rm e}^{{\rm i}\omega}\chi$ in $\sS^*\setminus\big\{\big(s^2-\mu\big)^2-\epsilon=0\big\}$ that have the same pair of forward/backward limit points at the singular points, and are ``homotopic'' to each other.
And they evolve continuously with $\omega$ as long as the $\omega$ is admissible. More details will be given in Sections~\ref{sec:RC-vectorfield} and~\ref{sec:RC-2.22} below.

For a fixed $m$, one associates to each complete real trajectory $\sigma$ of the vector field ${\rm e}^{{\rm i}\omega}\chi$ within $\sS^*$, which starts and terminates in two equilibrium points $s_{i+}$ and $s_{i-}$, a fundamental solution matrix $\Phi_\sigma$ of~\eqref{eq:RC-d7a} whose first column is the unique solution provided by Lemma~\ref{lemma:RC-boundedsolutions}(a) with the given asymptotics along the backward orbit,
and whose second column is the right asymptotics along the forward orbit.\footnote{The columns of such fundamental matrix $F_\Omega\Psi$ form a so called \emph{mixed basis} of the solution space, originally introduced by J.-P.~Ramis~\cite{Ra} and C.~Zhang~\cite{Zha}. See also \cite{HLR, LR}.}
Clearly, this fundamental solution matrix is independent of the trajectory $\sigma$ within the same zone of ${\rm e}^{{\rm i}\omega}\chi$ in $\sS^*$.
Moreover, it is also independent of the angle $\omega\in \big] {-} \frac{\pi}{2},\frac{\pi}{2}\big[$ as long as it is admissible.
Therefore we will construct domains $\Omega(\check\mu,\check\epsilon)$ as ramified unions of the topological closures of the zones of ${\rm e}^{{\rm i}\omega}\chi$ in $\sS^*$ over admissible $\omega$, and define the fundamental solutions $\Phi_\Omega$ on $\Omega(\check\mu,\check\epsilon)$ as the solution $\Phi_\sigma$ for any trajectory $\sigma$ inside the domain. We will then trace the evolution of these domains in dependence on $(\check\mu,\check\epsilon)$,
and define
\[\Omega=\coprod_{\check m}\Omega(\check\mu(\check m),\check\epsilon(\check m)),\]
as their ramified union in the $(\check s,\check m)$-space. We will describe the domains $\Omega$ obtained this way in Section~\ref{sec:RC-2.22}.

To fix the notation, from now on let
\begin{gather}\label{eq:RC-si}
s_1(\check\mu,\check\epsilon):=\sqrt{\check\mu+\sqrt{\check\epsilon}}, \qquad s_2(\check\mu,\check\epsilon):=\sqrt{\check\mu-\sqrt{\check\epsilon}}
\end{gather}
such that for $\arg\epsilon=\arg\mu=0$ and $\mu>\sqrt\epsilon>0$ they are given by the usual square root.

\begin{Proposition}[diagonalizing gauge transformations]\label{proposition:RC-3}
Let a parametric system $\Delta^{ s}(s,m)$ be as in \eqref{eq:RC-d7a} and its diagonal model $\widebar\Delta^s(s,m)$ be as in~\eqref{eq:RC-d7c}.
There are $4$ different parametric domains~$\Omega$ as defined above: a~symmetric pair of inner domains~$\Omega_I$ adjoint to the points $\{s_2,s_1\}$, $\Omega_I^\p$ adjoint to the points $\{-s_1,-s_2\}$, and a symmetric pair of outer domains~$\Omega_O$, $\Omega_O^\p$ both adjoint to the points $\{-s_1,s_1\}$. On these domains there exist unique diagonalizing gauge transformations
\begin{alignat*}{3}
& F_\bullet\in\GL_2(\Cal B(\Omega_\bullet)), \qquad && (F_\bullet)^*\Delta^{ s}=\widebar\Delta^s, & \\
& F_\bullet^\p\in\GL_2\big(\Cal B\big(\Omega_\bullet^\p \big)\big),\qquad && \big(F_\bullet^\p \big)^*\Delta^{ s}=\widebar\Delta^s,\qquad \bullet=O,I,
\end{alignat*}
$($see Definition~{\rm \ref{definition:RC-CalB}} and Notation~{\rm \ref{notation:RC-p})},
such that
\begin{alignat}{3}
& F_I(s_1,\check m)=\left(\begin{smallmatrix} 1 & 0 \\ 0 & \kappa_I(\check m) \end{smallmatrix}\right) ,\qquad &&
 F_I(s_2,\check m)=\left(\begin{smallmatrix} \kappa_I(\check m) & 0 \\ 0 & 1 \end{smallmatrix}\right),& \nonumber\\
& F_O(s_1,\check m)=\left(\begin{smallmatrix} 1 & 0 \\ 0 & \kappa_O(\check m) \end{smallmatrix}\right) ,\qquad &&
 F_O(-s_1,\check m)=\left(\begin{smallmatrix} \kappa_O(\check m) & 0 \\ 0 & 1 \end{smallmatrix}\right), & \label{eq:RC-kappaIO}
\end{alignat}
$\check m\in\rM$,
where the functions $ \kappa_\bullet=\det F_\bullet \in\Cal B\big(\rM\big)$, $\bullet=O,I$, are uniquely determined by $\Delta^{ s}$,
$\kappa_O(\check m)=1$ if $\mu(m)=\epsilon(m)=0$, and $\kappa_I(\check m)=1$ if $\epsilon(m)=0$.
\end{Proposition}

\begin{proof}The domains $\Omega_\bullet$ will be constructed in detail in Section~\ref{sec:RC-2.22}.
The first and the second columns of $\Phi_\bullet=F_\bullet\Psi_\bullet$ are respectively the unique vectors of Lemma~\ref{lemma:RC-boundedsolutions}(a) at the singular points $\{s_{i_-},s_{i_+}\}$ (in the same order as in the statement).
By their construction they are each bounded at their respective point. But they are also bounded at the other point.
In fact, remark that if $\phi$ is \emph{any} solution of~$\Delta_m^s$ then
$-{\rm i}\sqrt 2s^{\frac{1}{2}}{\rm e}^{-\theta}\phi$ is bounded along any incoming trajectory to~$s_{i_+}$
(resp.~$-{\rm i}\sqrt 2s^{\frac{1}{2}}{\rm e}^{\theta}\phi$ is bounded along any outgoing trajectory from~$s_{i_-}$).

The fact that $\kappa_\bullet=\det F_\bullet$ is constant in $s$ follows from
the Liouville--Ostrogradski formula:
\[\kappa_\bullet=\det F_\bullet=\det(SVF_\bullet\Psi_\bullet)\]
is constant since the trace of the matrix of the system \eqref{eq:RC-prenormal} is null.

We still need to prove that the gauge transformations are invertible for small $m$, i.e., that $\kappa_\bullet(\check m)\neq 0$.
For the outer domains and $\kappa_O$, this follows from the continuity of the construction in $\check m$ which persist well to the limit $m\to 0$, and the fact that $F_O(0,0)=I$.
For the inner domains, their construction gives us that $\kappa_I(\check m)=1$ for $\epsilon(m)=0$, $\mu(m)\neq 0$, but doesn't tell us if a limit $\lim\limits_{m\to 0}\kappa_I(\check m)$ exists. We'll prove it in Corollary~\ref{corollary:RC-kappaI}, until then we'll treat $\kappa_I$ as an analytic function of $\check m$ which may a priori have zeros at some points.
\end{proof}

\subsubsection[The vector field $\chi$]{The vector field $\boldsymbol{\chi}$}\label{sec:RC-vectorfield}
In this section we will study the real phase portrait of the vector fields ${\rm e}^{{\rm i}\omega}\chi$ \eqref{eq:RC-chi} in $\C^*=\C\setminus\{0\}$. And in the following Section~\ref{sec:RC-2.22} we will describe the effect of its restriction to $\sS^*$, and construct the domains $\Omega$.

\begin{Remark}[rotated vector field]\label{remark:RC-omega}
The change of coordinates
\[\ (s,\mu,\epsilon)\mapsto \big({\rm e}^{{\rm i}\omega} s,{\rm e}^{2{\rm i}\omega}\mu,{\rm e}^{4{\rm i}\omega}\epsilon\big),\qquad \omega\in\C,\]
transforms the vector field $\chi$ to ${\rm e}^{{\rm i}\omega}\chi$. This means we can restrict the discussion to $\omega=1$.
\end{Remark}

\begin{figure}[t]\centering
 \begin{tikzpicture}[scale=0.99]
\node at (0,0){\includegraphics{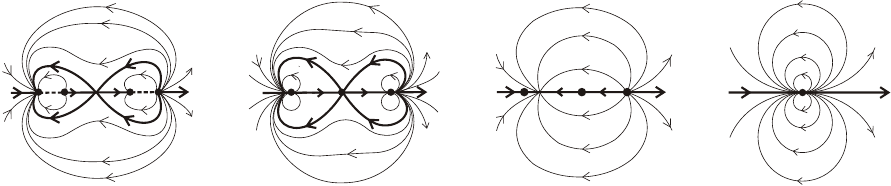}};
\node at (-7.6,0){\footnotesize $\infty$};
\node at (-4.2,0){\footnotesize $\infty$};
\node at (-3.3,0){\footnotesize $\infty$};
\node at (-0.1,0){\footnotesize $\infty$};
\node at (0.65,0){\footnotesize $\infty$};
\node at (4.0,0){\footnotesize $\infty$};
\node at (4.65,0){\footnotesize $\infty$};
\node at (7.8,0){\footnotesize $\infty$};
\node at (-5.95,-0.24){\footnotesize $0$};
\node at (-1.75,-0.24){\footnotesize $0$};
\node at (2.35,-0.22){\footnotesize $0$};
\node at (6.1,-0.20){\footnotesize $0$};
\node at (-6.0,-1.9){(a) $\mu>\sqrt{\epsilon}>0$};
\node at (-1.75,-1.9){(b) $\mu>\sqrt\epsilon=0$};
\node at (2.25,-1.9){(c) $\mu=\sqrt{\epsilon}>0$};
\node at (6.1,-1.9){(d) $\mu=\sqrt\epsilon=0$};
 \end{tikzpicture}
\caption{Examples of the real trajectories of the vector field $\chi$ \eqref{eq:RC-chi}, with the separatrices of $0$ and $\infty$ in bold, for selected values of $(\mu,\epsilon)$. See Fig.~\ref{figure:RC-2} below for other values. In (a) the dotted bold trajectories between $s_1=\sqrt{\mu+\sqrt\epsilon}$ and $s_2=\sqrt{\mu-\sqrt\epsilon}$, and $-s_2$ and $-s_1$, split each of the inner zones into 2 half-zones corresponding to those in (b).}	\label{figure:RC-chi}
\end{figure}

The vector field
\[\chi(s,\mu,\epsilon)=\frac{\big(s^2 -\mu\big)^2 -\epsilon}{2s^2}\frac{\partial}{\partial s}=\frac{1}{2}\big(\epsilon s^{-4}-\big(\mu s^{-2} - 1\big)^2\big)\frac{\partial}{\partial s^{-1}}\]
is a rational vector field on~$\C$, but becomes a polynomial vector field in the coordinate~$s^{-1}$ on
$\CP^1=\C\cup\{\infty\}$ with a regular point at $s=\infty$.
The real dynamics of complex polynomial vector fields on $\CP^1$ has been extensively studied in \cite{Ben,BD,DES} (see also \cite{KR,Mucino,Tahar,Tomasini}).
Some of the basic properties when applied to ${\rm e}^{{\rm i}\omega}\chi$ can be summarized as:
\begin{itshape}
\begin{itemize}\itemsep=0pt
\item See Fig.~{\rm \ref{figure:RC-chi}}.
\begin{itemize}\itemsep=0pt
	\item[--] For $\epsilon\neq\mu^2$ the vector field is of degree $4$ in $s^{-1}$, and the point $s=0$ is hyperbolic with~$6$ local separatrices $($alternating incoming/outgoing$)$, which can be either homoclinic or they terminate at an equilibrium point $($zero of $\big(s^2 -\epsilon\big)^2 -\mu)$.
	The~$6$ sectors at~$0$ in between of the separatrices are called ends.
	\item[--] For $\epsilon=\mu^2\neq 0$ the vector field is of degree $2$ in $s^{-1}$, and the point $s=0$ is regular, but shall be thought of as hyperbolic with $2$ local separatrices $($one incoming and one outgoing$)$ and with~$2$ ends.
	\item[--] For $\epsilon=\mu^2=0$ the vector field is of degree $2$ in $s^{-1}$, the point $s=0$ is a double equilibrium, and the point $s=\infty$ which is regular shall be thought of as hyperbolic with $2$ local separatrices $($one incoming and one outgoing$)$ and with $2$ ends.
\end{itemize}

\item The connected components of the complement in $\CP^1$ of all the separatrices of the hyperbolic point $s=0$, and of all the equilibria~$s_i$,
are called zones of ${\rm e}^{{\rm i}\omega}\chi$ in $\CP^1$.
The zones can be of~$3$ types
\begin{itemize}\itemsep=0pt
	\item[--] Center zone consisting of periodic trajectories around an equilibrium center. The image of a center zone in the coordinate $\theta$ is an infinite half-strip 	perpendicular to the line ${\rm e}^{{\rm i}\omega}\R$, whose two bounding rays are identified by a period shift.
	\item[--] $\alpha\omega$-zone consisting of trajectories that converge in forward, resp.\ backward time, to the same equilibrium, and these two equilibria are distinct. The image of an $\alpha\omega$-zone in the coordinate $\theta$ is an open infinite strip of a finite width, parallel with the line~${\rm e}^{{\rm i}\omega}\R$.
	\item[--] Sepal zone consisting of trajectories that converge in both forward and backward time to the same one equilibrium, which is necessary multiple.
	 The image of a sepal zone in the coordinate $\theta$ is an open half-plane with boundary parallel with the line ${\rm e}^{{\rm i}\omega}\R$.
\end{itemize}
In particular, there are no limit cycles.
\item The vector field ${\rm e}^{{\rm i}\omega}\chi$ is called rotationally stable if there are no homoclinic separatrices. In particular, there are no centers since the boundary of a center zone is formed by a~union of homoclinic separatrices.
This is equivalent to the real phase portrait of ${\rm e}^{{\rm i}\omega}\chi$ staying topologically equivalent under a small change in~$\omega$.
For each fixed $(\mu,\epsilon)$, the set of $\omega\in \big] {-} \frac{\pi}{2},\frac{\pi}{2}\big[$ for which ${\rm e}^{{\rm i}\omega}\chi$ is not rotationally stable is finite, in fact, in our situation there are no more then~$4$ such values.
\item When all the equilibria~$s_i$ are simple and $\neq 0$ and the vector field is rotationally stable, then there are exactly $3$ zones of ${\rm e}^{{\rm i}\omega}\chi$ in $\CP^1$ $($as the vector filed is of order~$4$ in $s^{-1})$, which are all of $\alpha\omega$ type, and each has two ends at~$0$.
For $\epsilon=\mu^2\neq 0$, then $\chi$ has two simple equilibria at $s=\pm\frac{\mu}{2}$ and $0$ is a regular point, and the whole $\CP^1\setminus\big\{{-}\frac{\mu}{2},\frac{\mu}{2}\big\}$ is a single $\alpha\omega$-zone. And for $\epsilon=\mu=0$, there are~$2$ sepal zones attached to the double equilibrium at~$0$.
\end{itemize}\end{itshape}

When considered in $\C^*=\CP^1\setminus\{0,\infty\}$, one shall consider also the point $s=\infty$ as a hyperbolic point with 2 separatrices, and for $\mu^2=\epsilon\neq 0$ the same also with $s=0$. We then talk about \emph{zones of ${\rm e}^{{\rm i}\omega}\chi$ in $\C^*$}.

If ${\rm e}^{{\rm i}\omega}\chi$ is rotationally stable as a polynomial vector field on $\CP^1$, then it cannot have any heteroclinic connection between $0$ and $\infty$ as the symmetry of ${\rm e}^{{\rm i}\omega}\chi$ would force it to have a pair of them that would form a homoclinic connection of $0$ passing through $\infty$, contradicting the assumption, and neither can $\infty$ have a homoclinic connection. Thus it is rotationally stable also in $\C^*$.

\begin{Definition}[half-zones]
Assume ${\rm e}^{{\rm i}\omega}\chi$ is rotationally stable.
\begin{itemize}\itemsep=0pt
	\item[--] For each $\alpha\omega$-zone of ${\rm e}^{{\rm i}\omega}\chi$ in $\CP^1$, there is a unique trajectory whose image by $\theta$ is the line that splits the strip, that is the image of the zone by $\theta$, lengthwise into two strips of equal widths. This trajectory
	splits the $\alpha\omega$-zone into 2 \emph{half-zones}.

	\item[--] Each sepal zone of ${\rm e}^{{\rm i}\omega}\chi$ in $\CP^1$ is already a \emph{half-zone}.
\end{itemize}
For $(\mu,\epsilon)\neq0$, each half-zone contains exactly one end at the hyperbolic point $s=0$.

For the reason of symmetry, for $(\mu,\epsilon)\neq0$, the trajectory through $s=\infty$ is always the splitting trajectory for the $\alpha\omega$-zone in $\CP^1$ containing $\infty$,
dividing it into a pair of \emph{outer half-zones}.
\end{Definition}

\begin{Proposition}For $\omega$ such that ${\rm e}^{{\rm i}\omega}\chi$ is rotationally stable, and $(\mu,\epsilon)\in\cM\times\cE$, $\mu^2\neq\epsilon$, $\epsilon\neq 0$, the vector field has exactly~$4$ different zones in $\C^*$, and $6$ different half-zones each having exactly~$1$ end at $s=0$. See Fig.~{\rm \ref{figure:RC-chi}}.
There is:
\begin{itemize}\itemsep=0pt
	\item a symmetric pair of outer half-zones bounded by the trajectory through $\infty$ and the separatrices of the origin,
	\item a symmetric pair of inner zones that are bounded solely by the separatrices of the origin, each of these inner zones is divided into $2$
	inner half-zones.
\end{itemize}
At the limit, when $\epsilon= 0$, $\mu\neq0$, all the half-zones persist, and when $\mu^2=\epsilon$ the $4$ inner half-zones become empty and only the $2$ outer half-zones persist.
\end{Proposition}

{\bf Bifurcation diagram of $\boldsymbol{\chi}$.}
Let us take a better look on how these half-zones evolve depending on the parameters $(\mu,\epsilon)$ and on $\omega$ (cf.\ Remark~\ref{remark:RC-omega}).
A bifurcation of the real phase portrait of the vector field ${\rm e}^{{\rm i}\omega}\chi$ can occur if it either becomes rotationally unstable, or when either an equilibrium or a hyperbolic point changes its multiplicity. This second bifurcation, occurring when $\epsilon\big(\mu^2-\epsilon\big)=0$, doesn't affect the decomposition of the real phase portrait into half-zones other than that some of the half-zone may become empty.
There are two possibilities how rotational instability of ${\rm e}^{{\rm i}\omega}\chi$ can occur: 1) through an appearance of a homoclinic separatrix of the origin in~$\CP^1$, either encircling a single singularity (let's denote this case $\Sigma_I$) or a pair of singularities (let's denote this case~$\Sigma_O$),
or~2) through an appearance of a heteroclinic separatrix connecting~$0$ and $\infty$: due to the symmetry of ${\rm e}^{{\rm i}\omega}\chi$ this latter bifurcation agrees exactly with the case~$\Sigma_O$.

\begin{figure}[b!]\centering
 \begin{tikzpicture}[scale=1]
\node at (0,0){\includegraphics[width=0.84\textwidth]{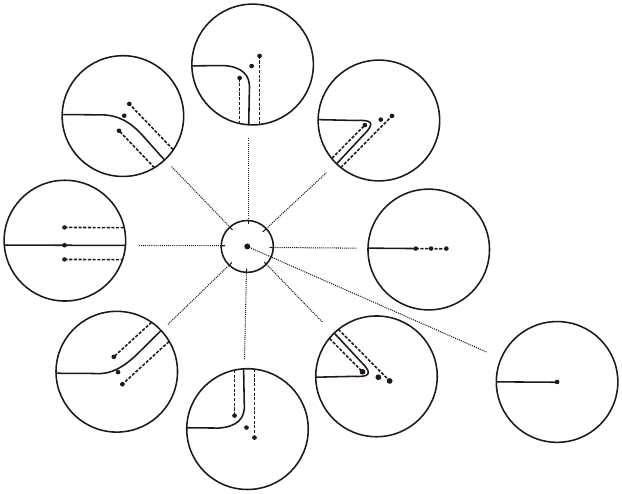}};
\node at (-1.55,0){$0$};
\node at (-0.7,0.2){$\epsilon$};
\node at (1.9,-0.3){$-\sqrt{\epsilon}$};
\node at (2.9,-0.3){$\sqrt{\epsilon}$};
\node at (4.1,0){(i)};
\node at (-1.4,-5.55){(iv)};
\node at (-1.3,5.5){(ii)};
\node at (-7.0,0){(iii)};
\node at (7.0,-2.9){(0)};
\node at (5.3,-3.2){$0$};
 \end{tikzpicture}
	\caption{Bifurcation curves in the $\mu$-plane for the vector field $\chi(s,\mu,\epsilon)$ (i.e., $\omega=0$) according to values of~$\epsilon$: dashed lines $\Sigma_{I,0}(\epsilon)$ correspond to change of stability of a singular point, solid line curve $\Sigma_{O,0}(\epsilon)$ corresponds to bifurcation of the trajectory passing through $\infty$.} 	\label{figure:RC-1}
\end{figure}

\begin{itemize}\itemsep=0pt 
\item[$\Sigma_I$:]
The bifurcation $\Sigma_I$ occurs when the stability of a zero point $s_i=\pm\sqrt{\mu\pm\sqrt\epsilon}$ of ${\rm e}^{{\rm i}\omega}\chi$ changes between attractive and repulsive, i.e., when the multiplier $\pm {\rm e}^{{\rm i}\omega}\frac{2\sqrt\epsilon}{s_i}$ of the
linearization $\pm {\rm e}^{{\rm i}\omega}\frac{2\sqrt\epsilon}{s_i}(s-s_i)\partial_{s}$ of vector field ${\rm e}^{{\rm i}\omega}\chi$ at the point $s_i$ becomes purely imaginary:
${\rm e}^{{\rm i}\omega}\frac{2\sqrt\epsilon}{\sqrt{\mu\pm\sqrt\epsilon}}\in {\rm i}\R$, which is equivalent to
\begin{gather*}
\mu\in\mp\sqrt\epsilon-{\rm e}^{2{\rm i}\omega}\epsilon \R_{\geq 0}=:\Sigma_{I,\omega}(\epsilon).
\end{gather*}
It is well known that a holomorphic vector field in $\C$ is analytically equivalent to its linearization near each simple zero
(see, e.g., \cite[Theorem~5.5]{IlYa}).
As a consequence, if $\mu\in\Sigma_{I,\omega}(\epsilon)$ (the dashed lines in Fig.~\ref{figure:RC-1}) then the real phase portrait of ${\rm e}^{{\rm i}\omega}\chi$ near the point~$s_i$ with purely imaginary multiplier is that of a~center. By Remark~\ref{remark:RC-omega}, $\mu\in\Sigma_{I,\omega}(\epsilon)\Longleftrightarrow {\rm e}^{2{\rm i}\omega}\mu\in\Sigma_{I,0}({\rm e}^{4{\rm i}\omega}\epsilon)$.

\item[$\Sigma_O$:]
The bifurcation $\Sigma_O$ occurs when the trajectory through infinity passes by the origin.
This means that $\theta(0,\mu,\epsilon)-\theta(\infty,\mu,\epsilon)\in {\rm e}^{{\rm i}\omega}\R$, where $\theta$ is as in \eqref{eq:RC-theta},
i.e., $ \tfrac{\sqrt{\mu+\sqrt\epsilon}\pm\sqrt{\mu-\sqrt\epsilon}}{2\sqrt\epsilon}\pi{\rm i} \in {\rm e}^{{\rm i}\omega}\R,$
which is equivalent to $ -\frac{\mu\pm\sqrt{\mu^2-\epsilon}}{\epsilon}=a\in {\rm e}^{2{\rm i}\omega}\R_{>0} $, that is
\begin{gather*} 
\mu\in\big\{ {-}\tfrac{1}{2}\big(a^{-1} +\epsilon a\big)\,|\, a\in {\rm e}^{2{\rm i}\omega}\R_{>0}\big\}=:\Sigma_{O,\omega}(\epsilon).
\end{gather*}
The set $\Sigma_{O,\omega}(\epsilon)$ is a branch of a hyperbola (the solid curve in Fig.~\ref{figure:RC-1}).
By Remark~\ref{remark:RC-omega}, $\mu\in\Sigma_{O,\omega}(\epsilon)\Longleftrightarrow {\rm e}^{2{\rm i}\omega}\mu\in\Sigma_{O,0}({\rm e}^{4{\rm i}\omega}\epsilon)$.
\end{itemize}

\begin{figure}[t]\centering
\leftline{(i) $\epsilon\in\R_{>0}$:}
 \begin{tikzpicture}[scale=1]
\node at (0,0){\includegraphics{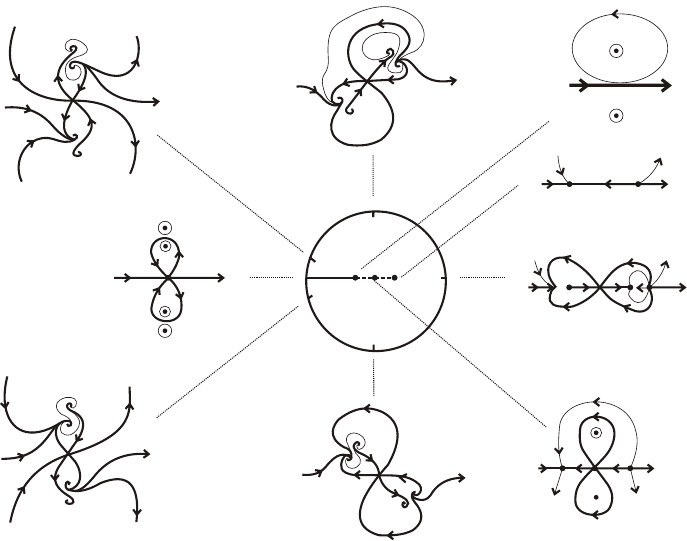}};
\node at (0,0.5){\footnotesize $\mu$};
\node at (-5.9,2.75){\footnotesize $\infty$};
\node at (-2.9,2.85){\footnotesize $\infty$};
\node at (-1.0,3.1){\footnotesize $\infty$};
\node at (2.1,3.2){\footnotesize $\infty$};
\node at (3.6,3.1){\footnotesize $\infty$};
\node at (5.8,3.1){\footnotesize $\infty$};
\node at (3.2,1.45){\footnotesize $\infty$};
\node at (5.7,1.45){\footnotesize $\infty$};
\node at (2.9,-0.3){\footnotesize $\infty$};
\node at (6.0,-0.3){\footnotesize $\infty$};
\node at (-1.85,-0.15){\footnotesize $\infty$};
\node at (-4.1,-0.15){\footnotesize $\infty$};
\node at (-3.05,-3.1){\footnotesize $\infty$};
\node at (-6.0,-3.2){\footnotesize $\infty$};
\node at (2.25,-3.5){\footnotesize $\infty$};
\node at (-0.9,-3.45){\footnotesize $\infty$};
\node at (5.5,-3.35){\footnotesize $\infty$};
\node at (3.05,-3.35){\footnotesize $\infty$};
 \end{tikzpicture}
\leftline{(ii) $\epsilon\in -{\rm i}\R_{>0}$:}
 \begin{tikzpicture}[scale=1]
\node at (0,0){\includegraphics{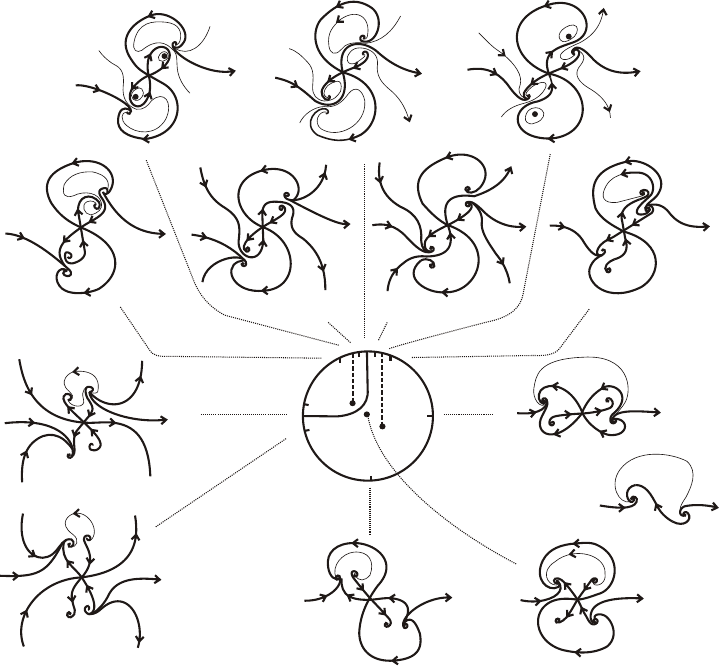}};
\node at (0,-1.7){\footnotesize $\mu$};
\node at (-5.0,4.2){\footnotesize $\infty$};
\node at (-1.95,4.4){\footnotesize $\infty$};
\node at (-1.6,4.3){\footnotesize $\infty$};
\node at (1.25,4.4){\footnotesize $\infty$};
\node at (1.65,4.45){\footnotesize $\infty$};
\node at (4.95,4.45){\footnotesize $\infty$};
\node at (-6.2,1.7){\footnotesize $\infty$};
\node at (-3.05,1.65){\footnotesize $\infty$};
\node at (0.05,1.85){\footnotesize $\infty$};
\node at (3.0,1.8){\footnotesize $\infty$};
\node at (6.1,1.8){\footnotesize $\infty$};
\node at (-6.2,-1.55){\footnotesize $\infty$};
\node at (-3.05,-1.5){\footnotesize $\infty$};
\node at (2.5,-1.35){\footnotesize $\infty$};
\node at (5.3,-1.35){\footnotesize $\infty$};
\node at (3.9,-2.95){\footnotesize $\infty$};
\node at (6.25,-2.95){\footnotesize $\infty$};
\node at (2.55,-4.5){\footnotesize $\infty$};
\node at (5.0,-4.5){\footnotesize $\infty$};
\node at (-1.1,-4.55){\footnotesize $\infty$};
\node at (1.75,-4.5){\footnotesize $\infty$};
\node at (-6.25,-4.1){\footnotesize $\infty$};
\node at (-3.2,-4.18){\footnotesize $\infty$};
 \end{tikzpicture}
\end{figure}

\begin{figure}[t!]\centering
\leftline{(iii) $\epsilon\in-\R_{>0}$:}
 \begin{tikzpicture}[scale=1]
\node at (0,0){\includegraphics{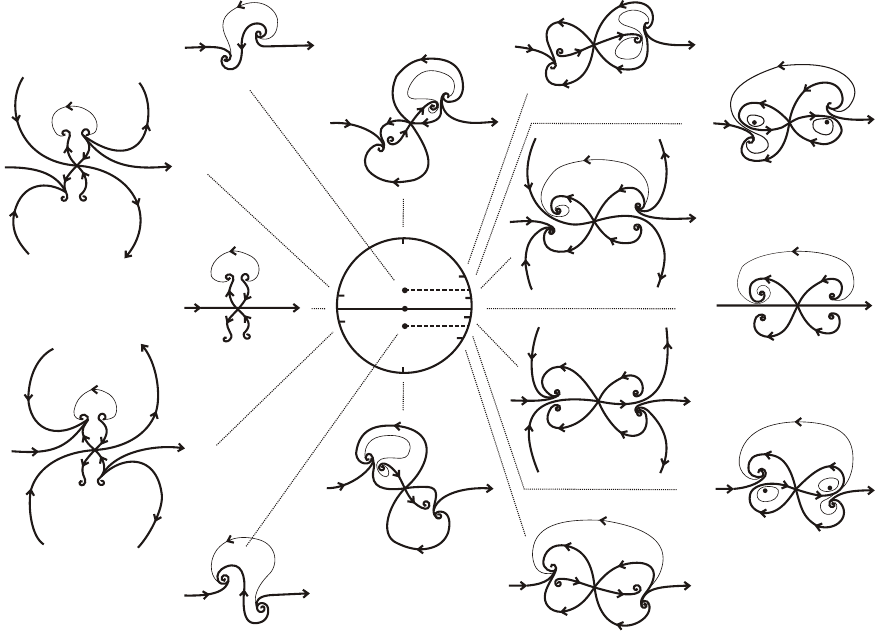}};
\node at (-1.0,0.65){\footnotesize $\mu$};
\node at (-4.5,4.55){\footnotesize $\infty$};
\node at (-1.9,4.55){\footnotesize $\infty$};
\node at (1.1,4.55){\footnotesize $\infty$};
\node at (4.55,4.6){\footnotesize $\infty$};
\node at (-7.5,2.5){\footnotesize $\infty$};
\node at (-4.3,2.5){\footnotesize $\infty$};
\node at (-2.0,3.25){\footnotesize $\infty$};
\node at (4.5,3.25){\footnotesize $\infty$};
\node at (4.6,1.65){\footnotesize $\infty$};
\node at (-4.45,0.13){\footnotesize $\infty$};
\node at (4.55,0.16){\footnotesize $\infty$};
\node at (4.50,-1.46){\footnotesize $\infty$};
\node at (4.50,-2.96){\footnotesize $\infty$};
\node at (4.50,-4.56){\footnotesize $\infty$};
\node at (1.00,-4.56){\footnotesize $\infty$};
\node at (-2.10,-2.92){\footnotesize $\infty$};
\node at (-2.00,-4.71){\footnotesize $\infty$};
\node at (-4.50,-4.76){\footnotesize $\infty$};
\node at (-4.10,-2.26){\footnotesize $\infty$};
\node at (-7.40,-2.31){\footnotesize $\infty$};
 \end{tikzpicture}
\leftline{(iv) $\epsilon\in {\rm i}\R_{>0}$:}
 \begin{tikzpicture}[scale=1]
\node at (0,0){\includegraphics{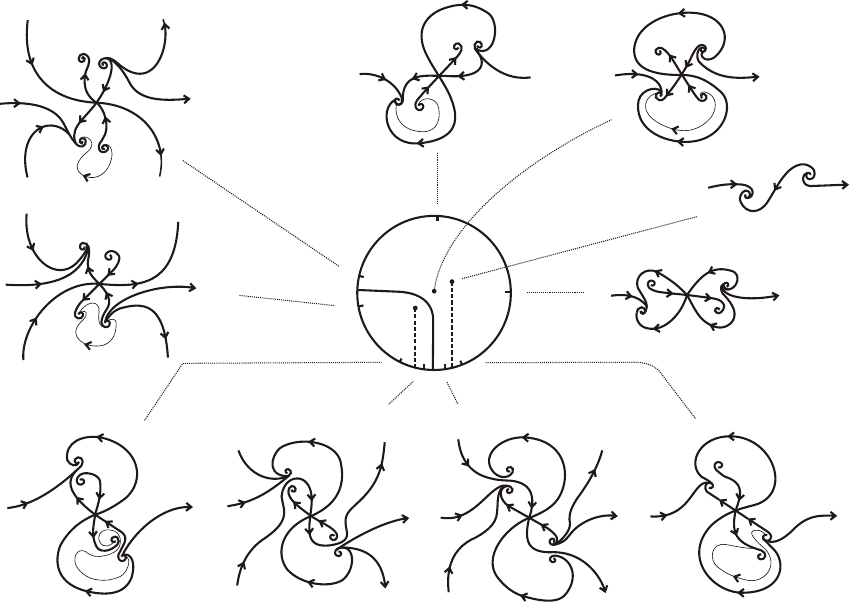}};
\node at (-0.5,0.85){\footnotesize $\mu$};
\node at (-7.4,3.35){\footnotesize $\infty$};
\node at (-3.75,3.4){\footnotesize $\infty$};
\node at (-1.3,3.85){\footnotesize $\infty$};
\node at (2.0,3.80){\footnotesize $\infty$};
\node at (3.05,3.80){\footnotesize $\infty$};
\node at (5.85,3.80){\footnotesize $\infty$};
\node at (4.6,1.9){\footnotesize $\infty$};
\node at (7.35,1.95){\footnotesize $\infty$};
\node at (2.95,0.10){\footnotesize $\infty$};
\node at (6.25,0.10){\footnotesize $\infty$};
\node at (-7.25,0.27){\footnotesize $\infty$};
\node at (-3.65,0.24){\footnotesize $\infty$};
\node at (-7.35,-3.5){\footnotesize $\infty$};
\node at (-3.65,-3.5){\footnotesize $\infty$};
\node at (-0.00,-3.65){\footnotesize $\infty$};
\node at (3.55,-3.65){\footnotesize $\infty$};
\node at (7.20,-3.65){\footnotesize $\infty$};
 \end{tikzpicture}
\end{figure}

\begin{figure}[t]\centering
\leftline{(0) $\epsilon=0$:}
 \begin{tikzpicture}[scale=1]
\node at (0,0){\includegraphics{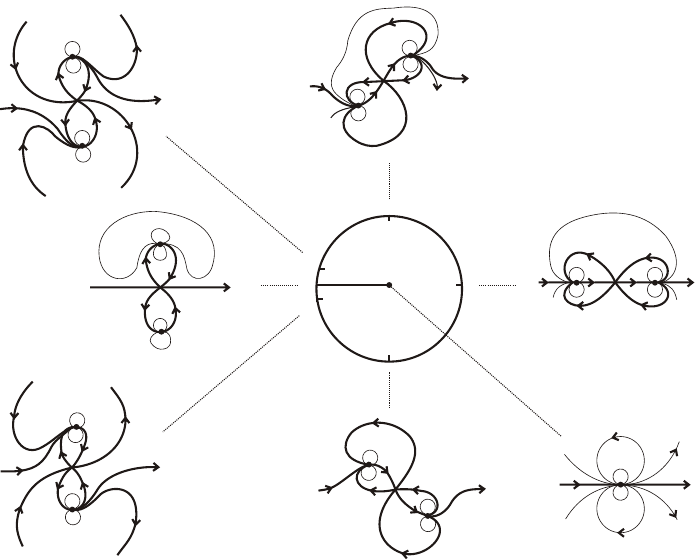}};
\node at (1.0,0.6){\footnotesize $\mu$};
\node at (-6.1,2.9){\footnotesize $\infty$};
\node at (-2.9,3.05){\footnotesize $\infty$};
\node at (-0.8,3.27){\footnotesize $\infty$};
\node at (2.25,3.40){\footnotesize $\infty$};
\node at (-4.55,-0.12){\footnotesize $\infty$};
\node at (-1.75,-0.12){\footnotesize $\infty$};
\node at (3.08,-0.03){\footnotesize $\infty$};
\node at (6.10,-0.03){\footnotesize $\infty$};
\node at (3.40,-3.48){\footnotesize $\infty$};
\node at (6.03,-3.48){\footnotesize $\infty$};
\node at (-0.70,-3.56){\footnotesize $\infty$};
\node at (2.56,-3.53){\footnotesize $\infty$};
\node at (-6.10,-3.23){\footnotesize $\infty$};
\node at (-2.96,-3.17){\footnotesize $\infty$};
 \end{tikzpicture}
\caption{The real phase portrait of the vector field $\chi$ according to $\mu$ for selected values of $\epsilon$ (see Fig.~\ref{figure:RC-1}). The separatrices of $0$ and $\infty$ are in bold. The splitting trajectories for the inner zones are not shown in the picture.}\label{figure:RC-2}
\end{figure}

\vspace*{15mm}

\subsubsection[Construction of the ramified domains $\Omega$ revisited]{Construction of the ramified domains $\boldsymbol{\Omega}$ revisited}\label{sec:RC-2.22}

Instead of in $\C^*$, let us consider now the real phase portrait of the vector field ${\rm e}^{{\rm i}\omega}\chi(s,\mu,\epsilon)$,
inside the pierced disc~$\sS^*$.
We have now zones of ${\rm e}^{{\rm i}\omega}\chi$ in $\sS^*$ (Definition~\ref{definition:RC-zones}) which are restrictions of those in $\C^*$, and half-zones in $\sS^*$ obtained by splitting the $\alpha\omega$-zones in two parts by the same trajectory as before. Instead of the values of $\omega$ for which ${\rm e}^{{\rm i}\omega}\chi$ is rotationally stable, we consider the admissible values of $\omega$ (Definition~\ref{definition:RC-zones}): let us remark that if the hyperbolic point at $s=0$ had a homoclinic separatrix this separatrix would reach it in both forward and backward direction in a finite time~$\theta$, thus leave~$\sS^*$.

Again, for a generic value of $(\mu,\epsilon)$ there are up to 4 connected zones or 6 connected half-zones in $\sS^*$: a symmetric pair of inner zones, each consisting of two half-zones, denote them
\[\sR_{I,\omega}(\check\mu,\check\epsilon)\supseteq\sR_{I+,\omega}(\check\mu,\check\epsilon)\cup \sR_{I-,\omega}(\check\mu,\check\epsilon),\qquad \sR_{I,\omega}^\p(\check\mu,\check\epsilon)\supseteq\sR_{I+,\omega}^\p(\check\mu,\check\epsilon)\cup \sR_{I-,\omega}^\p(\check\mu,\check\epsilon),\]
and a symmetric pair of outer zones, denote them
\[\sR_{O,\omega}(\check\mu,\check\epsilon),\qquad \sR_{O,\omega }^\p(\check\mu,\check\epsilon).\]
By Remark~\ref{remark:RC-omega}, $\sR_{\bullet,\omega}(\check\mu,\check\epsilon)={\rm e}^{-{\rm i}\omega}\sR_{\bullet,0}\big({\rm e}^{2{\rm i}\omega}\check\mu,{\rm e}^{4{\rm i}\omega}\check\epsilon\big)$.

\begin{figure}[t]\centering
 \begin{tikzpicture}[scale=1]
\node at (-3.5,0){\includegraphics[width=0.35\textwidth]{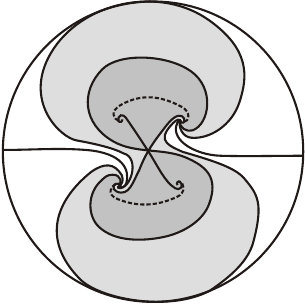}};
\node at (3.5,0){\includegraphics[width=0.35\textwidth]{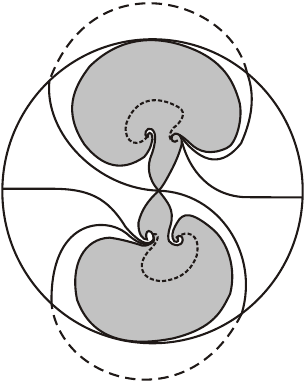}};
\node at (-3.5,-3.8){(a)};
\node at (3.5,-3.8){(b)};
\node at (-3.5,2.2){$\sR_{O}$};
\node at (-3.5,1.4){$\sR_{I+}$};
\node at (-3.5,0.7){$\sR_{I-}$};
\node at (-3.5,-0.65){$\sR_{I-}^\p$};
\node at (-3.5,-1.3){$\sR_{I+}^\p$};
\node at (-3.5,-2.2){$\sR_{O}^\p$};
\node at (3.5,2.2){$\sR_{I+}$};
\node at (3.7,1.4){$\sR_{I-}$};
\node at (3.75,-1.3){$\sR_{I-}^\p$};
\node at (3.5,-2.2){$\sR_{I+}^\p$};
 \end{tikzpicture}
\caption{The outer and inner half-zones $\sR_{O,0}(\check\mu,\check\epsilon)$ and $\sR_{I\pm,0}(\check\mu,\check\epsilon)$
(with $\omega=0$) inside the pointed disc $\sS^*$ for (a) $\epsilon\in {\rm i}\R_{>0}$, $\mu=0$,
 (b)~$\epsilon\in {\rm i}\R_{>0}$, $\mu$ close to $\Sigma_I(\epsilon)$: $\sR_{O,0}=\sR_{O,0 }^\p=\varnothing$. Compare with the corresponding vector fields in Fig.~\ref{figure:RC-2}(ii).}\label{figure:RC-zones}
\end{figure}

Let us agree that out of the two inner zones, $\sR_{I,\omega}(\check\mu,\check\epsilon)$ is the one consisting of trajectories from $s_1(\check\mu,\check\epsilon)$ to
$s_2(\check\mu,\check\epsilon)$~\eqref{eq:RC-si},
and that out of the two outer zones (both consisting of trajectories from $s_1(\check\mu,\check\epsilon)$ to $-s_1(\check\mu,\check\epsilon)$),
$\sR_{O,\omega}(\check\mu,\check\epsilon)$ is the upper one (see Fig.~\ref{figure:RC-zones}(a)).

The outer zones $\sR_{O,\omega}(\check\mu,\check\epsilon)$ can became empty: this happens whenever a separatrix
of the origin leaves the disc $\sS$ (see Fig.~\ref{figure:RC-zones}(b)).
Therefore a bifurcation of the zone $\sR_{O,\omega}$ occurs when a~separatrix of the origin touches the boundary of the disc from inside for the first time:
at that moment the zone ceases to exist as there is no trajectory
of ${\rm e}^{{\rm i}\omega}\chi$ joining $s_i(\check\mu,\check\epsilon)$ and $-s_j(\check\mu,\check\epsilon)$ inside the disc. We have:

\begin{Lemma}
	For a given $(\mu,\epsilon)\in\cM\times\cE$, a value of $\omega\in \big] {-} \frac{\pi}{2},\frac{\pi}{2}\big[$ is admissible if and only if
	the vector field ${\rm e}^{{\rm i}\omega}\chi$ has no centers, i.e., $\big({\rm e}^{-2{\rm i}\omega}\mu,{\rm e}^{-4{\rm i}\omega}\epsilon\big)\notin\Sigma_I$, and the outer zones $\sR_{O,\omega}(\check\mu,\check\epsilon)$, $\sR_{O,\omega}^\p(\check\mu,\check\epsilon)$ are non-empty. In this case the inner half-zones $\sR_{I\pm,\omega}(\check\mu,\check\epsilon)$, $\sR_{I\pm,\omega}^\p(\check\mu,\check\epsilon)$ agree with the inner half-zones of ${\rm e}^{{\rm i}\omega}\chi$ in $\C^*$.
\end{Lemma}

\begin{figure}[t]\centering
 \begin{tikzpicture}[scale=1]
\node at (0,0){\includegraphics[width=0.8\textwidth]{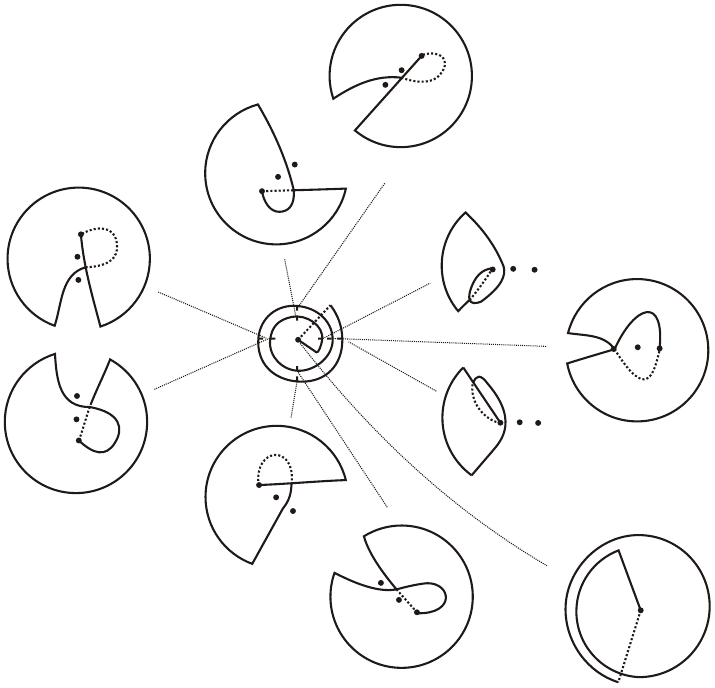}};
\node at (-6.0,2.9){$\check \cM(\tilde\epsilon)$};
\node at (-5.35,2.0){$\sqrt{\tilde\epsilon}$};
\node at (-6.0,0.0){$\check \cM(\bar\epsilon)$};
\node at (-5.35,-1.8){$\sqrt{\bar\epsilon}$};
\node at (-1.6,0.9){$\check \cE$};
\node at (4.2,-3.2){$\check \cM(0)$};
\node at (5.3,-4.8){$0$};
 \end{tikzpicture}
\caption{The ramified domains $\check \cM(\check\epsilon)$ for the parameter $\check\mu$ depending on $\check\epsilon\in\check \cE$.}\label{figure:RC-3}
\end{figure}

Corresponding to the inner and outer half-zones of the vector field $\chi$ we will construct domains $\Omega(\check\mu,\check\epsilon)$:
 two symmetric pairs of inner domains $\Omega_{I\pm}(\check\mu,\check\epsilon)$, $\Omega_{I\pm}^\p(\check\mu,\check\epsilon)$, and a symmetric pair of outer domains $\Omega_O(\check\mu,\check\epsilon)$, $\Omega_O^\p (\check\mu,\check\epsilon)$,
each obtained as a ramified union of the respective half-zones $\sR_{I\pm,\omega}(\check\mu,\check\epsilon)$, $\sR_{I\pm,\omega}^\p (\check\mu,\check\epsilon)$ and $\sR_{O,\omega}(\check\mu,\check\epsilon)$, $\sR_{O,\omega}^\p (\check\mu,\check\epsilon)$ over varying admissible $\omega$.
They will experience the same kind of bifurcations as their corresponding half-zones $\sR_{\bullet,\omega}$, but this time delayed by the effect of the variation of $\omega\in \big] {-} \frac{\pi}{2},\frac{\pi}{2}\big[$.
This will determine the set of ramified parameters $(\check\mu,\check\epsilon)$ for which they exist (Fig.~\ref{figure:RC-3}).

\begin{Lemma}\label{lemma:RC-6}
	If $\delta_\mu+\delta_\epsilon^\frac{1}{2}$ is small enough $($with respect to $\delta_s^2)$, then for each
	parameter $(\mu,\epsilon)$ there is an admissible $\omega\in \big] {-} \frac{\pi}{2},\frac{\pi}{2}\big[$.
\end{Lemma}

\begin{proof}
	We need to show that there exists $\omega\in \big] {-} \frac{\pi}{2},\frac{\pi}{2}\big[$, such that, if $\mu^2\neq\epsilon$, no separatrix of 0 of the vector field ${\rm e}^{{\rm i}\omega}\chi$ leaves the disc $\sS$ of radius $\delta_s$.
	The complement of the disc $\sS$ corresponds in the map $\theta$~\eqref{eq:RC-theta} to a~``disc'' centered at 0 of a radius uniformly bounded w.r.t.~$(\mu,\epsilon)$. For $|\mu|$, $|\epsilon|^{^\frac{1}{2}}$ small, all the preimages $\xi_0$ of the point $s=0$ by $\theta^{-1}$ are far enough (they depend continuously on the parameter and tend to $\infty$ as $(\mu,\epsilon)\to 0$)
	so that for some $\omega\in \big] {-} \frac{\pi}{2},\frac{\pi}{2}\big[$ none of the lines $\xi_0+{\rm e}^{-{\rm i}\omega}\R$ cross the~``disc''.
\end{proof}

{\bf The ramified parameter space.}
Let $\cM$, $\cE$ be as in \eqref{eq:RC-sEsM} small discs of radii $\delta_\mu$, $\delta_\epsilon$ in the~$\mu$- and~$\epsilon$-spaces.
Define a ramified sectoral cover $\check \cE$ of $\cE$ as
\[\check \cE=\{|\check\epsilon|<\delta_\epsilon,\, |\arg\check\epsilon|<4\pi\},\]
with each $\check\epsilon$ being projected to $\epsilon\in \cE$.
For each value of $\omega\in \big] {-} \frac{\pi}{2},\frac{\pi}{2}\big[$ and $\check\epsilon\in\check \cE$ such that
$|\arg\check\epsilon+4\omega|<2\pi$,
let $\cM_\omega(\check\epsilon)$ denote the connected component of the set{\samepage
\[\big\{\mu\in \cM\,|\, \omega\ \text{is admissible for } {\rm e}^{{\rm i}\omega}\chi(s,\mu,\epsilon)\big\}\]
that is attached to the point $\mu = \sqrt{\check\epsilon}$.
By Remark~\ref{remark:RC-omega}
$\check\mu\in\cM_\omega(\check\epsilon)\Longleftrightarrow {\rm e}^{2{\rm i}\omega}\check\mu\in\cM_0\big({\rm e}^{4{\rm i}\omega}\check\epsilon\big)$.}

As seen in Fig.~\ref{figure:RC-OI}, the component of $\cM\setminus\big({\rm e}^{-2{\rm i}\omega}\Sigma_O\big({\rm e}^{4{\rm i}\omega}\check\epsilon\big)\cup {\rm e}^{-2{\rm i}\omega}\Sigma_I\big({\rm e}^{4{\rm i}\omega}\check\epsilon\big)\big)$ attached to $\mu=\sqrt\epsilon$ undergoes a bifurcation for $|\arg\check\epsilon+4\omega|=2\pi$ when it ceases to exist, the corresponding bifurcation for $\cM_0\big({\rm e}^{4{\rm i}\omega}\check\epsilon\big)$
happens a bit earlier. For given $\check\epsilon\in\check\cE$, the set of $\omega$ for which $\cM_\omega(\check\epsilon)\neq 0$ is a proper subinterval of
$\{|\omega|<\frac{\pi}{2},\, |\arg\check\epsilon+4\omega|<2\pi\}$.

Define a domain $\check \cM(\check\epsilon)$ of ramified parameter $\check\mu$ as a ramified union
\[\check \cM(\check\epsilon)=\bigcup_{\substack{\omega\in {]} {-} \frac{\pi}{2},\frac{\pi}{2}[\\ |\arg\check\epsilon+4\omega|<2\pi}} \cM_\omega(\check\epsilon)\cup\big\{\sqrt\epsilon\big\}=\bigcup_{\substack{\omega\in {]} {-} \frac{\pi}{2},\frac{\pi}{2}[\\ |\arg\check\epsilon+4\omega|<2\pi}} {\rm e}^{-2{\rm i}\omega}\cM_0\big({\rm e}^{4{\rm i}\omega}\check\epsilon\big)\cup\big\{\sqrt\epsilon\big\},\]
with $\check\mu = \sqrt{\check\epsilon}$ as the ramification point included in $\check \cM(\check\epsilon)$.
See Fig.~\ref{figure:RC-3}.

\begin{Proposition}The union of $\check \cM(\check\epsilon)$ in the $(\check\mu,\check\epsilon)$-space
\begin{gather}\label{eq:RC-checksM}
\coprod_{\check\epsilon\in\check \cE} \check \cM(\check\epsilon)
\end{gather}
is a single simply connected ramified cover of $\cM \times \cE$ with ramification at $\big\{\epsilon\big(\mu^2-\epsilon\big)=0\big\}$.
On this domain, an admissible value of $\omega$ can be chosen for every $(\check\mu,\check\epsilon)$ that varies continuously.
\end{Proposition}

\begin{proof}The simple connectedness is by definition, the fact that it covers all $\cM \times \cE$ is a consequence of Lemma \ref{lemma:RC-6}.
\end{proof}

\begin{Lemma}\label{lemma:RC-arg}
If $\check\epsilon\in\check\cE$, $\check\mu\in\check\cM(\check\epsilon)$, then
\[\left|\arg\frac{\sqrt{\mu\pm\sqrt\epsilon}}{\sqrt\epsilon}\right|<\pi,\qquad \left|\arg\frac{\sqrt{\mu+\sqrt\epsilon}\pm\sqrt{\mu-\sqrt\epsilon}}{2\sqrt\epsilon}\right|<\pi.\]
\end{Lemma}

\begin{proof}From the definition, Fig.~\ref{figure:RC-1} and Remark~\ref{remark:RC-omega}, we see that if $\check\mu\in\check\cM_\omega(\check\epsilon)$ then
\[\left|\arg\frac{\sqrt{\mu\pm\sqrt\epsilon}}{\sqrt\epsilon}-\omega\right|<\frac{\pi}{2},\qquad
\left|\arg\frac{\sqrt{\mu+\sqrt\epsilon}\pm\sqrt{\mu-\sqrt\epsilon}}{2\sqrt\epsilon}-\omega\right|<\frac{\pi}{2}.\tag*{\qed} \]\renewcommand{\qed}{}
\end{proof}

We define a simply connected ramified domain $\rM$ over the $\check m$-space, covering a neighborhood~$\sM$ of 0 in the $m$-space with ramification at
 $\big\{m\in\sM\,|\, \epsilon(m)\big(\mu(m)-\epsilon(m)^2\big)=0\big\}$,
by lifting the map $m\mapsto(\mu(m),\epsilon(m))$ to the ramified cover~\eqref{eq:RC-checksM}
\begin{displaymath}
 \xymatrix{ \rM\ni \check m \ \ar@{|->}[r]\ar@{|->}[d] & \ (\check\mu(\check m),\check\epsilon(\check m))\ar@{|->}[d]\mathrlap{ \in \coprod_{\check\epsilon\in\check \cE} \check \cM(\check\epsilon)} \\
 \sM\ni m \ \ar@{|->}[r] & \ (\mu(m),\epsilon(m))\mathrlap{ \in\cM \times \cE.}}
\end{displaymath}

{\bf The ramified domains $\boldsymbol{\Omega(\check\mu,\check\epsilon)}$.}
For each $\check\epsilon\in\check\cE,\ \check\mu\in\check\cM(\check\epsilon)$ let
\begin{gather*}
 \Omega_\bullet(\check\mu,\check\epsilon)=\bigcup_{\substack{\omega \text{ such that}\\ \mu\in\cM_\omega(\check\epsilon)}} \overline{\sR}_{\bullet,\omega}(\check\mu,\check\epsilon),\qquad \bullet=O,I,\quad \text{resp.} \ O,I\pm
\end{gather*}
be a ramified union of the topological closures $\overline{\sR}_{\bullet,\omega}$ of the zones, resp.~half-zones, $\sR_{\bullet,\omega}$ in the $\check s$-variable.

\begin{figure}[t]\centering
 \begin{tikzpicture}[scale=1]
\node at (-5.2,0){\includegraphics[width=0.3\textwidth]{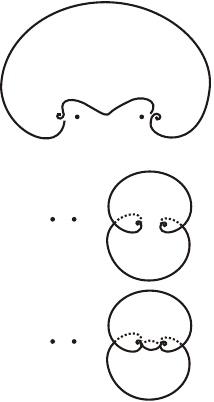}};
\node at (0,1.55){\includegraphics[width=0.3\textwidth]{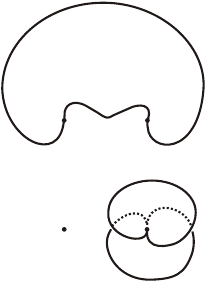}};
\node at (5.2,3.05){\includegraphics[width=0.3\textwidth]{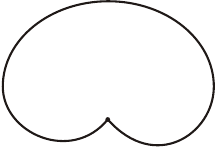}};
\node at (-5.2,-4.9){(a) $0<\epsilon<\mu^2$, $\mu>0$};
\node at (0.0,-4.9){(b) $\epsilon=0<\mu$};
\node at (5.2,-4.9){(c) $\epsilon=\mu=0$};
\node at (-5.2,3.4){$\Omega_O(\check\mu,\check\epsilon)$};
\node at (0,3.4){$\Omega_O(\check\mu,0)$};
\node at (5.2,3.4){$\Omega_O(0,0)$};
\node at (-6.65,1.75){$-s_1$};
\node at (-5.65,1.75){$-s_2$};
\node at (-4.6,1.75){$s_2$};
\node at (-3.8,1.75){$s_1$};
\node at (-1.35,2.0){$-s_1$};
\node at (1.3,1.95){$s_1$};
\node at (5.0,2.2){$0$};
\node at (-5.9,0.5){$\Omega_I(\check\mu,\check\epsilon)$};
\node at (-6.65,-0.8){$-s_1$};
\node at (-5.75,-0.8){$-s_2$};
\node at (-4.6,-0.9){$s_2$};
\node at (-3.8,-0.9){$s_1$};
\node at (-0.7,0.5){$\Omega_{I+}(\check\mu,0)$};
\node at (-0.9,-0.8){$-s_1{=}-s_2$};
\node at (1.55,-0.45){$s_1{=}s_2$};
\node at (-0.7,-1.5){$\Omega_{I-}(\check\mu,0)$};
\node at (-5.9,-2.3){$\Omega_{I+}(\check\mu,\check\epsilon)$};
\node at (-6.65,-3.5){$-s_1$};
\node at (-5.75,-3.5){$-s_2$};
\node at (-4.6,-3.6){$s_2$};
\node at (-3.8,-3.6){$s_1$};
\node at (-5.9,-4.3){$\Omega_{I-}(\check\mu,\check\epsilon)$};
 \end{tikzpicture}
\caption{Examples of the domains $\Omega_O(\check\mu,\check\epsilon)$ and $\Omega_I(\check\mu,\check\epsilon)$ for selected values of $\mu$, $\epsilon$.}\label{figure:RC-Omega}
\end{figure}

\emph{Following from their construction, the outer domain $\Omega_O(\check\mu,\check\epsilon)$ is connected nonempty for all $\check\epsilon\in\check\cM$, $\check\mu\in\check\cM(\check\epsilon)$,
while the inner domain $\Omega_I(\check\mu,\check\epsilon)=\Omega_{I + }(\check\mu,\check\epsilon)\cup\Omega_{I - }(\check\mu,\check\epsilon)$
becomes empty for $\check\mu^2=\check\epsilon$. For $\check\epsilon=0$, $\check\mu\neq 0$, the two parts $\Omega_{I\pm}(\check\mu,0)$ are disjoint except
for their common vertex at $ s_1(\check\mu,0)=s_2(\check\mu,0)= \sqrt{\check\mu}$
$($see Fig.~{\rm \ref{figure:RC-Omega})}.}

Finally, the ramified domain $\Omega_\bullet$ in the $(\check s,\check m)$-space is defined as the union of all $\Omega_\bullet(\check\mu,\check\epsilon)$, fibered over $\rM$,
\[\Omega_\bullet:=\coprod_{\check m\in\rM} \Omega_\bullet(\check\mu(\check m),\check\epsilon(\check m)).\]

\begin{Lemma}\label{lemma:RC-3}
For each $\check\mu\in\cM(\check\epsilon)$ the $6$ domains $\Omega_O(\check\mu,\check\epsilon)$, $\Omega_O^\p(\check\mu,\check\epsilon)$, $\Omega_{I\pm}(\check\mu,\check\epsilon)$, $\Omega_{I\pm}^\p(\check\mu,\check\epsilon)$
cover a full neighborhood of each equilibrium point~$s_i$ of $\chi(s,\mu,\epsilon)$ $(s_i\neq 0$ if $(\mu,\epsilon)\neq 0)$.
\end{Lemma}

\begin{proof}
If $s_i$ is a simple singularity, we can also suppose by the symmetry that it is attractive for ${\rm e}^{{\rm i}\omega}\chi$,
$\omega\in {]}\omega_1,\omega_2[\subseteq \big] {-} \frac{\pi}{2},\frac{\pi}{2}\big[$.
The map~$\theta$~\eqref{eq:RC-theta} is logarithmic on a small neighborhood of $s_i$ with a period $\frac{s_i\pi{\rm i}}{\sqrt\epsilon}$ and the image by $\theta$ of the inner domain attached to $s_i$ contains by construction a sector at $\infty$ where $\arg\theta\in {]} {-} \omega_2,-\omega_1[$, and hence also half-strips of any direction ${\rm e}^{-{\rm i}\omega}$ and of any width. Therefore it covers a full neighborhood of~$s_i$.

If $\epsilon=0$, $\mu\neq 0$, and $s_i$ is a double singularity, then the inner half-zones $\sR_{I\pm,\omega}$ contain each a~small disc attached to $s_i$ and centered in the direction
$\pm {\rm i}{\rm e}^{{\rm i}\omega}\R_{>0}$; their union over $\omega$ varying in some interval $]\omega_1,\omega_2[$ covers a full neighborhood of~$s_i$.

If $(\mu,\epsilon)=0$ ans $s_i=0$, then $\sR_{O,\omega}$ (resp.~$\sR_{O,\omega}^\p$) contains a small disc attached to~$s_i$ and centered in the direction
${\rm i}{\rm e}^{{\rm i}\omega}\R_{>0}$ (resp.\ $-{\rm i}{\rm e}^{{\rm i}\omega}\R_{>0}$), and their union over $\omega$ varying in some interval $]\omega_1,\omega_2[$ again covers a full neighborhood of $s_i=0$.
\end{proof}

\subsection[Connection matrices and proof of Theorem~\ref{theorem:RC-I}(a)]{Connection matrices and proof of Theorem~\ref{theorem:RC-I}(a)}

For the following discussion we will want to fix a branch $\Psi_\bullet$ of the fundamental solution $\Psi$ \eqref{eq:RC-Psi} of the diagonal system $\widebar\Delta^{s}$ on interior of each of the domains $\Omega_\bullet$. However, no single branch of~$\Psi$ converges as $\check\epsilon\to 0$ on the interior of both sectoral components $\Omega_{I + }(\check\mu,0)$ and $\Omega_{I - }(\check\mu,0)$ of $\Omega_I(\check\mu,0)$, $\check\mu\neq 0$.
This is one of the main reasons for splitting the inner domain $\Omega_I(\check\mu,\check\epsilon)$ in the two parts
$\Omega_{I + }(\check\mu,\check\epsilon)$ and $\Omega_{I - }(\check\mu,\check\epsilon)$.

\begin{Definition}Let $\Phi_1$, $\Phi_2$ be two fundamental matrix solutions of a linear system on two domains $U_1$, $U_2$ with connected non-empty intersection $U_1\cap U_2$.
We call the matrix $C=\Phi_1^{-1}\Phi_2$ \emph{a connection matrix} between $\Phi_1$ and $\Phi_2$ and represent it schematically as
\[\Phi_1\xrightarrow{\ C \ }\Phi_2.\]
\end{Definition}

\subsubsection[Choice of the fundamental solutions $\Psi_\bullet$]{Choice of the fundamental solutions $\boldsymbol{\Psi_\bullet}$}

\begin{figure}[t]\centering
 \begin{tikzpicture}[scale=1]
\node at (0,0){\includegraphics[width=0.72\textwidth]{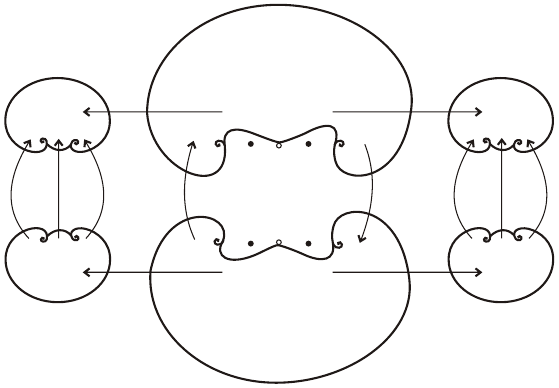}};
\node at (0,1.9){$\Psi_O$};
\node at (0,0.75){$0$};
\node at (0.55,0.75){$s_2$};
\node at (1.0,0.75){$s_1$};
\node at (0,-0.75){$0$};
\node at (0.5,-0.75){$s_2$};
\node at (0.95,-0.75){$s_1$};
\node at (0,-1.9){$\Psi_O^\p$};
\node at (-3.1,1.9){$I$};
\node at (3.1,1.9){$I$};
\node at (-4.5,1.7){$\Psi_{I-}^\p$};
\node at (4.6,1.7){$\Psi_{I+}$};
\node at (-3.1,-1.4){$I$};
\node at (3.1,-1.4){$I$};
\node at (-3.35,0.1){$N$};
\node at (3.35,0.1){$N$};
\node at (-4.2,0.1){$N_1$};
\node at (4.9,0.1){$N_1$};
\node at (-5.7,0.1){$I$};
\node at (5.7,0.1){$I$};
\node at (-2.1,0.1){$I$};
\node at (2.1,0.1){$I$};
\node at (-4.2,0.7){$-s_2$};
\node at (-4.9,0.7){$-s_1$};
\node at (-4.2,-0.5){$-s_2$};
\node at (-4.9,-0.5){$-s_1$};
\node at (4.3,0.7){$s_2$};
\node at (4.9,0.7){$s_1$};
\node at (4.3,-0.5){$s_2$};
\node at (4.9,-0.5){$s_1$};
\node at (-4.5,-1.6){$\Psi_{I+}^\p$};
\node at (4.6,-1.6){$\Psi_{I-}$};
 \end{tikzpicture}
\caption{The connection matrices between the fundamental solutions $\Psi_\bullet$ for each fixed parameter $(\check\mu,\check\epsilon)$, with $\check\mu^2\neq\check\epsilon\neq 0$, where $N_1$ and $N$ are given by \eqref{eq:RC-N1N2} and \eqref{eq:RC-N}. If $\check\epsilon=0$ then $s_1(\check\mu,0)=s_2(\check\mu,0)$ and the matrices $N_1$ are missing from the picture. If $\check\mu^2 = \check\epsilon$ then only the fundamental solutions~$\Psi_O$ and~$\Psi_O^\p$ persist.}\label{figure:RC-5}
\end{figure}

On the interior $\mathring\Omega_\bullet$ of each of the domains $\Omega_\bullet$, $\bullet=O,{I + },{I - }$, we fix a branch $\Psi_\bullet(\check s,\check\mu,\check\epsilon)$
of the fundamental solution $\Psi(\check s,\check\mu,\check\epsilon)$~\eqref{eq:RC-Psi} of the diagonal system
$\widebar\Delta^{s}$ so that the connection matrices between them are as in Fig.~\ref{figure:RC-5}.

The monodromy matrices of $\Psi(s,\check\mu,\check\epsilon)$ around the points $s_1(\check\mu,\check\epsilon)$, resp.\ $s_2(\check\mu,\check\epsilon)$, in the positive direction
are independent of the choice of the branch of $\Psi$, and are given by
\begin{gather}\label{eq:RC-N1N2}
N_1(\check\mu,\check\epsilon)=
\begin{pmatrix} {\rm e}^{\frac{s_1(\check\mu,\check\epsilon)}{\sqrt{\check\epsilon}}\pi{\rm i}} & 0 \\ 0 & {\rm e}^{-\frac{s_1(\check\mu,\check\epsilon)}{\sqrt{\check\epsilon}}\pi{\rm i}} \end{pmatrix},
\qquad
N_2(\check\mu,\check\epsilon)=
\begin{pmatrix} {\rm e}^{-\frac{s_2(\check\mu,\check\epsilon)}{\sqrt{\check\epsilon}}\pi{\rm i}} & 0 \\ 0 & {\rm e}^{\frac{s_2(\check\mu,\check\epsilon)}{\sqrt{\check\epsilon}}\pi{\rm i}} \end{pmatrix},
\end{gather}
They satisfy
\[N_i^\p=N_i^{-1},\qquad i=1,2.\]
The monodromy matrix of $\Psi$ around both of the points $s_1(\check\mu,\check\epsilon), s_2(\check\mu,\check\epsilon)$ is equal to
\begin{gather}\label{eq:RC-N}
N(\check\mu,\check\epsilon)=N_1(\check\mu,\check\epsilon)N_2(\check\mu,\check\epsilon)=
\begin{pmatrix} {\rm e}^{\frac{s_1-s_2}{\sqrt{\check\epsilon}}\pi{\rm i}} & 0 \\ 0 & {\rm e}^{-\frac{s_1-s_2}{\sqrt{\check\epsilon}}\pi{\rm i}} \end{pmatrix}.
\end{gather}
At the limit when $\check\epsilon\to 0$ we get
$ N(\check\mu,0)=\left(\begin{smallmatrix} {\rm e}^{\frac{1}{\sqrt{\check\mu}}\pi{\rm i}} & 0 \\ 0 & {\rm e}^{-\frac{1}{\sqrt{\check\mu}}\pi{\rm i}} \end{smallmatrix}\right)$,
which is for $\mu\neq 0$ the monodromy matrix of $\Psi$ around the double zero $s_1(\check\mu,0)=s_2(\check\mu,0)$. On the other hand,
none of the matrices $N_1(\check\mu,\check\epsilon)$, $N_2(\check\mu,\check\epsilon)$ has a limit as $\check\epsilon\to 0$.

\subsubsection[Connection matrices of the fundamental solutions $\Phi_\bullet=F_\bullet\Psi_\bullet$]{Connection matrices of the fundamental solutions $\boldsymbol{\Phi_\bullet=F_\bullet\Psi_\bullet}$}

Let $F_\bullet(\check s,\check m)$ and $F_\bullet^\p(\check s,\check m)$ be the constructed diagonalizing gauge transformations on the domains $\Omega_\bullet$ and $\Omega_\bullet^\p={\rm e}^{\pi{\rm i}}\Omega_\bullet$, $\bullet=O,I$, in the ramified $(\check s,\check m)$-space. Namely, the restrictions~$F_{I\pm}$ of $F_{I}$ on the two sub-domains $\Omega_{I\pm}$ of~$\Omega_I$ agree on the common part between the two singularities~$s_1$,~$s_2$.
Then{\samepage
\[F_\bullet\Psi_\bullet\qquad\text{and}\qquad F_\bullet^\p\Psi_\bullet^\p, \qquad \bullet=O,I\pm,\]
are fundamental solution matrices of $\Delta^{ s}$~\eqref{eq:RC-d7a}.}

Whenever a point $(s,m)\in\sS\times\sM$ is covered more than once, then there is a connection matrix between these fundamental solutions:
Either there can be two domains $\Omega(\check\mu,\check\epsilon)$ with the same $(\check\mu,\check\epsilon)$,
or with two different ramified parameters $(\check\mu,\check\epsilon)$ corresponding to the same $(\mu,\epsilon)$.
The collection of all these connection matrices carries all the information about the analytic equivalence class of the system $\Delta$.

In Lemmas~\ref{lemma:RC-4} and~\ref{lemma:RC-5} and Proposition~\ref{proposition:RC-5} we give a semi-explicit expression of all the connection matrices.

\begin{Proposition}\label{proposition:RC-4}
Let $\Delta(x,m)$, $\Delta'(x,m)$ be two parametric systems and let
$\Delta^{ s}(s,m)$, $\Delta'^{s}(s,m)$ be their transforms in the $s$-coordinate \eqref{eq:RC-d7a}.
Let $F_\bullet$, $F'_\bullet$ be normalizing gauge transformations for $\Delta^{ s}$, $\Delta'^{s}$:
\[(F_\bullet)^*\Delta^{ s}=\widebar\Delta^{s}=(F_\bullet')^*\Delta'^{s}\]
on the domains $\bullet=O,{I + },{I - }$ defined above.
If all the connection matrices associated to the fundamental solutions $F_\bullet\Psi_\bullet$ of $\Delta^{ s}$ agree with those associated to the fundamental solutions $F'_\bullet\Psi_\bullet$ of $\Delta'^{s}$,
then the two parametric families of systems $\Delta$, $\Delta'$ are analytically equivalent.
\end{Proposition}

\begin{proof}
Let $H(s,m):=F'_\bullet(\check s,\check m) F_\bullet(\check s,\check m)^{-1}$. Since all the connection matrices are equal, $H$~is a~well defined non-ramified invertible matrix function
defined on the union of the projections of the domains $\Omega_\bullet$ to $(s,m)$-space, $\bullet=O,{I + },{I - }$.
It is bounded on a neighborhood of each singularity $s_i\neq 0$, hence $H$ can be analytically extended on
$(\sS\setminus\{0\})\times \cM \times \cE$, where~$\sS$,~$\cM$,~$\cE$ are as in~\eqref{eq:RC-sEsM}.
It satisfies $H=H^\p$:
if $s$ is in the projection of $\Omega_\bullet(\check\mu,\check\epsilon) $ and
$ H(s,\check m)=F'_\bullet(\check s,\check m) F_\bullet(\check s,\check m)^{-1} $
then $-s$ is in the projection of $\Omega_\bullet^\p $ and
\begin{align*} H(-s,\check m)&={F'_\bullet}^\p\big({\rm e}^{\pi{\rm i}}\check s,\check m\big) \big(F_\bullet^\p\big({\rm e}^{\pi{\rm i}}\check s,\check m\big)\big)^{-1} \\
&=\left(\begin{smallmatrix} 0 & 1 \\ 1 & 0 \end{smallmatrix}\right) F'_\bullet(\check s,\check m) F_\bullet(\check s,\check m)^{-1} \left(\begin{smallmatrix} 0 & 1 \\ 1 & 0 \end{smallmatrix}\right)\\
&=\left(\begin{smallmatrix} 0 & 1 \\ 1 & 0 \end{smallmatrix}\right) H(\check s,\check m) \left(\begin{smallmatrix} 0 & 1 \\ 1 & 0 \end{smallmatrix}\right).
\end{align*}
Hence the function
$G(x,m):=S(s,m) VH(s,m) V^{-1}S^{-1}(s,m)$, with $S$, $V$ as in~\eqref{eq:RC-SV}, is well defined.

The fundamental solutions $Y_\bullet(\check x,\check m)=S(s) VF_\bullet(\check s,\check m) \Psi_\bullet(\check s,\check m)$ of the system $\Delta(x,m)$, and
$Y'_\bullet(\check x,\check m)=S(s) VF'_\bullet(\check s,\check m) \Psi_\bullet(\check s,\check m)$ of $\Delta'(x,m)$,
can for $\mu^2\neq\epsilon$ be analytically extended on a~neighborhood of the point $x = -\mu$ (i.e., $s = 0$) which is non-singular for these systems.
As $G=Y'_\bullet Y_\bullet^{-1}$, it means that $ G^*\Delta'=\Delta$ and that
$G$ is an invertible analytic matrix function on $(\sX \times \sM) \setminus\big\{x=-\mu,\, \epsilon=\mu^2\big\}$, where
$\sX:=\big\{|x|\leq\delta_s^2 -\delta_\mu\big\}$.
Since the problematic points are in a set of codimension~2, by Hartog's theorem $G$ is analytic on the whole neighborhood $\sX \times \sM$ of~0.
\end{proof}

\begin{figure}[t]\centering
 \begin{tikzpicture}[scale=1]
\node at (0,0){\includegraphics[width=0.72\textwidth]{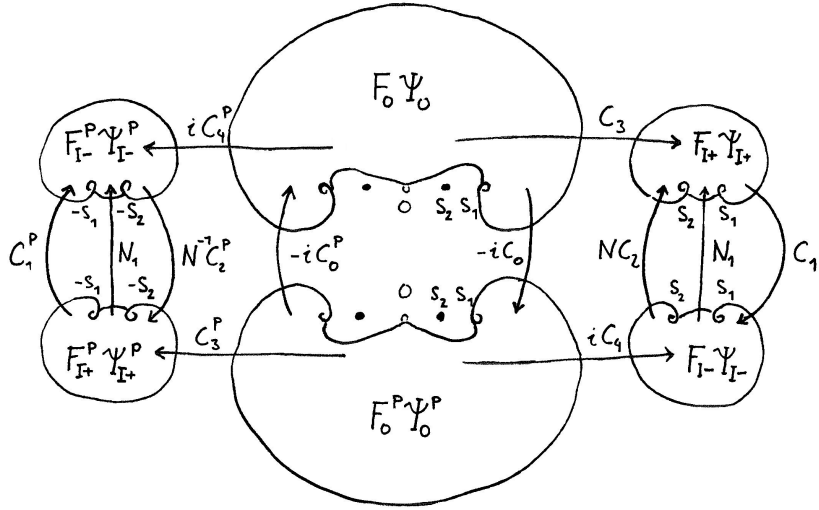}};
\node at (0,1.9){$F_O\Psi_O$};
\node at (0,0.75){$0$};
\node at (0.55,0.75){$s_2$};
\node at (1.0,0.75){$s_1$};
\node at (0,-0.75){$0$};
\node at (0.5,-0.75){$s_2$};
\node at (0.95,-0.75){$s_1$};
\node at (0,-1.9){$F_O^\p\Psi_O^\p$};
\node at (3.1,2.0){$C_3$};
\node at (-3.1,2.0){$C_4^\p$};
\node at (-4.7,1.7){$F_{I-}^\p\Psi_{I-}^\p$};
\node at (4.85,1.7){$F_{I+}\Psi_{I+}$};
\node at (-3.1,-1.3){$C_3^\p$};
\node at (3.1,-1.3){$C_4$};
\node at (-2.9,0.1){$N^{-1}C_2^\p$};
\node at (3.15,0.1){$NC_2$};
\node at (-4.2,0.1){$N_1$};
\node at (4.9,0.1){$N_1$};
\node at (-5.8,0.1){$C_1^\p$};
\node at (5.8,0.1){$C_1$};
\node at (-1.6,0.1){$C_0^\p$};
\node at (1.6,0.1){$C_0$};
\node at (-4.2,0.7){$-s_2$};
\node at (-4.9,0.7){$-s_1$};
\node at (-4.2,-0.5){$-s_2$};
\node at (-4.9,-0.5){$-s_1$};
\node at (4.3,0.7){$s_2$};
\node at (4.9,0.7){$s_1$};
\node at (4.3,-0.5){$s_2$};
\node at (4.9,-0.5){$s_1$};
\node at (-4.7,-1.6){$F_{I+}^\p\Psi_{I+}^\p$};
\node at (4.85,-1.6){$F_{I-}\Psi_{I-}$};
 \end{tikzpicture}
\caption{The connection matrices between the fundamental solutions $F_\bullet\Psi_\bullet$ for a fixed parameter $(\check\mu,\check\epsilon)$, $\check\mu^2\neq\check\epsilon$.
For $\check\mu^2=\check\epsilon$, only the fundamental solutions $F_O\Psi_O$ and $F_O^\p\Psi_O^\p$ persist, with the two corresponding connection matrices $C_0$, $C_0^\p$.
(Picture with $0<\check\epsilon<\check\mu^2$ as in Figs.~\ref{figure:RC-5} and \ref{figure:RC-Omega}(a)).}\label{figure:RC-6}
\end{figure}

\begin{Lemma} \label{lemma:RC-4}Let $F_\bullet$ be the normalizing gauge transformations from Proposition~{\rm \ref{proposition:RC-3}} satis\-fying~\eqref{eq:RC-kappaIO} with the uniquely determined functions $\kappa_\bullet$, and let $\Psi_\bullet$ be as Fig.~{\rm \ref{figure:RC-5}}. Then for each fixed $\check m\in\sM$ the connection matrices between the solutions $F_\bullet\Psi_\bullet$ on the domains $\Omega_\bullet(\check\mu,\check\epsilon)$ are given in Fig.~{\rm \ref{figure:RC-6}} with the matrices $C_0(\check m),\dots,C_4(\check m)$ equal to
\begin{alignat}{3}
& C_0=\left(\begin{smallmatrix} 1 & {\rm i}\gamma \\ 0 & 1 \end{smallmatrix} \right), \qquad &&
 C_3=\left(\begin{smallmatrix} 1 & {\rm i}\kappa^{-1} {\rm e}^{ -2a\pi{\rm i}} \\ 0 & \kappa^{-1} \end{smallmatrix} \right), & \nonumber\\
 & C_1=\left(\begin{smallmatrix} 1 & {\rm i}\kappa^{-1}(\gamma- {\rm e}^{ 2a\pi{\rm i}} - {\rm e}^{ -2a\pi{\rm i}}) \\ 0 & 1\end{smallmatrix} \right) , \qquad &&
 C_4=\left(\begin{smallmatrix} 1 & -{\rm i}\kappa^{-1} {\rm e}^{ 2a\pi{\rm i}} \\ 0 & \kappa^{-1} \end{smallmatrix} \right),&\nonumber\\
 & C_2=\left(\begin{smallmatrix} 1 & 0 \\ -{\rm i}\kappa {\rm e}^{ 2a\pi{\rm i}} & 1 \end{smallmatrix} \right) ,&&&\label{eq:RC-Ci}
\end{alignat}
where
\begin{gather}
a(\check m) :=
 \begin{cases}
 \dfrac{s_1(\check\mu,\check\epsilon)-s_2(\check\mu,\check\epsilon)}{2\sqrt{\check\epsilon}}=\frac{\sqrt{\mu+\sqrt\epsilon}-\sqrt{\mu-\sqrt\epsilon}}{2\sqrt\epsilon} & \text{if} \ \check\epsilon\neq 0,\vspace{1mm}\\
 \dfrac{1}{2\sqrt{\check\mu}} & \text{if}\ \check\epsilon=0\ \text{and}\ \check\mu\neq 0,
 \end{cases}\label{eq:RC-a} \\
\kappa(\check m) :=\frac{\kappa_O(\check m)}{\kappa_I(\check m)} \label{eq:RC-kappa}
\end{gather}
and $\gamma(m)$, the analytic invariant of the system $\Delta(x,m)$, is the trace of monodromy~\eqref{eq:RC-gamma}.
\end{Lemma}

\begin{proof}From Lemma \ref{lemma:RC-boundedsolutions}(b) we know that a connection matrix on an intersection domain that is adjacent to the point $s_1(\check\mu,\check\epsilon)$ (resp.~$s_2(\check\mu,\check\epsilon)$)
must be upper triangular (resp.\ lower triangular),
with the diagonal terms determined by the values of the corresponding pair of gauge transformations $F_\bullet(s_1(\check\mu,\check\epsilon),\check m)$
(resp.\ $F_\bullet(s_2(\check\mu,\check\epsilon),\check m)$). Hence we have
\[C_0=\left(\begin{smallmatrix} 1 & c_0 \\ 0 & 1 \end{smallmatrix}\right),\qquad
C_1=\left(\begin{smallmatrix} 1 & c_1 \\ 0 & 1 \end{smallmatrix}\right),\qquad
C_2=\left(\begin{smallmatrix} 1 & 0 \\ c_2 & 1 \end{smallmatrix}\right),\qquad
C_3=\left(\begin{smallmatrix} 1 & c_3 \\ 0 & \kappa^{-1} \end{smallmatrix}\right),\qquad
C_4=\left(\begin{smallmatrix} 1 & c_4 \\ 0 & \kappa^{-1} \end{smallmatrix}\right),\]
for some $c_0(\check m),\dots,c_4(\check m)$.

Let $M(\check m)$ be the monodromy matrix of the fundamental solution
\[ Y_O(\check x,\check m)=S(s) VF_O(\check s,\check m) \Psi_O(\check s,\check\mu,\check\epsilon) \]
of the system $\Delta$ around the two singular points $x=\pm\sqrt{\check\epsilon}$ in the positive direction. On the one hand we have
\begin{gather*}
M =\Psi_O(\check s)^{-1}F_O(\check s)^{-1}V^{-1}S(\check s)^{-1}\cdot S\big({\rm e}^{\pi{\rm i}}\check s\big) VF_O\big({\rm e}^{\pi{\rm i}}\check s\big) \Psi_O\big({\rm e}^{\pi{\rm i}}\check s\big)\\
\hphantom{M}{} =\Psi_O(\check s)^{-1}F_O(\check s)^{-1}V^{-1}S(\check s)^{-1}\cdot S\big({\rm e}^{\pi{\rm i}}\check s\big)\left(\begin{smallmatrix} -{\rm i} & 0 \\ 0 & {\rm i} \end{smallmatrix}\right)
VF_O^\p\big({\rm e}^{\pi{\rm i}}\check s\big) \Psi_O^\p\big({\rm e}^{\pi{\rm i}}\check s\big) C_0^\p\\
\hphantom{M}{} =-{\rm i}\left(\begin{smallmatrix} 0 & 1 \\ 1 & 0 \end{smallmatrix}\right) C_0^\p=-{\rm i} C_0 \left(\begin{smallmatrix} 0 & 1 \\ 1 & 0 \end{smallmatrix}\right),
\end{gather*}
using that $\left(\begin{smallmatrix} 1 & 0 \\ 0 & -1 \end{smallmatrix}\right)V=V\left(\begin{smallmatrix} 0 & 1 \\ 1 & 0 \end{smallmatrix}\right)$.
On the other hand, as apparent from Fig.~\ref{figure:RC-6},
\[M=C_3C_1NC_2C_3^{-1},\]
where $ N=\left(\begin{smallmatrix} {\rm e}^{ 2a\pi{\rm i}} & 0 \\ 0 & {\rm e}^{ -2a\pi{\rm i}} \end{smallmatrix}\right)$. Therefore
\begin{gather}\label{eq:RC-monodromy1}
-{\rm i} C_0 \left(\begin{smallmatrix} 0 & 1 \\ 1 & 0 \end{smallmatrix}\right) = C_3C_1NC_2C_3^{-1} =M, \\
\left(\begin{smallmatrix} -{\rm i}c_0 & -{\rm i} \\ -{\rm i} & 0 \end{smallmatrix}\right)=
\left(\begin{smallmatrix} {\rm e}^{ 2a\pi{\rm i}} + {\rm e}^{ -2a\pi{\rm i}}c_2(c_1 + c_3) & \kappa {\rm e}^{ -2a\pi{\rm i}}(c_1 + c_3)(1 - c_2c_3) - \kappa {\rm e}^{ 2a\pi{\rm i}}c_3 \\
\kappa^{-1}{\rm e}^{ -2a\pi{\rm i}}c_2 & {\rm e}^{ -2a\pi{\rm i}}(1 - c_2c_3) \end{smallmatrix}\right), \nonumber
\end{gather}
which implies that
\begin{gather}\label{eq:RC-gammac1c2}
\gamma=\tr M=-{\rm i} c_0={\rm e}^{ 2a\pi{\rm i}} + {\rm e}^{ -2a\pi{\rm i}} + {\rm e}^{ -2a\pi{\rm i}}c_1c_2,\\
c_2c_3=1,\qquad\text{and}\qquad c_2=-{\rm i}\kappa {\rm e}^{ 2a\pi{\rm i}}, \qquad c_3={\rm i}\kappa^{-1} {\rm e}^{ -2a\pi{\rm i}}.\nonumber\end{gather}
From Fig.~\ref{figure:RC-6} one also sees that
\begin{gather}\label{eq:RC-monodromy2}
C_3C_1=C_0C_4,
\end{gather}
which gives the matrix $C_4$.
\end{proof}

The matrices $C_0(\check m),\dots,C_4(\check m)$ of Lemma~\ref{lemma:RC-4} determine for each fixed $\check m\in\rM$ all the relations between the set of fundamental solutions
$F_\bullet(\cdot,\check m)\Psi_\bullet(\cdot,\check m)$ and $F_\bullet^\p(\cdot,\check m)\Psi_\bullet^\p(\cdot,\check m)$, $\bullet=O,I +,I -$.
We will now look at the situation of two different $\check m\in\rM$ corresponding to the same value of~$m$. One finds that the corresponding connection matrices can be expressed in terms of the values of $C_0,\dots,C_4$ for the two ramified parameters, while certain cocycle relations must be satisfied.

\begin{Lemma} \label{lemma:RC-5} Let $F_\bullet$, $\Psi_\bullet$ be as in Lemma~{\rm \ref{lemma:RC-4}}.
We will use the following kind of notation: If $\Bar m,\Barbar m\in\rM$ $($resp.\ $\Tilde m$, $\tilde{\tilde m}\in\rM)$ are two values of the ramified parameter $\check m$,
we write $\Bar X=X(\Bar m)$, $\Barbar X=X(\Barbar m)$ $($resp.\ $\Tilde X=X(\Tilde m)$, $\tilde{\tilde X}=X(\tilde{\tilde m}))$ for any object $X$ depending on~$\check m$.
\begin{itemize}\itemsep=0pt
\item[$(a)$] Let $\Bar m, \Barbar m\in\rM$ be two values of the ramified parameter that project to the same $m$, such that
\[\Bar\epsilon=\Barbar\epsilon=:\check\epsilon \qquad\text{and}\qquad \Barbar\mu-\sqrt{\check\epsilon}={\rm e}^{2\pi{\rm i}}\big(\Bar\mu-\sqrt{\check\epsilon}\big),\]
i.e., $\Barbar\mu$ is $\Bar\mu$ plus one positive turn around the ramification point $\sqrt{\check\epsilon}$ in $\check{\Cal M}(\check\epsilon)$.
So
\[\Barbar s_1=\Bar s_1,\qquad \Barbar s_2={\rm e}^{\pi{\rm i}}\Bar s_2.\]
Hence
\begin{gather*}
\Barbar\Omega_O=\Bar\Omega_O,\qquad \Barbar F_O=\Bar F_O,\qquad \Barbar \Psi_O=\Bar \Psi_O,
\end{gather*}
and we have
\begin{gather}\label{eq:RC-11}
\Barbar\kappa_O=\Bar\kappa_O=\frac{\Barbar\kappa \Bar\kappa }
{1-{\rm e}^{-2\frac{\Bar s_2}{\sqrt{\check\epsilon}}\pi{\rm i}}}.
\end{gather}

\item[$(b)$] Let $\Tilde m, \tilde{\tilde m}\in\rM$ be two values of the ramified parameter that project to the same $m$ such that
\[(\tilde{\tilde\mu},\tilde{\tilde\epsilon})={\rm e}^{2\pi{\rm i}}(\Tilde\mu,\Tilde\epsilon),\]
or more precisely, for $|\mu|\gg\sqrt{|\epsilon|}$, $(\tilde{\tilde\mu},\tilde{\tilde\epsilon})$ is obtained from $(\Tilde\mu,\Tilde\epsilon)$
by simultaneously turning both $\check\epsilon$ and $\check\mu$.
So
\[\tilde{\tilde s}_1={\rm e}^{\pi{\rm i}}\Tilde s_2,\qquad \tilde{\tilde s}_2={\rm e}^{\pi{\rm i}}\Tilde s_1,\qquad\text{and}\qquad \tilde{\tilde N}=\Tilde N^{-1}.\]
Hence
\begin{alignat*}{4}
& \tilde{\tilde\Omega}_{I +}=\Tilde\Omega_{I -}^\p ,\qquad && \tilde{\tilde F}_{I +}=\Tilde F_{I -}^\p , \qquad && \tilde{\tilde\Psi}_{I +}=\Tilde\Psi_{I -}^\p ,& \\
& \tilde{\tilde\Omega}_{I -}=\Tilde\Omega_{I +}^\p ,\qquad && \tilde{\tilde F}_{I -}=\Tilde F_{I +}^\p , \qquad && \tilde{\tilde\Psi}_{I -}=\Tilde\Psi_{I +}^\p \Tilde N^{-1}.&
\end{alignat*}
Therefore
\begin{gather}\label{eq:RC-C12}
\tilde{\tilde C}_1=\Tilde N^{-1}\Tilde C_2^\p \Tilde N,\qquad \tilde{\tilde C}_2=\Tilde C_1^\p,
\end{gather}
and we have
\begin{gather}\label{eq:RC-11b}
\tilde{\tilde\kappa}_{I}=\Tilde\kappa_I,\\
\label{eq:RC-12}
\gamma={\rm e}^{2\Tilde a\pi{\rm i}}+{\rm e}^{-2\Tilde a\pi{\rm i}}-\Tilde\kappa \tilde{\tilde\kappa} {\rm e}^{-2\Tilde a\pi{\rm i}},
\end{gather}
where $ \check a$ and $\check\kappa $ are defined in \eqref{eq:RC-a} and \eqref{eq:RC-kappa}, $\tilde{\tilde a}=-\Tilde a$.
\end{itemize}
\end{Lemma}

\begin{proof} (a) For each $\check\epsilon\in\check\cE$ the ramification of the $\check\mu$-parameter domain $\check{\Cal M}(\check\epsilon)$ corresponds to the bifurcation $\Sigma_I(\epsilon)$: the difference between $(\Bar\mu,\check\epsilon)$ and $(\Tilde\mu,\check\epsilon)$ is that of crossing the line $\Sigma_I(\epsilon)$. Since this bifurcation affects only the inner zones of the field $\chi$, it therefore affects only the internal domains $\Omega_{I \pm}$, $\Omega_{I \pm}^\p$, while the outer domains are not affected. Therefore $\Barbar\Omega_O=\Bar\Omega_O$ and consequently $\Barbar F_O=\Bar F_O$.

To obtain the assertion~\eqref{eq:RC-11}, it is enough to prove it for generic values of $(\mu,\epsilon)$, and extend it to the other values by continuity. So we can assume that $\epsilon\neq0$, $\mu^2\neq\epsilon$, and moreover that both of the points $s_1(\check\mu,\check\epsilon)$, $s_2(\check\mu,\check\epsilon)$ are non-resonant. In that case, aside from the transformations $F_\bullet(\check s,\check m)$, $\bullet=O,I +,I -$, we have also unique local normalizing transformations $F_i(\check s,\check m)$ defined on a neighborhood $\Omega_i(\check\mu,\check\epsilon)$, $i=1,2$, of $s_i(\check\mu,\check\epsilon)$ not containing any other singularity $s_j(\check\mu,\check\epsilon)$ nor the origin, with $F_i(\check s_i(\check\mu,\check\epsilon),\check m)=I$. They satisfy
\[\Barbar F_1=\Bar F_1,\qquad \Barbar F_2=\Bar F_2^\p.\]
Let $A_i$ be the connection matrix between $F_i\Psi_{I +}$ and $F_{I +}\Psi_{I +}$:
\[F_{I +}\Psi_{I +}=F_i\Psi_{I +}A_i,\]
see Fig.~\ref{figure:RC-7}.
It is easy to see that the monodromy of $F_{I +}\Psi_{I +}$ around the point $s_1$ (resp.\ $s_2$) is equal to
\[C_1 N_1=A_1^{-1}N_1A_1,\qquad \big(\text{resp.} \ N_2C_2=A_2^{-1}N_2A_2\big),\]
from which one can calculate using Lemma \ref{lemma:RC-4} that
\begin{gather}\label{eq:RC-A1A2}
A_1=\begin{pmatrix} 1 & \frac{1}{e_1^2-1}c_1 \\ 0 & \kappa_I \end{pmatrix},\qquad
A_2=\begin{pmatrix} \kappa_I & 0 \\ \frac{1}{1-e_2^{-2}}c_2 & 1 \end{pmatrix},
\end{gather}
with
\[e_1(\check\mu,\check\epsilon):={\rm e}^{\frac{s_1(\check\mu,\check\epsilon)}{\sqrt{\check\epsilon}}\pi{\rm i}},\qquad
 e_2(\check\mu,\check\epsilon):={\rm e}^{\frac{s_2(\check\mu,\check\epsilon)}{\sqrt{\check\epsilon}}\pi{\rm i}},\]
and $c_1={\rm i}\kappa^{-1}\big(\gamma-\frac{e_1}{e_2}-\frac{e_2}{e_1}\big)$, and $c_2=-{\rm i}\kappa \frac{e_1}{e_2}$.

\begin{figure}[t]\centering
\begin{tikzpicture}[scale=1]
\node at (0,0){\includegraphics[width=0.81\textwidth]{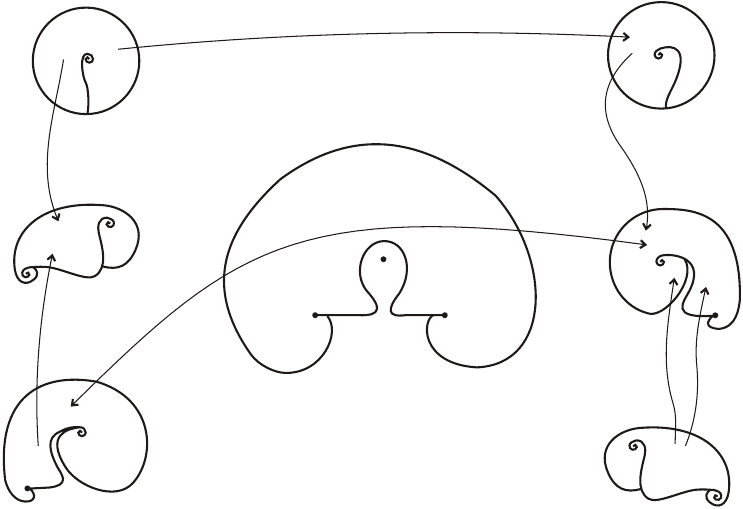}};
\node at (0,4.1){$\Bar N_1$};
\node at (5.3,3.9){$\Barbar F_2 \Barbar\Psi_{I+}$};
\node at (-5.0,3.8){$\Bar F_2^\p \Bar\Psi_{I+}^\p$};
\node at (-4.6,3.1){$-\Bar s_2$};
\node at (4.8,3.4){$\Barbar s_2$};
\node at (5.6,3.2){$\Barbar N_2$};
\node at (4.0,1.9){$\Barbar A_2$};
\node at (-5.9,1.9){$\Bar A_2^\p$};
\node at (-5.4,0.25){$\Bar F_{I+}^\p \Bar\Psi_{I+}^\p$};
\node at (5.6,0.25){$\Barbar F_{I+} \Barbar\Psi_{I+}$};
\node at (3.6,0.6){$\Barbar C_3$};
\node at (4.8,-0.2){$\Barbar s_2$};
\node at (6.2,-0.8){$\Barbar s_1$};
\node at (6.5,-4.2){$\Barbar s_1$};
\node at (0,1.1){$\Bar F_0 \Bar\Psi_0=\Barbar F_0\Barbar\Psi_0$};
\node at (4.7,-2.1){$\Barbar N \Barbar C_2$};
\node at (-3.5,-1.5){$\Bar C_4^\p$};
\node at (-6.1,-1.5){$\Bar N_1^\p$};
\node at (6.0,-2.1){$\Barbar N_1$};
\node at (-6.3,-0.7){$-\Bar s_1$};
\node at (-6.0,-3.8){$-\Bar s_1$};
\node at (-5.0,-2.9){$-\Bar s_2$};
\node at (-4.7,-3.5){$\Bar F_{I-}^\p \Bar\Psi_{I-}^\p$};
\node at (-1.4,-1.1){$-\Bar s_1$};
\node at (1.6,-1.1){$\Barbar s_1$};
\node at (5.5,-3.8){$\Barbar F_{I-} \Barbar\Psi_{I-}$};
\end{tikzpicture}
\caption{Connection matrices between fundamental solutions $F_\bullet\Psi_\bullet$ of Lemma~\ref{lemma:RC-4}, with $\check\epsilon$ fixed and $\Bar\mu\neq\Barbar\mu$.
(Picture with $\arg\check\epsilon=0$.) The corresponding diagram for the diagonal solutions $\Psi_\bullet$ of $\Bar\Delta^{s}$ is obtained by erasing all the $F$'s and replacing the matrices $A_i$, $C_i$ by identity matrices. The top arrow in the diagram here $\Bar F_2^\p\Bar\Psi_{I +}^\p \xrightarrow{\ \Bar N_1 \ }\Barbar F_2\Barbar\Psi_{I +}$ follows from the corresponding arrow $\Bar\Psi_{I +}^\p \xrightarrow{\ \Bar N_1 \ } \Barbar\Psi_{I +}$ which one can easily read in the corresponding diagram for the diagonal solutions.}\label{figure:RC-7}
\end{figure}

Knowing that $\Barbar F_2=\Bar F_2^\p$ one can see from Fig.~\ref{figure:RC-7} that
\begin{gather*}
 \Barbar A_2\Barbar C_3^{-1} =\Bar N_1^{-1} \Bar A_2^\p \big(\Bar N_1^\p\big)^{-1}\big(\Bar C_4^\p \big)^{-1},
\end{gather*}
where $\big(\Bar N_1^\p\big)^{-1}=\Bar N_1$, i.e.,
\begin{gather}\label{eq:RC-A2C3}
\begin{pmatrix} \Barbar\kappa_I & -{\rm i}\Barbar\kappa_I \frac{\Barbar e_2}{\Barbar e_1} \vspace{2mm}\\
{\rm i}\Barbar\kappa \frac{\Barbar e_1\Barbar e_2}{1-\Barbar e_2^2} & \Barbar\kappa \frac{1}{1-\Barbar e_2^2} \end{pmatrix}
=\begin{pmatrix} \Bar\kappa \frac{\Bar e_2}{\Bar e_2-\Bar e_2^{-1}} & -{\rm i}\Bar\kappa \frac{\Bar e_1^{-1}}{\Bar e_2-\Bar e_2^{-1}} \vspace{2mm}\\
{\rm i}\Bar\kappa_I \frac{\Bar e_1}{\Bar e_2} & \Bar\kappa_I\end{pmatrix}.
\end{gather}
This is satisfied if and only if
\begin{gather*}
\Barbar\kappa_I\Bar\kappa_I=\kappa_O\tfrac{\Bar e_2}{\Bar e_2-\Bar e_2^{-1}},
\end{gather*}
which is equivalent to \eqref{eq:RC-11}.
Similarly, one would find that
\[\Barbar A_1\Barbar C_3^{-1}=\Bar A_1\Bar C_3^{-1},\]
which is satisfied without imposing any new condition, since
\[A_1C_3^{-1}=\begin{pmatrix} 1 & \frac{{\rm i}\gamma e_1^{-1}-{\rm i}e_2-{\rm i}e_2^{-1}}{e_1-e_1^{-1}} \\ 0 & \kappa_O \end{pmatrix}.\]

\begin{figure}[t]\centering
\begin{tikzpicture}[scale=1]
\node at (0,0){\includegraphics[width=0.98\textwidth]{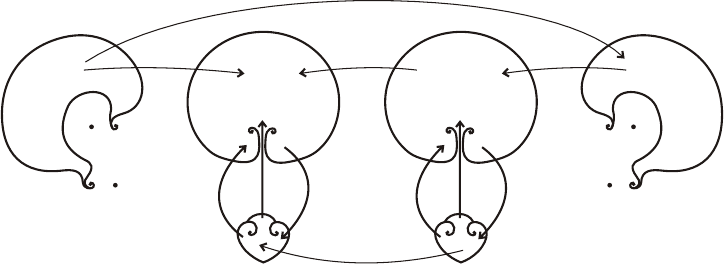}};
\node at (-4.3,2.6){$I$};
\node at (-4.3,1.7){$I$};
\node at (0,1.7){$I$};
\node at (4.3,1.7){$I$};
\node at (-6.4,1.5){$\tilde{\tilde\Psi}_{O}$};
\node at (-1.9,1.5){$\tilde{\tilde\Psi}_{I+}$};
\node at (2.3,1.5){$\tilde\Psi_{I-}^\p$};
\node at (6.3,1.5){$\tilde\Psi_{O}^\p$};
\node at (0,-2.4){$\tilde{\tilde N}{}^{-1}=\tilde N$};
\node at (-6.2,0.0){$\tilde{\tilde s}_2$};
\node at (-5.1,0.0){$\tilde{\tilde s}_1$};
\node at (-6.0,-1.5){$-\tilde{\tilde s}_1$};
\node at (-5.1,-1.5){$-\tilde{\tilde s}_2$};
\node at (-2.7,0.0){$\tilde{\tilde s}_2$};
\node at (-1.65,0.0){$\tilde{\tilde s}_1$};
\node at (1.5,0.0){$-\tilde s_1$};
\node at (2.8,0.0){$-\tilde s_2$};
\node at (5.0,-0.1){$-\tilde s_1$};
\node at (5.7,-0.1){$-\tilde s_2$};
\node at (5.0,-1.4){$\tilde s_2$};
\node at (5.85,-1.4){$\tilde s_1$};
\node at (-1.0,-1.3){$I$};
\node at (1.0,-1.3){$I$};
\node at (-1.85,-1.3){$\tilde{\tilde N}_1$};
\node at (1.85,-1.3){$\tilde N_1$};
\node at (-2.85,-1.3){$\tilde{\tilde N}$};
\node at (2.85,-1.3){$\tilde N$};
\node at (-2.9,-2.6){$\tilde{\tilde\Psi}_{I-}$};
\node at (3.1,-2.6){$\tilde\Psi_{I+}^\p$};
\end{tikzpicture}

\caption{Connection matrices between fundamental solutions $\Psi_\bullet$ of Lemma~\ref{lemma:RC-4} with $(\Tilde{\Tilde\mu},\Tilde{\Tilde\epsilon})\allowbreak ={\rm e}^{2\pi{\rm i}} (\Tilde\mu,\Tilde\epsilon)$.}\label{figure:RC-8}
\end{figure}

(b) Similarly to (a), the passage between $(\Tilde\mu,\Tilde\epsilon)$, $\Tilde\mu\in\check\cM(\Tilde\epsilon),$ and
$(\tilde{\tilde\mu},\tilde{\tilde\epsilon})={\rm e}^{2\pi{\rm i}}(\Tilde\mu,\Tilde\epsilon)$, $\tilde{\tilde\mu}\in\check\cM(\tilde{\tilde\epsilon})$,
is that of crossing the curve
$\Sigma_O(\epsilon)$, which affects only the outer zones, and hence the outer domains.
The inner domains rotate together with their vertices $s_1(\check\mu,\check\epsilon),\ s_2(\check\mu,\check\epsilon)$, therefore
$\tilde{\tilde\Omega}_{I +}=\Tilde\Omega_{I -}^\p$ and $\tilde{\tilde\Omega}_{I -}=\Tilde\Omega_{I +}^\p$.
So we have
\[\tilde{\tilde F}_{I +}=\Tilde F_{I -}^\p,\qquad \tilde{\tilde F}_{I -}=\Tilde F_{I +}^\p.\]
One can see from Fig.~\ref{figure:RC-8} that the fundamental solutions $\Psi_{I \pm}$ of the diagonal system $\widebar\Delta^{ s}$ satisfy
\[\tilde{\tilde\Psi}_{I +}=\Tilde\Psi_{I -}^\p,\qquad \tilde{\tilde\Psi}_{I -}=\Tilde\Psi_{I +}^\p\Tilde N.\]
This then implies \eqref{eq:RC-C12}, i.e.,
\begin{gather*}
 \Tilde C_1=\left(\begin{smallmatrix} 1 & -{\rm i}\Tilde{\Tilde\kappa} {\rm e}^{ -2\Tilde a\pi{\rm i}} \\ 0 & 1 \end{smallmatrix}\right),\qquad \!
\tilde{\tilde C}_1=\left(\begin{smallmatrix} 1 & -{\rm i}\Tilde\kappa {\rm e}^{2\Tilde{\Tilde a} \pi{\rm i}} \\ 0 & 1 \end{smallmatrix}\right), \qquad\!
 \Tilde C_2=\left(\begin{smallmatrix} 1 & 0 \\ -{\rm i}\Tilde\kappa {\rm e}^{ 2\Tilde a\pi{\rm i}} & 1 \end{smallmatrix}\right), \qquad \!
\tilde{\tilde C}_2=\left(\begin{smallmatrix} 1 & 0 \\ -{\rm i}\Tilde{\Tilde\kappa} {\rm e}^{ 2\Tilde{\Tilde a}\pi{\rm i}} & 1 \end{smallmatrix}\right),
\end{gather*}
as $\tilde{\tilde a}=-\Tilde a$. Then \eqref{eq:RC-12} follows from \eqref{eq:RC-gammac1c2}.
\end{proof}

The following proposition gives a semi-explicit formula for the determinants $\kappa$, $\kappa_I$, $\kappa_O$, analogical to the Gauss--Kummer formula for the hypergeometric equation \cite{IKSY, Koh}, and similar to the connection formulas of \cite{BJL}.

\begin{Proposition}\label{proposition:RC-5} \quad
\begin{enumerate}\itemsep=0pt
\item[$(a)$] Let $\Delta$ be a parametric system, $\Delta^{ s}$ its transform \eqref{eq:RC-d7a}, let $F_\bullet$ be the normalizing gauge transformations from Proposition~{\rm \ref{proposition:RC-3}} determined by the condition~\eqref{eq:RC-kappaIO} and let $\Psi_\bullet$ be as in Fig.~{\rm \ref{figure:RC-5}}.
The collection of all the connection matrices between the fundamental solutions $F_\bullet\Psi_\bullet$ is uniquely determined by $\kappa=\frac{\kappa_O}{\kappa_I}$ and by the invariant~$\gamma$, satisfying the relation~\eqref{eq:RC-12}.

\item[$(b)$] Let $\gamma(m)$ be a germ of analytic function and assume that there exists an analytic germ $Q(m)$ such that
\[\gamma=2\cos 2\pi Q.\]
Let
\[a(\check\mu,\check\epsilon):=\frac{s_1-s_2}{2\sqrt{\check\epsilon}},\qquad
b(\check\mu,\check\epsilon):=\frac{s_1+s_2}{2\sqrt{\check\epsilon}},\]
with $s_1$, $s_2$ as in \eqref{eq:RC-si}.
Then any triple of functions $\kappa_I$, $\kappa_O$, $\kappa=\frac{\kappa_O}{\kappa_I}\in\Cal B(\rM)$
with $\kappa_O(\check m)=1$ if $(\check\mu,\check\epsilon)=0$, $\check\mu\in\check{\Cal M}(0)$, satisfying the relations~\eqref{eq:RC-11},~\eqref{eq:RC-11b} and~\eqref{eq:RC-12} of Lemma~{\rm \ref{lemma:RC-5}} are equal to
\begin{gather}\label{eq:RC-kappaxI}
\kappa_I =\begin{cases}
	\frac{\sqrt{\frac{s_1s_2}{\check\epsilon}} \Gamma\big(\frac{s_1}{\sqrt{\check\epsilon}}\big) \Gamma\big(\frac{s_2}{\sqrt{\check\epsilon}}\big) }
	{\Gamma(1 + b - Q) \Gamma(b + Q) } {\rm e}^{ 2b\log b - \frac{s_1}{\sqrt{\check\epsilon}}\log\frac{s_1}{\sqrt{\check\epsilon}} - \frac{s_2}{\sqrt{\check\epsilon}}\log\frac{s_2}{\sqrt{\check\epsilon}} + f_I},&\text{if }\epsilon\neq 0,\\
	1,&\text{if }\epsilon=0,\end{cases}	\\
\label{eq:RC-kappaxO}
\kappa_O = \begin{cases}
	2\pi \frac{\frac{s_1}{\sqrt{\check\epsilon}}\Gamma\big(\frac{s_1}{\sqrt{\check\epsilon}}\big)^2}
	{\Gamma(1 + a - Q) \Gamma(1 + b - Q) \Gamma(a + Q) \Gamma(b + Q) }\\
\quad {}\times
	{\rm e}^{ 2a\log a + 2b\log b - 2\frac{s_1}{\sqrt{\check\epsilon}}\log\frac{s_1}{\sqrt{\check\epsilon}} + f_O}, &
	\text{if }\epsilon\neq 0,\\
	2\pi \frac{1}{\Gamma(1 + a - Q) \Gamma(a + Q) } 	{\rm e}^{ 2a\log a -2a + f},\ \text{where}\ a=\frac{1}{2\sqrt{\mu}}, &\text{if }\epsilon=0,\ \mu\neq 0,\\
	1, &\text{if }(\mu,\epsilon)=0,\end{cases}	\\
\label{eq:RC-kappax}
\kappa = \begin{cases}
	 2\pi \frac{\sqrt{\frac{s_1}{s_2}} \Gamma\big(\frac{s_1}{\sqrt{\check\epsilon}}\big) \Gamma\big(\frac{s_2}{\sqrt{\check\epsilon}}\big)^{-1}}
	{\Gamma(1 + a - Q) \Gamma(a + Q) }
	{\rm e}^{ 2a\log a - \frac{s_1}{\sqrt{\check\epsilon}}\log\frac{s_1}{\sqrt{\check\epsilon}} + \frac{s_2}{\sqrt{\check\epsilon}}\log\frac{s_2}{\sqrt{\check\epsilon}} + f},
	&\text{if }\epsilon\neq 0,\\
	2\pi \frac{1}{\Gamma(1 + a - Q) \Gamma(a + Q)} 	{\rm e}^{ 2a\log a -2a + f},\ \text{where}\ a=\frac{1}{2\sqrt{\mu}}, &\text{if }\epsilon=0,\ \mu\neq0,\\
	1, &\text{if }(\mu,\epsilon)=0,\end{cases}
\end{gather}
where $\Gamma$ is the gamma function and
\begin{gather*}
f=(s_1 + s_2) g\big(s_1s_2, s_1^2 + s_2^2,m\big),\\
 f_I=(s_1 - s_2) g\big( {-}s_1s_2, s_1^2 + s_2^2,m\big),\qquad f_O=f+f_I,\end{gather*}
for a unique analytic germ~$g$.
\end{enumerate}
\end{Proposition}

\begin{proof}{(a)} All the connection matrices between the fundamental solutions $F_\bullet\Psi_\bullet$ can be determined from Lemmas \ref{lemma:RC-4} and \ref{lemma:RC-5}.

{(b)} Denote $\sigma\colon \Bar m \mapsto \Barbar m$ the continuation map from Lemma~\ref{lemma:RC-5}(a),
and $\rho\colon \Tilde m \mapsto \tilde{\tilde m}$ the continuation map from Lemma~\ref{lemma:RC-5}(b).
Hence,
\begin{gather*} \sqrt\epsilon \circ \sigma=\sqrt\epsilon,\qquad s_1 \circ \sigma=s_1,\qquad s_2 \circ \sigma={\rm e}^{\pi{\rm i}}s_2, \qquad a \circ \sigma=b, \qquad b \circ \sigma=a,\\
 \tfrac{s_1}{\sqrt\epsilon} \circ \rho=\tfrac{s_2}{\sqrt\epsilon},\qquad \tfrac{s_2}{\sqrt\epsilon} \circ \rho=\tfrac{s_1}{\sqrt\epsilon}, \qquad a \circ \rho={\rm e}^{-\pi{\rm i}}a,\qquad b \circ \rho=b.
\end{gather*}
One can easily verify that the functions $\kappa_I$, $\kappa_O$, $\kappa$ of \eqref{eq:RC-kappaxI}--\eqref{eq:RC-kappax} satisfy $\kappa=\frac{\kappa_O}{\kappa_I}$ and the identities
\begin{gather*}
\eqref{eq:RC-11}\colon \quad \kappa_O \circ \sigma=\kappa_O=\frac{\kappa (\kappa \circ \sigma) {\rm e}^{\frac{s_2}{\sqrt{\check\epsilon}}\pi{\rm i}}}{2i \sin \frac{s_2}{\sqrt{\check\epsilon}}\pi},\\
\eqref{eq:RC-11b}\colon \quad \kappa_I \circ \rho=\kappa_I,\\
\eqref{eq:RC-12}\colon \quad 2\cos 2\pi Q=2\cos 2\pi a-\kappa (\kappa \circ \rho) {\rm e}^{- 2 a\pi{\rm i}},
\end{gather*}
using the standard reflection formula $\Gamma(z)\Gamma(1-z)=\frac{\pi}{\sin\pi z}$.
It follows from the Stirling formula:
\begin{gather*}
\Gamma(1+z)\sim \sqrt{2\pi z} \big(\tfrac{z}{e}\big)^z\big( 1+O\big(\tfrac{1}{z}\big)\big) \ \ \text{in the sector at infinity},
\end{gather*}
where $|\arg z|<\pi-\eta$ for any $0<\eta<\pi$, and Lemma~\ref{lemma:RC-arg} that
\[ \lim_{\substack{\check\epsilon\to 0\\ \check m\in\rM}}\kappa_I(\check m)=1
\qquad\text{and}\qquad
\lim_{\substack{(\check\mu,\check\epsilon)\to 0\\ \check m\in\rM}}\kappa_O(\check m)=1.\]

On the other hand if $\kappa_I$, $\kappa_O$, $\kappa$ are some functions satisfying the assumptions of the proposition,
let $\kappa'_I$, $\kappa'_O$, $\kappa'$ be given by \eqref{eq:RC-kappaxI}--\eqref{eq:RC-kappax} with $f_I=f_O=f=0$,
then it follows that the functions
\[f_I:=\log\frac{\kappa_I}{\kappa'_I},\qquad f_O:=\log\frac{\kappa_O}{\kappa'_O},\qquad f:=\log\frac{\kappa}{\kappa'},\]
satisfy
\begin{gather*} \begin{split} & f_O=f+f_I,\qquad
f_I \circ \rho=f_I,\qquad
f \circ \rho=-f,\qquad
f_O \circ \sigma=f_O,\\
& f \circ \sigma=f_I,\qquad
f_I \circ \sigma=f.\end{split}\end{gather*}
This implies, in particular, that $f_O\circ\sigma^2=f_O=f_O\circ\rho^2$, hence that $f_O$ is non-ramified as a~function of $(s_1,s_2)$, and therefore $f_O$ is an analytic function of $(s_1,s_2)$. Since one can express
\[f=\frac{1}{2}(f-f \circ \rho)=\frac{1}{2}(f_O-f_O \circ \rho) \qquad\text{and}\qquad f_I=\frac{1}{2}(f_O+f_O \circ \rho),\]
they too are germs of \emph{analytic} functions of $(s_1,s_2)$.
Moreover, for $\mu\neq 0$, $0=\lim\limits_{\check\epsilon\to 0} f_I=\lim\limits_{s_1-s_2\to 0} f_I$, so we can write
\[f_I=(s_1-s_2)\cdot g,\qquad\text{and}\qquad f=(s_1+s_2)\cdot (g\circ\sigma),\qquad f_O=f_I+f,\]
with $g$ that is $\rho$-invariant, thus a germ of an analytic function of
\[s_1s_2=\sqrt{\check\mu^2-\check\epsilon} \qquad\text{and}\qquad s_1^2+s_2^2=2\check\mu,\]
which are algebraically independent and form a Hilbert basis of the space of polynomials of $(s_1,s_2)$ that are invariant to the action of~$\rho$.
\end{proof}

\begin{Corollary}\label{corollary:RC-kappaI}The determinants $\kappa_\bullet(\check m)\neq 0$ for $m$ small, and $\lim\limits_{\substack{m\to 0\\ \check m\in\rM}}\kappa_\bullet(\check m)=1$, $\bullet=O,I$.
\end{Corollary}

\begin{proof}From the formulas \eqref{eq:RC-kappaxI} and \eqref{eq:RC-kappaxO} of Proposition~\ref{proposition:RC-5} and Lemma~\ref{lemma:RC-arg}.
\end{proof}

\begin{Remark}The set of points $\check m$ for which
\[\frac{s_1}{\sqrt{\check\epsilon}} \in-\N^*,\qquad \text{or} \qquad
\frac{s_2}{\sqrt{\check\epsilon}} \in-\N^*, \qquad \text{or} \qquad
a \in \pm Q-\N^*, \qquad \text{or}\qquad
b \in \pm Q-\N^*, \]
outlines a natural boundary for the set $\rM$.
\end{Remark}

\subsubsection{Proofs of the main theorems}

Let $ \check x(\check s,m)=\check s^2-\check\mu(m)$, a one-to-one map from the ramified coordinate $\check s$ to a ramified coordinate~$\check x$, be a lifting of the map
$x(s,m)=s^2-\mu(m)$~\eqref{eq:RC-xs}. Then by Lemma~\ref{lemma:RC-3} the ramified images of $\Omega_O$, $\Omega_{I\pm}$ in the $\check x$-coordinate cover for each $(\check\mu,\check\epsilon)$ a full neighborhood of each singular point $x=\pm\sqrt{\epsilon}\neq-\mu$ and of $x=0$ if $(\mu,\epsilon)=0$. Define $\rXO(\check\mu,\check\epsilon)$, $\rXI{}_\pm(\check\mu,\check\epsilon)$, depending continuously on~$\check m\in\sM$, as simply connected ramified extensions of these images, in such a way that they agree with them near these singularities, are open away of the singularities, and the union of their projections covers either all~$\sX$, or $\sX\setminus\{-\mu\}$ if $\epsilon=\mu^2\neq 0$.

\begin{proof}[Proof of Theorem~\ref{theorem:RC-IV}] The fundamental solution of a system $\Delta(x,m)$ corresponding to the fundamental solution $\Phi_\bullet=F_\bullet\Psi_\bullet$ of the associated system in the $s$-coordinate is given by
\[Y_\bullet(\check x,\check m)=S(s)VF_\bullet(\check s,\check m)\Psi_\bullet(\check s,m)=H_\bullet(\check x,\check m) S(s)V\Psi_\bullet(\check s,m),\]
where $H_\bullet(\check x,\check m)=S(s)VF_\bullet(\check s,\check m)V^{-1}S(s)^{-1}$. The function $\Theta_\bullet(\check x,m)=\theta_\bullet(\check s,m)$~\eqref{eq:RC-Theta} is a~branch of~\eqref{eq:RC-theta}, chosen according to the Fig.~\ref{figure:RC-8}.
	
The connection matrices between the fundamental solutions $Y_\bullet$ on the different domains $\rX_\bullet$ in Fig.~\ref{figure:RC-YOI}, are obtained from Fig.~\ref{figure:RC-6} and from the formulas~\eqref{eq:RC-monodromy1} and~\eqref{eq:RC-monodromy2}.
\end{proof}

\begin{proof}[Proof of Proposition~\ref{proposition:RC-IVa}] This is a consequence of Lemma~\ref{lemma:RC-boundedsolutions}.
\end{proof}

\begin{proof}[Proof of Proposition~\ref{proposition:RC-IVb}] This is how the fundamental solution matrices $Y_\bullet$, $\bullet=O,I\pm$, are constructed. We have $\kappa_O(\check m)=\det Y_O(\check x,\check m)$ and $\kappa_I(\check m)=\det Y_{I\pm}(\check x,\check m)$. In the case of $\kappa_O$, we know that $\kappa_O(0)=1$ and $\kappa_O(\check m)$ is non-vanishing because of its continuity in $\check m\in\rM$. This argument no longer works for~$\kappa_I$, but we know that it is non-vanishing by Corollary~\ref{corollary:RC-kappaI}.
\end{proof}	

\begin{proof}[Proof of Corollary~\ref{theorem:RC-III}] 	Since the fundamental solutions $Y_\bullet=SVF_\bullet\Psi_\bullet$ of $\Delta(x,m)$ and $SVE_\bullet\Psi_\bullet$ of $\widehat\Delta(x,m)$ (see Section~\ref{sec:RC-2.3}) are analytic away from the singularities $x=\pm\sqrt\epsilon$, the gauge transformations \eqref{eq:RC-T}
	\[T_\bullet(\check x,\check m)=S(s)VF_\bullet(\check s,\check m) E_\bullet(\check s,\check m)^{-1}V^{-1}S(s)^{-1},\qquad \bullet=O,I\]
	extend then on $\rX_\bullet$, $\bullet=O,I\pm$, as normalizing transformation for the parametric system $\Delta(x,m)$:
	\[T_\bullet\in\GL_2\big({\mathcal B}\big(\rX_\bullet\big)\big),\qquad T_\bullet^*\Delta=\widehat\Delta.
	\tag*{\qed}	\]
	\renewcommand{\qed}{}
\end{proof}

\begin{proof}[Proof of Theorem~\ref{theorem:RC-I}(a)] Let $\Delta(x,m)$, $\Delta'(x,m)$ be two parametric families of systems, let $\Delta^{ s}(s,m)$, $\Delta'^{s}(s,m)$ be their transforms \eqref{eq:RC-d7a}, and let $F_\bullet$, $F'_\bullet$ be the normalizing transformations from Proposition \eqref{proposition:RC-3} determined by the condition \eqref{eq:RC-kappaIO} with $\kappa_\bullet$, $\kappa'_\bullet$. Suppose that their invariants $\gamma=\gamma'$ are the same. We want to show that the two families of systems $\Delta, \Delta'$ are then analytically equivalent. We know that $ \kappa_I(\check m)=1=\kappa'_I(\check m) $ {when} $\check\epsilon=0$, and $ \kappa_O(\check m)=1=\kappa'_O(\check m) $ {when} $(\check\mu,\check\epsilon)=0$. Let $\delta(\check m)$, depending continuously on the parameter $\check m\in\rM$, be such that
\[\frac{\kappa'_O}{\kappa'_I}=\delta^2 \frac{\kappa_O}{\kappa_I},\qquad \delta(0)=1.\]
The relation \eqref{eq:RC-12} implies that $ \delta(\tilde{\tilde m})\cdot\delta(\Tilde m)=1.$
Put
\[F''_O=\delta^{-1} F'_O,\qquad F''_I=F'_I\left(\begin{smallmatrix} \delta^{-1} & 0 \\ 0\ & \delta \end{smallmatrix}\right).\]
They are also normalizing transformations for the system $\Delta'^{s}$: $(F''_\bullet)^*\Delta'^{s}=\Delta'^{s}$.
It is easily verified that the connection matrices between the fundamental solutions $F''_\bullet\Psi_\bullet$ are exactly the same as those between the fundamental solutions $F_\bullet\Psi_\bullet $ (with $\Psi_\bullet$ as in Fig.~\ref{figure:RC-6}), and one concludes by Proposition \ref{proposition:RC-4}.
\end{proof}

\begin{proof}[Proof of Theorem~\ref{theorem:RC-V}] (i) For $\epsilon(m)=0$, the transformation $T_{I+, m}$ converges to $T_{O,m}$, i.e., $|T_{I+, m}(s)-T_{O, m}(s)|\to 0$, $s\in S_{I+, m}$, if and only if $F_{I}(\cdot, m)$ converges to $F_{O}(\cdot,m)$, which happens if and only if the matrix $C_3(m)\to I$.

(ii) To show that the transformation $T_{2,m}$ converges to $T_{O, m}$, we need to show that the corresponding transformation $F_{2}(\cdot, m)$ converges to $F_{O}(\cdot,m)$.
It will be enough to show that the difference of fundamental solutions $F_2\Psi_O-F_O\Psi_O$ converges to~$0$ for each fixed~$s$. We know from the proof of Lemma \ref{lemma:RC-5}(a), Fig.~\ref{figure:RC-7}, that $F_O\Psi_O=F_2\Psi_OA_2C_3^{-1}$, where $A_2$ is given by~\eqref{eq:RC-A1A2} and~$A_2C_3^{-1}$ has
been calculated in~\eqref{eq:RC-A2C3}
\[A_2C_3^{-1}=\begin{pmatrix} \kappa_I & -{\rm i}\kappa_I \frac{e_2}{e_1} \vspace{1mm}\\ {\rm i}\kappa \frac{e_1 e_2}{1-e_2^2} & \kappa \frac{1}{1-e_2^2} \end{pmatrix}, \qquad
\text{where} \quad e_j={\rm e}^{\frac{s_j\pi{\rm i}}{\sqrt\epsilon}},\quad j=1,2.\]
We need that $A_2C_3^{-1}\to I$, which happens if and only if $\frac{e_2}{e_1}\to 0$ and $e_1e_2\to 0$ as $\epsilon(m)\to 0$,
i.e., $\Im\big(\frac{s_2-s_1}{\sqrt\epsilon}\big)>0$ and $\Im\big(\frac{s_2+s_1}{\sqrt\epsilon}\big)>0$.
For $\mu=O(\epsilon)$, we have $s_1=\epsilon^\frac{1}{4}+O\big(\epsilon^\frac{3}{4}\big)$, $s_2=\pm {\rm i}\epsilon^\frac{1}{4}+O\big(\epsilon^\frac{3}{4}\big)$, hence
$\frac{s_2-s_1}{\sqrt\epsilon}=\frac{-1\mp {\rm i}}{s_2}+O\big(\epsilon^\frac{1}{4}\big)$, $\frac{s_2+s_1}{\sqrt\epsilon}=\frac{-1\pm {\rm i}}{s_2}+O\big(\epsilon^\frac{1}{4}\big)$. Therefore the condition of convergence is satisfied if $\arg s_2\in \big(\frac{\pi}{4},\frac{3\pi}{4}\big)$, i.e., if $\arg x_2\in\big(\frac{\pi}{2},\frac{3\pi}{2}\big)$.
\end{proof}

\subsection*{Acknowledgment}
This article has been originally written during my doctoral studies at Universit\'e de Montr\'eal under the direction of Christiane Rousseau~-- I'd like to thank her for her support and encouragement. I'd also like to thank Alexey Glutsyuk and the anonymous referees for their numerous suggestions that helped to improve this paper.

\pdfbookmark[1]{References}{ref}
\LastPageEnding


\begin{thebibliography}{99}
\footnotesize\itemsep=0pt

\bibitem{BaVa}
Babbitt D.G., Varadarajan V.S., Local moduli for meromorphic differential
 equations, \textit{Ast\'erisque} \textbf{169--170} (1989), 217~pages.

\bibitem{Balser}
Balser W., Formal power series and linear systems of meromorphic ordinary
 differential equations, \textit{Universitext}, \href{https://doi.org/10.1007/b97608}{Springer-Verlag}, New York, 2000.

\bibitem{BJL}
Balser W., Jurkat W.B., Lutz D.A., Birkhoff invariants and {S}tokes'
 multipliers for meromorphic linear differential equations, \href{https://doi.org/10.1016/0022-247X(79)90217-8}{\textit{J.~Math.
 Anal. Appl.}} \textbf{71} (1979), 48--94.

\bibitem{BJL12-I}
Balser W., Jurkat W.B., Lutz D.A., A general theory of invariants for
 meromorphic differential equations. {I}.~{F}ormal invariants,
 \textit{Funkcial. Ekvac.} \textbf{22} (1979), 197--221.

\bibitem{BJL12-II}
Balser W., Jurkat W.B., Lutz D.A., A general theory of invariants for
 meromorphic differential equations. {II}.~{P}roper invariants,
 \textit{Funkcial. Ekvac.} \textbf{22} (1979), 257--283.

\bibitem{Ben}
Benzinger H.E., Plane autonomous systems with rational vector fields,
 \href{https://doi.org/10.2307/2001769}{\textit{Trans. Amer. Math. Soc.}} \textbf{326} (1991), 465--483.

\bibitem{Bo92}
Bolibrukh A.A., On sufficient conditions for the positive solvability of the
 {R}iemann--{H}ilbert problem, \href{https://doi.org/10.1007/BF02102113}{\textit{Math. Notes}} \textbf{51} (1992),
 110--117.

\bibitem{Bo94}
Bolibrukh A.A., On analytic transformation to {B}irkhoff standard form,
 \textit{Math. Dokl.} \textbf{49} (1994), 150--153.

\bibitem{BD}
Branner B., Dias K., Classification of complex polynomial vector fields in one
 complex variable, \href{https://doi.org/10.1080/10236190903251746}{\textit{J.~Difference Equ. Appl.}} \textbf{16} (2010),
 463--517, \href{https://arxiv.org/abs/0905.2293}{arXiv:0905.2293}.

\bibitem{CL}
Coddington E.A., Levinson N., Theory of ordinary differential equations,
 McGraw-Hill Book Company, Inc., New York~-- Toronto~-- London, 1955.

\bibitem{DES}
Douady A., Estrada F., Sentenac P., Champs de vecteurs polyn\^omiaux sur
 ${\mathbb C}$, {u}npublished manuscript, 2005.

\bibitem{Duv}
Duval A., Biconfluence et groupe de {G}alois, \textit{J.~Fac. Sci. Univ. Tokyo
 Sect. IA Math.} \textbf{38} (1991), 211--223.

\bibitem{G0}
Glutsuk A.A., Stokes operators via limit monodromy of generic perturbation,
 \href{https://doi.org/10.1023/A:1021744801409}{\textit{J.~Dynam. Control Systems}} \textbf{5} (1999), 101--135.

\bibitem{Glu}
Glutsyuk A.A., Resonant confluence of singular points and {S}tokes phenomena,
 \href{https://doi.org/10.1023/B:JODS.0000024125.05337.9e}{\textit{J.~Dynam. Control Systems}} \textbf{10} (2004), 253--302.

\bibitem{HLR}
Hurtubise J., Lambert C., Rousseau C., Complete system of analytic invariants
 for unfolded differential linear systems with an irregular singularity of
 {P}oincar\'e rank~{$k$}, \href{https://doi.org/10.17323/1609-4514-2014-14-2-309-338}{\textit{Mosc. Math.~J.}} \textbf{14} (2014),
 309--338.

\bibitem{HR}
Hurtubise J., Rousseau C., Moduli space for generic unfolded differential
 linear systems, \href{https://doi.org/10.1016/j.aim.2016.11.037}{\textit{Adv. Math.}} \textbf{307} (2017), 1268--1323,
 \href{https://arxiv.org/abs/1508.06616}{arXiv:1508.06616}.

\bibitem{Il}
Ilyashenko Yu., Realization of irreducible monodromy by {F}uchsian systems and
 reduction to the {B}irkhoff standard form (by {A}ndrey {B}olibrukh), in
 Differential Equations and Quantum Groups, \textit{IRMA Lect. Math. Theor.
 Phys.}, Vol.~9, \href{https://doi.org/10.4171/020-1/1}{Eur. Math. Soc.}, Z\"{u}rich, 2007, 1--8.

\bibitem{IlYa}
Ilyashenko Yu., Yakovenko S., Lectures on analytic differential equations,
 \textit{Graduate Studies in Mathematics}, Vol.~86, Amer. Math. Soc.,
 Providence, RI, 2008.

\bibitem{IKSY}
Iwasaki K., Kimura H., Shimomura S., Yoshida M., From {G}auss to {P}ainlev\'e.
 A modern theory of special functions, \textit{Aspects of Mathematics}, Vol.~E16, \href{https://doi.org/10.1007/978-3-322-90163-7}{Friedr. Vieweg \& Sohn}, Braunschweig, 1991.

\bibitem{JLP1}
Jurkat W., Lutz D., Peyerimhoff A., Birkhoff invariants and effective
 calcualtions for meromorphic linear differential equations.~{I},
 \href{https://doi.org/10.1016/0022-247X(76)90122-0}{\textit{J.~Math. Anal. Appl.}} \textbf{53} (1976), 438--470.

\bibitem{JLP2}
Jurkat W.B., Lutz D.A., Peyerimhoff A., Birkhoff invariants and effective
 calculations for meromorphic linear differential equations.~{II},
 \textit{Houston~J. Math.} \textbf{2} (1976), 207--238.

\bibitem{KT}
Kawai T., Takei Y., Algebraic analysis of singular perturbation theory,
 \textit{Translations of Mathematical Monographs}, Vol.~227, Amer. Math. Soc.,
 Providence, RI, 2005.

\bibitem{Kl}
Klime\v{s} M., Wild monodromy of the fifth {P}ainlev\'e equation and its action
 on the wild character variety: approach of confluence, \href{https://arxiv.org/abs/1609.05185}{arXiv:1609.05185}.

\bibitem{KR}
Klime\v{s} M., Rousseau C., Generic {$2$}-parameter perturbations of parabolic
 singular points of vector fields in~{$\mathbb{C}$}, \href{https://doi.org/10.1090/ecgd/325}{\textit{Conform. Geom.
 Dyn.}} \textbf{22} (2018), 141--184, \href{https://arxiv.org/abs/1710.00883}{arXiv:1710.00883}.

\bibitem{Koh}
Kohno M., Global analysis in linear differential equations, \textit{Mathematics
 and its Applications}, Vol.~471, \href{https://doi.org/10.1007/978-94-011-4605-0}{Kluwer Academic Publishers}, Dordrecht, 1999.

\bibitem{Ko}
Kostov V.P., Normal forms of unfoldings of non-{F}uchsian systems,
 \textit{C.~R.~Acad. Sci. Paris S\'er.~I Math.} \textbf{318} (1994), 623--628.

\bibitem{LR1}
Lambert C., Rousseau C., The {S}tokes phenomenon in the confluence of the
 hypergeometric equation using {R}iccati equation, \href{https://doi.org/10.1016/j.jde.2008.02.012}{\textit{J.~Differential
 Equations}} \textbf{244} (2008), 2641--2664, \href{https://arxiv.org/abs/0706.1773}{arXiv:0706.1773}.

\bibitem{LR}
Lambert C., Rousseau C., Complete system of analytic invariants for unfolded
 differential linear systems with an irregular singularity of {P}oincar\'e
 rank~1, \href{https://doi.org/10.17323/1609-4514-2012-12-1-77-138}{\textit{Mosc. Math.~J.}} \textbf{12} (2012), 77--138,
 \href{https://arxiv.org/abs/1105.2269}{arXiv:1105.2269}.

\bibitem{Mucino}
Muci\~{n}o Raymundo J., Valero-Vald\'es C., Bifurcations of meromorphic vector
 fields on the {R}iemann sphere, \href{https://doi.org/10.1017/S0143385700009883}{\textit{Ergodic Theory Dynam. Systems}}
 \textbf{15} (1995), 1211--1222.

\bibitem{Mullin}
Mullin F.E., On the regular perturbation of the subdominant solution to second
 order linear ordinary differential equations with polynomial coefficients,
 \textit{Funkcial. Ekvac.} \textbf{11} (1968), 1--38.

\bibitem{Parise}
Parise L., Confluence de singularit\'es r\'egulieres d'\'equations
 diff\'erentielles en une singularit\'e irr\'eguliere. Mod\`ele de Garnier,
 Ph.D.~Thesis, {IRMA} Strasbourg, 2001, available at
 \url{http://www-irma.u-strasbg.fr/annexes/publications/pdf/01020.pdf}.

\bibitem{Ra}
Ramis J.-P., Confluence et r\'esurgence, \textit{J.~Fac. Sci. Univ. Tokyo
 Sect.~IA Math.} \textbf{36} (1989), 703--716.

\bibitem{Sabbah}
Sabbah C., Isomonodromic deformations and {F}robenius manifolds. An
 introduction, \textit{Universitext}, Springer-Verlag London, Ltd., London, EDP
 Sciences, Les Ulis, 2007.

\bibitem{Sch}
Sch\"{a}fke R., Confluence of several regular singular points into an irregular
 singular one, \href{https://doi.org/10.1023/A:1022888516938}{\textit{J.~Dynam. Control Systems}} \textbf{4} (1998), 401--424.

\bibitem{Sibuya}
Sibuya Y., Global theory of a second order linear ordinary differential
 equation with a polynomial coefficient, \textit{North-Holland Mathematics
 Studies}, Vol.~18, North-Holland Publishing Co., Amsterdam~-- Oxford, 1975.

\bibitem{Sib}
Sibuya Y., Linear differential equations in the complex domain: problems of
 analytic continuation, \textit{Translations of Mathematical Monographs},
 Vol.~82, Amer. Math. Soc., Providence, RI, 1990.

\bibitem{Tahar}
Tahar G., Counting saddle connections in flat surfaces with poles of higher
 order, \href{https://doi.org/10.1007/s10711-017-0313-2}{\textit{Geom. Dedicata}} \textbf{196} (2018), 145--186,
 \href{https://arxiv.org/abs/1606.03705}{arXiv:1606.03705}.

\bibitem{Tomasini}
Tomasini J., Topological enumeration of complex polynomial vector fields,
 \href{https://doi.org/10.1017/etds.2013.100}{\textit{Ergodic Theory Dynam. Systems}} \textbf{35} (2015), 1315--1344,
 \href{https://arxiv.org/abs/1307.3850}{arXiv:1307.3850}.

\bibitem{Wa}
Wasow W., Asymptotic expansions for ordinary differential equations, \textit{Pure and
 Applied Mathematics}, Vol.~14, John Wiley \& Sons, Inc., New York~--
 London~-- Sydney, 1965.

\bibitem{Zha}
Zhang C., Confluence et ph\'enom\`ene de {S}tokes, \textit{J.~Math. Sci. Univ.
 Tokyo} \textbf{3} (1996), 91--107.

\end{thebibliography}
\end{document}